\documentclass[12pt]{article}
\usepackage{amstext,amsthm,amsmath,amsfonts,amssymb,amscd,latexsym,txfonts,stmaryrd,enumerate}
\usepackage{latexsym,stmaryrd}
\usepackage{color}
\usepackage{epsfig}
\usepackage{cite}
\usepackage[all,cmtip]{xy}
\usepackage{appendix}
\usepackage[colorlinks,citecolor={red}]{hyperref}
\usepackage[letterpaper, total={6in, 9in}]{geometry}
\usepackage{mathtools}
\usepackage{graphicx}
 \xymatrixcolsep{8mm} 

\newcommand{\Supp}{\operatorname{Supp}}
\newcommand{\Cb}{\mathcal{C}^\bullet}
\newcommand{\soc}{\operatorname{soc}}
\newcommand{\Hom}{\operatorname{Hom}}

\newcommand{\Ann}{\operatorname{Ann}}

\newcommand{\Spec}{\operatorname{Spec}}

\newcommand{\Sym}{\operatorname{Sym}}
\newcommand{\C}{{\mathbb C}}
\renewcommand{\P}{{\mathbb P}}

\newcommand{\Z}{{\mathbb Z}}

\newcommand{\N}{{\mathbb N}}
\newcommand{\Q}{{\mathbb Q}}

\newcommand{\F}{{\mathbb F}}

\newcommand{\m}{{\mathfrak{m}}}
\renewcommand{\L}{{\mathcal{L}}}

\def\ns{\footnotesize \it}
\def\cha{\mathrm{char}\,}

\theoremstyle{theorem}
\newtheorem{theorem}{Theorem}[section]
\newtheorem{theoremA}{Theorem}
\newtheorem*{theorem1*}{Theorem 1}
\newtheorem*{theorem2*}{Theorem 2}
\newtheorem{proposition}[theorem]{Proposition}
\newtheorem{lemma}[theorem]{Lemma}
\newtheorem{corollary}[theorem]{Corollary}
\theoremstyle{definition}

\newtheorem{example}[theorem]{Example}
\newtheorem{definition}[theorem]{Definition}
\newtheorem{remark}[theorem]{Remark}
\newtheorem{remarks}[theorem]{Remarks}

\newtheorem{fact}[theorem]{Fact}

\newtheorem*{problem*}{Problem}
\newtheorem*{question*}{Question}
\newtheorem*{remarks*}{Remarks}
\newtheorem*{claim*}{Claim}
\newtheorem*{remark*}{Remark}
\newtheorem*{conjecture*}{Conjecture}
\newtheorem*{fact*}{Fact}
\newtheorem*{ack}{Acknowledgment}

\newcommand\blfootnote[1]{%
  \begingroup
  \renewcommand\thefootnote{}\footnote{#1}%
  \addtocounter{footnote}{-1}%
  \endgroup
}

\title{Connected Sums of Graded Artinian Gorenstein Algebras and Lefschetz Properties}

\author{Anthony Iarrobino\\[.05in]
{\ns Department of Mathematics, Northeastern University, Boston, MA 02115,
 USA.} \\ {\ns a.iarrobino@northeastern.edu}
 \\[.2in]
Chris McDaniel\\[0.05in]
{\ns Endicott College, 376 Hale St
Beverly, MA 01915, USA.} \\ {\ns cmcdanie@endicott.edu}
\\[.2in]
Alexandra Seceleanu\\[0.05in]
{\ns Department of Mathematics, University of Nebraska-Lincoln,}\\
{\ns 203 Avery Hall,  Lincoln NE 68588, USA.} \\ {\ns aseceleanu@unl.edu}
}

\date{}

\begin{document}
		\maketitle
		
\begin{abstract}
In their paper \cite{AAM}, H. Ananthnarayan, L. Avramov, and W.F. Moore introduced a connected sum construction for local Gorenstein rings $A,B$ over a local Gorenstein ring $T$, which, in the graded Artinian case, can be viewed as an algebraic analogue of the topological construction of the same name.  We give two alternative descriptions of this algebraic connected sum:  the first uses algebraic analogues of Thom classes of vector bundles and Gysin homomorphisms, the second is in terms of Macaulay dual generators.  We also investigate the extent to which the connected sum of $A,B$ over an Artinian Gorenstein algebra $T$ preserves the weak or strong Lefschetz property, thus providing new classes of rings which satisfy these properties.
\end{abstract}

\blfootnote{\noindent \textbf{Keywords}: Artinian algebra, cohomology, connected sum, Gorenstein, Hilbert function, Lefschetz property, Gysin map}
  \blfootnote{\noindent\textbf{2010 Mathematics Subject Classification}: Primary: 13H10;  Secondary: 13E10,14C99, 14D06. }
  
\section{Introduction.}
Let $A$ and $B$ be two graded Artinian Gorenstein (AG) algebras of the same socle degree $d$, let $T$ be an AG algebra of socle degree $k<d$, and suppose there are surjective maps $\pi_A\colon A\rightarrow T$, and $\pi_B\colon B\rightarrow T$.  From this data, one forms the fibered product algebra $A\times_TB$ as the categorical pullback of $\pi_A,\pi_B$; the connected sum algebra $A\#_TB$ is the quotient of $A\times_TB$ by a certain principal ideal $\langle (\tau_A,\tau_B)\rangle\subset A\times_TB$.  The connected sum is again an AG algebra.

This algebraic connected sum operation for local Gorenstein algebras $A,B$ over a local Cohen-Macaulay algebra $T$ was introduced by H. Ananthnarayan, L. Avramov, and W.F. Moore (A-A-M) in their 2012 paper \cite{AAM}.  In the present paper, we focus on the graded Artinian case, and give a slightly different description of the A-A-M construction, taking our cues from topology.  

AG algebras can be viewed as algebraic analogues of cohomology rings (in even degrees) of smooth compact connected orientable manifolds, i.e. $A_i\cong H^{2i}(M,\F)$.  In this analogy one can view the connected sum of two AG algebras $A$ and $B$ over another AG algebra $T$ as the cohomology ring of a connected sum manifold, obtained   
by gluing two $2d$-dimensional manifolds $M_1$ and $M_2$ along diffeomorphic tubular neighborhoods of a common $2k$-dimensional submanifold $N$ (Theorem \ref{thm:TopCS}).  Here the elements $\tau_1$ and $\tau_2$ correspond to the respective Thom classes of the normal bundles of the submanifold $N$ in each of $M_1$ and $M_2$.  
Details on the connection between the algebraic and topological connected sums are provided in the Appendix.  For further references, see \cite{MS1} (also \cite{BT}) for discussion of Thom classes, and \cite{Kos} for a discussion of connected sums of manifolds.  
   
The connected sum of two AG algebras $A,B$ over the ground field $T=\F$ (corresponding to the connected sum of two manifolds over a point) has been studied by several authors \cite{ACLY}, \cite{ACLY2}, \cite{BBKT}, \cite{MS}, \cite{SS}, \cite{SS3}.  Moreover, connected sums of combinatorial objects such as simplicial complexes and polytopes have also led to algebraic connected sum operations on Stanley-Riesner rings \cite{BN},\cite{MM}.  While \cite{AAM} has already established several fundamental properties of the general connected sum construction (some of which we also prove here), we focus on two aspects which have not previously been considered in this generality. Specifically we study Macaulay-Matlis duality and Lefschetz properties of both fibered products and connected sums of AG algebras $A,B$ over other AG algebras $T$.  Our main results are summarized below.

Let $Q$ be a polynomial ring, let $R$ be the dual ring of $Q$ (so $R$ is a divided power algebra) and let $I\subset Q$ be an ideal such that the quotient $A=Q/I$ is an AG algebra of socle degree $d$.  A classical result of Macaulay states that there is a homogeneous polynomial $F\in R$ of degree $d$, unique up to scalar multiple, such that $I=\Ann(F)$.  The polynomial $F$ is called the Macaulay dual generator for $A$.  Since the connected sum of two AG algebras over an AG algebra is again an AG algebra, it seems natural to ask how their Macaulay dual generators are related.  The following result (Theorem \ref{theorem1}) characterizes such connected sums in terms of their dual generators:

\begin{theoremA}
	\label{thm:A}
	Let $F,G\in R_d$ be two linearly independent homogeneous forms of degree $d$, and suppose that there exists $\tau\in Q_{d-k}$ (for some $k<d$) satisfying
	\begin{enumerate}[(a)]
		\item $\tau\circ F=\tau\circ G\neq 0$, and 
		\item $\Ann(\tau\circ F=\tau\circ G)=\Ann(F)+\Ann(G)$.
	\end{enumerate} 
In this case, set
$$A=\frac{Q}{\Ann{F}}, \ B=\frac{Q}{\Ann(G)}, \ T=\frac{Q}{\Ann(\tau\circ F=\tau\circ G)},$$
and let $\pi_A\colon A\rightarrow T$ and $\pi_B\colon B\rightarrow T$ be the natural projection maps.  Then the Thom classes of $\pi_A$ and $\pi_B$ are given by $\tau_A=\tau+\Ann(F)$ and $\tau_B=\tau+\Ann(G)$, and we have algebra isomorphisms
$$A\times_TB\cong \frac{Q}{\Ann(F)\cap\Ann(G)}, \ \ A\#_TB\cong \frac{Q}{\Ann(F-G)}.$$
Conversely, every connected sum $A\#_T B$ of graded AG algebras of the same socle degree over graded AG algebra $T$ arises in this way.
\end{theoremA}   

For Thom classes see Definition \ref{def:ThomClass}, and for the connection of Thom classes with Macaulay dual generators, see Remark \ref{rmk:MDOrient} (b).\par382

Theorem \ref{thm:A} is new (as far as the authors can tell), although in the case $T=\F$, Theorem~\ref{thm:A} yields the well known result that $C=Q/\Ann(H)$ is a connected sum if and only if $H=F-G$ where $F$ and $G$ are polynomials in distinct sets of variables, possibly after a change of coordinates (Corollary \ref{cor:CSoverF}).  As another application of Theorem \ref{thm:A}, we characterize connected sums over arbitrary $T$ in the special case where $F$ and $G$ are monomials (Proposition \ref{prop:TFAE}).  In general, determining whether a homogeneous polynomial $H$ is the Macaulay dual of a connected sum over some $T$ seems to be a rather difficult problem.  Over $T=\F$, this problem has been studied by several authors, e.g. \cite{ACLY}, \cite{BBKT}, \cite{SS3}, and is related to the Waring rank of polynomials, see  \cite{BBKT}.

The strong Lefschetz property (SLP) for graded AG algebras is an algebraic version of a property of cohomology rings of smooth complex projective varieties stemming from the hard Lefschetz theorem in algebraic geometry.  We say a graded Artinian algebra $A=\bigoplus_{i=0}^dA_i$ satisfies SLP if there is a linear form $\ell\in A_1$ for which the multiplication maps $\times\ell^k\colon A_i\rightarrow A_{i+k}$ have full rank for each degree $i$ and each exponent $k$.  The weak Lefschetz property (WLP) requires only that the multiplication maps $\times\ell\colon A_i\rightarrow A_{i+1}$ have full rank for each $i$. These properties are especially intriguing for AG algebras, as they most closely resemble cohomology rings. 

We show that if $A$ and $B$ are two graded AG algebras satisfying SLP, then over $T=\F$, both their fibered product $A\times_\F B$ and their connected sum $A\#_\F B$ have SLP as well (Propositions \ref{prop:SLPFP} and \ref{prop:SLPCS}).  That connected sums over $\F$ preserve SLP has been previously shown by J. Watanabe and his coauthors \cite[Proposition 3.77]{HMMNWW} using higher Hessians.  While we do not use them here, we believe that higher Hessians may be a powerful tool for establishing SLP for connected sums over other $T$. On the other hand our examples show that, even if $A$ and $B$ both have SLP, the connected sum over general $T$ can fail to have SLP.  The corresponding fact from topology is that the connected sum of two complex projective manifolds over a complex submanifold may fail to be a projective manifold itself, or even homotopy equivalent to one (Remark~\ref{rem:CSNP}).  

However we prove that connected sums which generalize the cohomology rings of blowups of manifolds at submanifolds with trivial normal bundles retain the SLP (Theorem \ref{thm:blowup}). The description for these rings is given by Theorem \ref{thm:C}, which establishes the strong Lefschetz property for some families of Artinian Gorenstein algebras where it was previously unknown (Corollaries \ref{heightthreecor} and \ref{closurecor}).

\begin{theoremA}
\label{thm:C}
Let $T, A$ and $B=T[x]/(x^{d-k+1})$ be AG algebras with socle degrees $k, d, d,$ respectively and let $\pi_A\colon A\rightarrow T$ and $\pi_B\colon B\rightarrow T$ be surjective ring homomorphism such that the Thom class $\tau_A$ satisfies $\pi_A(\tau_A)=0$ and $\pi_B(t)=t$ for $t\in T$, while $\pi_B(x)=0$.  If $A$ and $T$ both satisfy the SLP,  then the fibered product $A\times_T B$ also satisfies the SLP and moreover, if the field $\F$ is algebraically closed, then the connected sum $A\#_T B$  also satisfies the SLP.
\end{theoremA}

 We also show that  connected sums and fibered products retain the WLP to some extent (Theorem~\ref{theorem2}). This gives another family of algebras for which the weak Lefschetz property was not previously known.

\begin{theoremA}
	\label{thm:B}
Let $A$ and $B$ be standard graded AG algebras of socle degree $d$ satisfying the SLP, and let $T$ be a graded AG algebra of socle degree $k$, with $k<\lfloor \frac{d-1}{2} \rfloor$,
endowed with $\F$-algebra homomorphisms $\pi_A:A\to T$ and $\pi_B:B\to T$. Then the resulting fibered product  $A\times_TB$ and the connected sum $A\#_T B$ both satisfy the WLP.
\end{theoremA}

Our Proposition \ref{ex:NonSLP} shows that we need quite restrictive hypotheses for the conclusion of Theorem \ref{thm:B}.  While the connected sum can yield new families that do have WLP or SLP (Theorem \ref{thm:blowup}), it can also lead to new families having other Jordan types.
\vskip 0.2cm

{\it This paper is organized as follows.}  In Section \ref{sect:Pre} we introduce algebraic versions of Thom classes and Gysin maps, named after their related topological objects, which are applied in Section \ref{sect:fpcs} to give an alternate description of the Ananthnarayan-Avramov-Moore construction.  In Section \ref{sect:MDG}, we review the basic tenets of Macaulay duality, prove Theorem \ref{thm:A}, and compute several examples.  In Section \ref{sect:LP} we study Lefschetz properties and find several new classes of rings which satisfy the SLP or the WLP:  the connected sums described in Theorem \ref{thm:C}, which generalize the cohomology rings of blowups of manifolds at a point; and the connected sums and fibered products over $\F$, the  rings described by Theorem~\ref{thm:B}. In the Appendix we describe the topological connected sum construction on smooth manifolds, and prove Theorem \ref{thm:TopCS}, which gives sufficient conditions for the cohomology ring of the (topological) connected sum of two manifolds to be the (algebraic) connected sum of the cohomology rings of the two manifolds.

\paragraph{Notation.} Throughout this paper we use the following notation.  We will assume unless otherwise stated, that all graded objects $M$ are graded over the non-negative integers $\N$, and $M_i$ denotes the $i^{th}$ graded component. For a graded object $M$ we will write $M(n)$ to be the graded object $M$ shifted up by $n$, meaning that $M(n)_d=M_{d+n}$.  All maps between graded objects will be graded of degree zero unless otherwise stated. 
All graded algebras $A$ are assumed to be commutative and connected over an arbitrary fixed field $\F$, meaning that $A_0=\F$. However, in Section \ref{sect:LP} about Lefschetz properties, we will assume that the field $\F$ is infinite of characteristic zero or of characteristic greater than the common socle degree $d$ of $A,B$. Our graded algebras are not necessarily standard graded, meaning that $A$ is not necessarily generated by $A_1$ as an algebra.  Given a graded algebra $A$, its homogeneous maximal ideal will be denoted by $\mathfrak{m}_A=\bigoplus_{i\geq 1}A_i$ or $A_+$. An algebra $A$ is called  \emph{Artinian} if it is a finite dimensional vector space over $\F$. The \emph{socle} of an Artinian algebra $A$ is the ideal $(0:A_+)$; its \emph{socle degree} is the largest integer $d$ such that $A_d\neq 0$: it is sometimes called the \emph{formal dimension} of $A$ \cite{MS,SS}. The \emph{type} of $A$ is the vector space dimension of its socle. The  \emph{Hilbert series} of an Artinian algebra $A$ is the generating function $H(A,t)=\sum_{i\in \N}\dim_\F(A_i)t^i$. In the topology literature this notion appears under the name of Poincar\'e series. By the \emph{Hilbert function} $H(A)$ of an Artinian algebra $A$, we mean the sequence of coefficients of its Hilbert series e.g. if $H(A,t)=1+2t+3t^2$, then $H(A)=(1,2,3)$.  We write $H(A)[n]$ to mean the coefficient sequence for the shifted Hilbert series $t^nH(A,t)$, e.g. if $H(A,t)=1+2t+3t^2$ then $H(A)[2]=(0,0,1,2,3)$. Given a set $S$ of elements in a vector space over the field $\F$ we denote by $\operatorname{span}_\F S$ their span.

{\hypersetup{linkcolor=black}
\tableofcontents}

\section{Preliminaries.}
\label{sect:Pre}
\subsection{Oriented AG Algebras.}
Let $A$ be a graded Artinian algebra with socle degree $d$.  We say that $A$ is \emph{Gorenstein} if its socle $(0:\mathfrak{m}_A)$ is one dimensional as an $\F$-vector space, which must then necessarily be $A_d$.  Equivalently, $A$ is Gorenstein if its $d^{th}$ graded piece $A_d$ is one dimensional, and for any non-zero map of graded vector spaces $\int_A\colon A\rightarrow \F(-d)$ (which can be obtained by fixing a vector space isomorphism $\int_A\colon A_d\rightarrow \F$ and extending it by zero to all of $A$) the bilinear pairing defined by multiplication in $A$ 
$$\xymatrixrowsep{.5pc}\xymatrix{A\times A\ar[r] & \F\\ 
	(a,a')\ar@{|->}[r] & \int_Aa\cdot a'}$$
is non-degenerate.  We call the pair $\left(A,\int_A\right)$ an \emph{oriented AG algebra}, where $\int_A$ is referred to as the {\em orientation}. Another notation found in the literature is $\langle \cdot,\cdot\rangle_\phi$ where $\phi$ is the isomorphism $\phi:A_d\to \F$. 

\subsubsection{Thom Classes}

Throughout this section suppose that $\left(A,\int_A\right)$ and $\left(T,\int_T\right)$ are oriented AG algebras with socle degrees $d$ and $k$, respectively, and suppose that $\pi\colon A\rightarrow T$ is a graded map between them.

\begin{lemma}
	\label{lem:TCExist}
	There exists a unique homogeneous element $\tau=\tau_{\pi}\in A_{d-k}$ such that 
	$$\int_A\tau\cdot a=\int_T\pi(a), \ \ \forall a\in A.$$
\end{lemma}
\begin{proof}
	Since the pairing $(-,-)\colon A\times A\rightarrow\F$ is non-degenerate, and $A$ is finite dimensional as an $\F$-vector space, the map $A\ni t\mapsto\int_A t\cdot(-)\in\operatorname{Hom}_\F(A,\F)$ is an isomorphism of graded vector spaces.  Therefore there exists a unique $\tau\in A$ corresponding to the homomorphism $\int_T\circ\pi\in\operatorname{Hom}_\F(A,\F)$, i.e. 
	$$\int_A\tau\cdot a=\int_T\pi(a), \ \ \forall a\in A.$$
	Since $\int_T\circ\pi\colon A\rightarrow \F$ is a graded map which vanishes on every graded component except $A_k$, $\tau$ must be homogeneous of degree $d-k$.
\end{proof}
\begin{definition}
	\label{def:ThomClass}
	The element $\tau\in A_{d-k}$ above is called the \emph{Thom class} for $\pi\colon A\rightarrow T$.  Note that it depends not only on the map $\pi$, but also on the orientations chosen for $A$ and $T$. 
\end{definition}

\begin{remark}
	\label{rem:MS}
	 D.M. Meyer and L. Smith refer to our Thom classes as transition elements \cite[p. 14]{MS}.  In fact they refer to a transition element as a Thom class in positive characteristic in certain cases when the two algebras carry an action of the Steenrod algebra \cite[p. 56]{MS}.  We will see in Section \ref{sect:MDG} (Remark \ref{rmk:MDOrient} (b)) that the Thom class has an alternative characterization in terms of the Macaulay dual generators of $A$ and $T$.
	
\end{remark}

\subsubsection{Gysin Maps.}
The map $\pi\colon A\rightarrow T$ gives $T$ an $A$-module structure.
\begin{lemma}[Gysin map]
	\label{lem:GysinMap}
	There exists an unique $A$-module map $\iota=\iota_{\pi}\colon T(k-d)\rightarrow A$ satisfying $\iota(1_T)=\tau_\pi$.
\end{lemma}
\begin{proof}
	For each $t\in T$, there is a unique corresponding homomorphism $\int_Tt\cdot(-)\colon T\rightarrow \F$ which pulls back via $\pi$ to give the homomorphism $\int_Tt\cdot \pi(-)\colon A\rightarrow \F$, which in turn corresponds to a unique element $\iota(t)\in A$.  In other words, the map $\iota\colon T(d-k)\rightarrow A$ is defined by the condition that 
	$$\int_A\iota(t)\cdot a=\int_Tt\cdot\pi(a), \ \ \forall a\in A.$$
	We need only check that it is a map of $A$-modules.  For $t_1,t_2\in T$ and fixed $a\in A$ we have for every $a'\in A$:
	\begin{align*}
	\int_A\iota(t_1+at_2)\cdot a'= & \int_T(t_1+\pi(a)t_2)\cdot\pi(a')\\
	= & \int_Tt_1\pi(a')+\int_T\pi(a)t_2\cdot\pi(a')\\
	= & \int_Tt_1\pi(a')+\int_Tt_2\pi(a\cdot a')\\
	= & \int_A\iota(t_1)\cdot a' + \int_Aa\cdot \iota(t_2)\cdot a'
	\end{align*}
	which shows that $\iota(t_1+at_2)=\iota(t_1)+a\cdot\iota(t_2)$, and hence that $\iota\colon T(k-d)\rightarrow A$ is an $A$-module map.  Clearly we must then have $\iota(1_T)=\tau_{\pi}$.
\end{proof}

\begin{definition}
	\label{def:GysinMap}
	The map $\iota\colon T(k-d)\rightarrow A$ we term the \emph{Gysin map} associated to $\pi\colon A\rightarrow T$.  Note that it also depends not only on $\pi$ but also on the chosen orientations of $A$ and $T$.\footnote{See Remark \ref{rem:PDGysin} for an explanation of this nomenclature.}
\end{definition}

\begin{lemma}
	\label{lem:InjSurj}
	The map $\pi\colon A\rightarrow T$ is surjective if and only if $\iota\colon T(k-d)\rightarrow A$ is injective.
\end{lemma}
\begin{proof}
	Assume that $\pi$ is surjective, and suppose that $\iota(t)=0$ for some $t\in T$.  Then $\int_A\iota(t)\cdot a=\int_Tt\cdot\pi(a)=0$ for all $a\in A$.  Since $\pi$ is surjective we must therefore have that $\int_Tt\cdot t'=0 $ for all $t'\in T$, and hence $t=0$ by the non-degeneracy of the pairing on $T$.  
	
	Conversely, assume that $\pi$ is not surjective, and let $S\subseteq T$ be the image of $\pi\colon A\rightarrow T$.  Since $S\neq T$, there is a non-zero $\phi\in\Hom_\F(T,\F)$ such that $\phi(S)=0$.  Let $t\in T$ be the non-zero element for which $\phi=\int_Tt\cdot(-)$.  Then we have $\int_A\iota(t)\cdot a=\int_Tt\cdot\pi(a)=0$ for all $a\in A$.  It follows from the non-degeneracy of the pairing on $A$ that $\iota(t)$ must be zero, hence $\iota\colon T(k-d)\rightarrow A$ is not injective.
\end{proof}

\begin{lemma}
	\label{lem:GysinMult}
	Assume that $\pi\colon A\rightarrow T$ is surjective.  Then the Gysin map $\iota\colon T(k-d)\rightarrow A$ coincides with the multiplication map $\times\tau\colon T(k-d)\rightarrow A$, $t\mapsto \tau\cdot t$ where $\tau\cdot t$ is interpreted as the product in $A$ given by $\tau\cdot a$, where $a$ is any $\pi$ lift of $t$. 
\end{lemma}
\begin{proof}
	First note that the multiplication map is well defined since for $a\in\ker(\pi)$ and any $a'\in A$ we have
	\begin{align*}
	\int_A\left(\tau\cdot a\right)\cdot a'= & \int_A\tau\cdot\left(a\cdot a'\right)\\
	= & \int_T\pi(a\cdot a')=0
	\end{align*}
	To see that the multiplication map coincides with the Gysin map note that for every $a'\in A$ we have 
\begin{align*}
\int_A\left(\tau\cdot t\right)\cdot a'= & \int_Tt\cdot\pi(a')\\
= & \int_A\iota(t)\cdot a'.
\end{align*}
Thus $\tau\cdot t=\iota(t)$ by non-degeneracy of the pairing on $A$.
\end{proof}

In the following remark we give an alternate interpretation of the Thom class and Gysin map using certain dualizing functors on AG algebras.

\begin{remark}
\label{rk:dualGysin}
Applying the dualizing functor $^\vee=\Hom_\F(-,\F)$ to the surjective map $\pi\colon A\rightarrow T$ yields an injective graded map $\pi^*:T^\vee \to A^\vee$. Since both $T$ and $A$ are AG algebras of socle degrees $k$ and $d$ respectively, there are graded isomorphisms $ T(k) \cong T^\vee $ and $ A(d)  \cong A^\vee$ given by $t\mapsto \int_T t\cdot -$ and $a\mapsto \int_A a\cdot -$. 
This can be summarized using the following commutative diagram, where $j$ is defined to be the composite of the inverse of the rightmost map and the other two maps.
\begin{equation*}
	\xymatrix{ T^\vee \ar[r]^-{\pi^*} & A^\vee\\  T(k) \ar[u]^-{\cong} \ar[r]_-{j} & A(d) \ar[u]^-{\cong}.\\}
\end{equation*}
The commutativity of the diagram yields  
$\int_A j(1_T)\cdot a=\int_T1_T\cdot \pi(a),  \forall a\in A$.

Comparing with Lemma \ref{lem:TCExist} gives that $\tau_\pi=j(1_T)$, that is, the Thom class is determined by the map $j$ induced by $\pi^*$.
The Gysin map introduced in Lemma \ref{lem:GysinMap} is thus a graded shift of $j$, namely $\iota_\pi=j(-d):T(k-d)\to A$.\par
Our Gysin map corresponds to the map $\iota_S:V\to S$ in \cite[\S 2, Fig. 2.0.1]{AAM} in the special case
that $V$ is the dualizing module of an AG algebra. 
Next we compute the image of the Gysin map for future reference. The diagram above can be extended to the following
\begin{equation*}
	\xymatrix{ 0\ar[r] &T^\vee \ar[r]^-{\pi^*} & A^\vee  \ar[r] & \ker(\pi)^\vee \ar[r] & 0\\  0\ar[r] &T(k) \ar[u]^-{\cong} \ar[r]_-{j} & A(d) \ar[u]^-{\cong} \ar[r] & A(d)/\rm{Im}(j) \ar[u]^-{\cong} \ar[r] & 0.\\}
\end{equation*}
By \cite[Proposition 3.6.16]{BH} the following functors are equal $-^\vee=\Hom_\F(-,\F)=\Hom_A(-, A^\vee)$. In light of this, we compute
$$\ker(\pi)^\vee=\Hom_A(\ker(\pi),A^\vee)\cong\Hom_A(\ker(\pi),A(d))=A(d)/(0:_{A(d)} \ker(\pi)),$$
which leads to ${\rm{Im}}(j)=(0:_{A(d)} \ker(\pi))$ and thus ${\rm{Im}}(\iota_\pi)=(0:_{A} \ker(\pi))$.

 \end{remark}

\subsection{Oriented Artinian level (AL) algebras.}
The class of Artinian level algebras generalizes the class of Artinian Gorenstein algebras.
Recall that the \emph{socle} of an Artinian algebra $A$ is the vector subspace 
$$\soc(A)=\left\{a\in A \ \left| \ m\cdot a=0, \ \forall \ m\in\mathfrak{m}_A\right.\right\}=(0:\mathfrak{m}_A).$$
\begin{definition}[Artinian level algebra]\label{leveldef}
The \emph{type} of an Artinian algebra $A$ is the vector space dimension of its socle. We say that a graded Artinian algebra $A$ is \emph{Artinian level}  (AL) if all elements of its socle have the same degree, i.e. $\soc(A)\subseteq A_d$, where $d=\max\{i \ | \ A_i \neq 0\}$.  In particular, $A$ is graded Artinian Gorenstein (AG) if and only if $A$ is level of type one. 
 \end{definition}
Let $A=\bigoplus_{i=0}^dA_i$ be an Artinian algebra with $A_d\neq 0$.  Then $A_d$ is a finite dimensional $\F$-vector space, say of dimension $n$.  We may choose a vector space isomorphism $A_d\cong\F^n$ and extend it by zero to an $\F$-linear map $\int_A\colon A(d)\rightarrow \F^n$.  We call the map $\int_A\colon A(d)\rightarrow\F^n$ an orientation on $A$.  Every orientation defines a \emph{generalized bilinear pairing} 
\begin{equation}
\label{eq:GBP}
\xymatrixrowsep{.5pc}\xymatrix{A\times A\ar[r] & \F^n\\
	(a,a')\ar@{|->}[r] & \int_Aa\cdot a'.\\}
\end{equation}
If we let $\pi_i\colon\F^n\rightarrow \F$ denote the projection map onto the $i^{th}$ coordinate, then the composite $\pi_i\left(a,a'\right)$ gives a bilinear pairing on $A$ in the usual sense.  We say that the generalized bilinear pairing is \emph{non-degenerate} if for every non-zero $a\in A$, there exists $a'\in A$ such that $\int_Aa\cdot a'\neq 0\in\F^n$.  
\begin{lemma}
	\label{lem:LevelND}
	The generalized bilinear pairing in Equation \eqref{eq:GBP} is non-de\-ge\-ne\-rate if and only if $A$ is level.  
\end{lemma}
\begin{proof}
	Suppose that $A$ is an Artinian level algebra.  We show that the pairing in Equation \eqref{eq:GBP} restricted to the subspace $A_i\times A_{d-i}\subset A\times A$ is non-degenerate for each degree $0\leq i\leq d$ by downward induction on $i$.  If $i=d$, then if $a\in A_d$ is non-zero then $\int_A a\cdot 1\neq 0$ since $\int_A\colon A_d\rightarrow \F^n$ is an isomorphism.  Next fix $i<d$ and assume that the pairing is non-degenerate for all $j$ such that $i<j\leq d$.  Fix $a\in A_i$ non-zero.  Since $A$ is level, we must have $\soc(A)=A_d$.  Therefore since $i<d$, $a\in A_i$ cannot be in the socle, hence there is a homogeneous element $x\in \mathfrak{m}_A=A_+$ for which $a\cdot x\neq 0$. 
	  Thus $\deg(a\cdot x)=j\geq i+1$ hence by our inductive hypothesis, there exists $y\in A_{d-j}$ for which $\int_A(a\cdot x)\cdot y\neq 0$ in $\F$.  If we take $a'=xy\in A_{d-i}$ we see that $\int_Aa\cdot a'\neq 0\in \F^n$, so the generalized pairing is non-degenerate.
	
	Conversely assume that the generalized bilinear pairing in Equation \eqref{eq:GBP} is non-degenerate.  Fix a non-zero homogeneous element $a\in \soc(A)$.  Suppose that $a\in A_i$ for some $i$.  By non-degeneracy of the generalized bilinear pairing, there exists $a'\in A_{d-i}$ for which $\int_Aa\cdot a'\neq 0\in \F^n$.  On the other hand if $d-i>0$ then $a'\in A_+$ from which it follows that $a\cdot a'=0$.  Thus $d-i=0$ and hence $i=d$.  This shows that $\soc(A)\subset A_d$ and thus $A$ is level.
\end{proof}

\begin{definition}
	\label{def:AL}
	A pair $\left(A,\int_A\right)$ consisting of a graded Artinian level algebra of socle degree $d$ and type $n$, and an orientation $\int_A\colon A(d)\rightarrow \F^n$ we shall call an \emph{oriented AL algebra}.  
\end{definition}

\subsubsection{Generalized Thom Classes.}
\begin{lemma}
	\label{lem:LevelThom}
	Suppose that $\left(L,\int_L\right)$ is an oriented AL algebra with socle degree $d$ of type $n$ and that $\left(K,\int_K\right)$ is an oriented AG algebra with socle degree $k$, and suppose that $\pi\colon L\rightarrow K$ is any algebra map between them.  Then there is a linear functional $\psi\colon \F^n\rightarrow\F$ and a homogeneous element $\tau\in L_{d-k}$ such that 
	$$\int_K\pi(y)=\psi\left(\int_L\tau\cdot y\right), \ \forall y\in L.$$  
\end{lemma}
\begin{proof}
	By Lemma \ref{lem:LevelND}, we see that the map 
	$$\xymatrixrowsep{.5pc}\xymatrix{L\ar[r] & \Hom_\F(L,\F^n)\\
		x\ar@{|->}[r] & \int_Lx\cdot(-)\\}$$
	is injective.  Hence there must be some linear map $\psi\colon\F^n\rightarrow\F$ such that the composition
	$$\xymatrix{L\ar[r]& \Hom_\F(L,\F^n)\ar[r]^-{\psi_*} & \Hom_\F(L,\F)}$$
	is injective, hence also an isomorphism.  Then since $\int_K\circ\pi\in\Hom_\F(L,\F)$, we deduce that there must be $\tau\in L_{d-k}$ such that 
	$$\int_K\circ\pi=\psi_*\circ\int_L\tau\cdot(-),$$
	and the result follows.
\end{proof}
Lemma \ref{cor:MDLevelGor} below is related, but uses Macaulay dual generators.
\begin{definition}
	\label{def:GenThom}
	The pair $\left(\psi,\tau\right)$ is a \emph{generalized Thom class} for the map $\pi\colon L\rightarrow K$.
\end{definition}

\section{Fibered Products and Connected sums.}
\label{sect:fpcs}

\subsection{Fibered Products.}

The fibered product is a particular instance of a general categorical construction termed pullback. The categories of (graded) rings, (graded) $\F$-algebras, and (graded) Artinian algebras 
are all closed under pullback (but not (graded) AG algebras, as we shall see in Lemma \ref{lem:Level}) and in each of these categories the pullback is called the fibered product. The existence of fibered products in the category of $\F$-algebras is closely related to  a dual concept to pushout in the dual category of (finite) affine schemes \cite[p.~87]{Ha}. In this paper we will work exclusively in the category of graded $\F$-algebras.  Here is
the formal definition.\par

\begin{definition}
\label{def:FiberedProduct}
	Given graded $\F$-algebras $A$, $B$, and $T$, and graded $\F$-algebra maps $\pi_A\colon A\rightarrow T$ and $\pi_B\colon B\rightarrow T$, the \emph{fibered product} of $A$ and $B$ over $T$ (with respect to $\pi_A$ and $\pi_B$) is the graded $\F$-subalgebra of $A\oplus B$
	$$A\times_T B=\left\{(a,b)\in A\oplus B \ \left| \ \pi_A(a)=\pi_B(b)\right.\right\}.$$
\end{definition}  
Let $\rho_1\colon A\times_TB\rightarrow A$ and $\rho_2\colon A\times_T B\rightarrow B$ be the natural projection maps.  It is well known that pullbacks satisfy the following 
universal property cf. \cite{Og}.
\begin{lemma}
	\label{lem:UnivProp}
	The fibered product $A\times_TB$ satisfies the following universal property:  If $C$ is another $\F$-algebra with maps $\phi_1\colon C\rightarrow A$ and $\phi_2\colon C\rightarrow B$ such that $\pi_A\circ\phi_1(c)=\pi_B\circ\phi_2(c)$ for all $c\in C$, then there is a unique $\F$-algebra homomorphism $\Phi\colon C\rightarrow A\times_TB$ which makes the diagram below commute:
	\begin{equation}
	\label{eq:UP}
	\xymatrix{C \ar@{-->}[dr]^{\Phi} \ar@/^1pc/[drr]^-{\phi_1}\ar@/_1pc/[ddr]_-{\phi_2}& & \\ & A\times_TB\ar[r]^-{\rho_1}\ar[d]_-{\rho_2} & A\ar[d]^-{\pi_A}\\ & B\ar[r]_-{\pi_B} & T.\\}
	\end{equation} 
\end{lemma}

{Note that if the maps $\pi_A$ and $\pi_B$ are surjective, then the fibered product fits into a short exact sequence of graded vector spaces:
\begin{equation}
\label{eq:SESFP}
\xymatrixcolsep{3pc}\xymatrix{0\ar[r] & A\times_T B\ar[r]^-{\rho_1\oplus \rho_2} & A\oplus B\ar[r]^-{\pi_A-\pi_B} & T\ar[r] & 0.\\}
\end{equation}
In this case we get a convenient formula for the Hilbert function:
\begin{lemma}[See also {\cite[Remark 3.2]{AAM}}]
	\label{lem:HFFP}
Assume that the maps $\pi_A$ and $\pi_B$ are surjective.  Then the Hilbert series of the fibered product $A\times_T B$ satisfies
	$$H(A\times_TB,t)=H(A,t)+H(B,t)-H(T,t).$$
\end{lemma}
\begin{proof}
	The result follows immediately from the exact sequence \eqref{eq:SESFP}.
\end{proof}
In fact in this case Lemma \ref{lem:UnivProp} implies that \eqref{eq:SESFP} characterizes the fibered product.
\begin{lemma}
	\label{lem:FPid}
	Assume that $\pi_A$ and $\pi_B$ are surjective, and suppose that $C$ is a graded $\F$-algebra with algebra maps $\phi_1\colon C\rightarrow A$ and $\phi_2\colon C\rightarrow B$.  Assume that 
	\begin{enumerate}[(a).]
		\item $\pi_A\circ\phi_1=\pi_B\circ \phi_2$, and that
		\item the sequence of maps 
		$$\xymatrixcolsep{3pc}\xymatrix{0\ar[r] & C\ar[r]^-{\phi_1\oplus \phi_2} & A\oplus B\ar[r]^-{\pi_A-\pi_B} & T\ar[r] & 0\\}$$
		is exact.
	\end{enumerate}
	Then $C\cong A\times_TB$ as $\F$-algebras.
\end{lemma}
\begin{proof}
	By the universal property (Lemma \ref{lem:UnivProp}), there exists an $\F$-algebra map $\Phi\colon C\rightarrow A\times_TB$ such that $\rho_i\circ \Phi=\phi_i$ for $i=1,2$.  Therefore, in the sequence of maps 
	$$\xymatrixcolsep{3pc}\xymatrixrowsep{2pc}\xymatrix{0\ar[r] & C\ar[d]_-{\Phi}\ar[r]^-{\phi_1\oplus \phi_2} & A\oplus B\ar[r]^-{\pi_A-\pi_B} & T\ar[r] & 0\\
	0\ar[r] & A\times_TB\ar[ur]_-{\rho_1\oplus\rho_2} & & &}$$
the triangle must commute.  But since $\phi_1\oplus\phi_2$ is injective, $\Phi\colon C\rightarrow A\times_TB$ must be injective too.  Finally by our exactness assumption, $C$ and $A\times_TB$ have the same Hilbert series, 
hence $\Phi$ must be an isomorphism.
\end{proof}

\begin{lemma}
	\label{lem:Level}
Assume that the maps $\pi_A$ and $\pi_B$ are surjective.  Then the fibered product $C = A \times_T B$ of graded AL algebras $A, B$ of the same socle degree $d$ and  types (socle dimensions) $s$, $t$, respectively, over a graded AG algebra $T$ of socle degree $k\leq d$ is an graded AL algebra of socle degree $d$ and type $s+t$ if $k<d$, and of type $s+t-1$ if $k=d$.
\end{lemma}
\begin{proof}
	Clearly $A\times_TB$ is Artinian, since it is a subalgebra of an Artinian algebra $A\oplus B$.  Suppose that $(a,b)\in A\times_TB$ is a non-zero element in the socle.  Then for any $(x,y)\in\mathfrak{m}_{A\times_TB}=\left(A\times_TB\right)_+$ we have $(x,y)\cdot(a,b)=(0,0)$.  Note that for any $x\in\mathfrak{m}_A=A_+$, $\pi_A(x)\in \mathfrak{m}_T=T_+$, hence by our surjectivity assumptions there exists $y\in \mathfrak{m}_B=B_+$ such that $(x,y)\in\left(A\times_TB\right)_+$.  Since $(x,y)\cdot(a,b)=(0,0)$, we deduce that $xa=0$.  Since $x\in A_+$ was arbitrary, this shows that $a\in A$ is in the socle of $A$.  A similar argument shows that $b\in B$ is in the socle of $B$.  Since the socle degrees of $A$ and $B$ are both equal to $d$ (again by our assumptions), we deduce that the degree of $(a,b)$ must also be equal to $d$.  Since $(a,b)$ was an arbitrary non-zero socle element, this argument shows that the socle of $A\times_TB$ is concentrated in degree $d$, and thus $C=A\times_TB$ is level. In particular this implies that $\operatorname{soc}(A\times_TB)=C_d$, and the socle dimensions can be computed using Lemma \ref{lem:HFFP}.
\end{proof}

\subsection{Connected Sums.}

For the rest of this section let  $\left(A,\int_A\right)$, $\left(B,\int_B\right)$, and $\left(T,\int_T\right)$ be oriented AG algebras with socle degrees $d$, $d$, and $k$, respectively.  Suppose we have surjective graded $\F$-algebra homomorphisms $\pi_A\colon A\rightarrow T$ and  $\pi_B\colon B\rightarrow T$, with Thom classes $\tau_A\in A_{d-k}$ and $\tau_B\in B_{d-k}$ respectively.

For the following definition we assume that $\pi_A(\tau_A)=\pi_B(\tau_B)$, so that $(\tau_A,\tau_B)\in A\times_TB$. This condition is clearly satisfied if the socle degrees $d$ of $A, B$, and $k$ of $T$, respectively, satisfy the inequality $d>2k$ since in that case $\pi_A(\tau_A)=\pi_B(\tau_B)=0$.
\begin{definition}
	\label{def:ConnectedSum}
	The connected sum of the oriented AG algebras $A$ and $B$ over $T$ (with respect to maps $\pi_A$ and $\pi_B$) is the quotient ring of the fibered product $A\times_TB$ by the principal ideal generated by the pair of Thom classes $(\tau_A,\tau_B)$, i.e. 
	$$A\#_TB=\frac{A\times_TB}{\left\langle \left(\tau_A,\tau_B\right)\right\rangle}.$$
	Note that this depends on $\pi_A$, $\pi_B$ and the orientations on $A$, $B$.
\end{definition}
Our definition for the connected sum is a special case of the notion by the same name defined in diagram (2.0.1) of \cite{AAM}. Indeed, the situation presented in this paper corresponds to the case where, in the notation of  \cite[(2.0.1)]{AAM}, $V$ is the canonical module of $T$. See Remark \ref{rk:dualGysin} for details on the maps appearing in that diagram.

The connected sum is also characterized by a short exact sequence of vector spaces: 

\begin{lemma}
	\label{lem:CSid}
	There is a short exact sequence of graded vector spaces
	\begin{equation}
	\label{eq:SESCS}
	\xymatrixcolsep{3pc}\xymatrix{0\ar[r] & T(k-d) \ar[r]^-{\times\tau_A\oplus\times\tau_B}& A\times_TB\ar[r] & A\#_TB\ar[r] & 0\\}
	\end{equation}	
	where the non-trivial map on the left is the direct sum of Gysin maps for $\pi_A$ and $\pi_B$, and the map on the right is the natural quotient map.  Moreover, if $C$ is another Artinian algebra with a surjective map $\phi\colon A\times_TB\rightarrow C$ making the sequence 
	$$\xymatrixcolsep{3pc}\xymatrix{0\ar[r] & T(k-d) \ar[r]^-{\times\tau_A\oplus\times\tau_B}& A\times_TB\ar[r]^-\phi & C\ar[r] & 0\\}$$
	exact, then $C\cong A\#_TB$ as $\F$-algebras.
\end{lemma}
\begin{proof}
	The image of the sum of Gysin maps $\times\tau_A\oplus\times\tau_B\colon T(k-d)\rightarrow A\oplus B$ is contained in the image of $\rho_1\oplus\rho_2\colon A\times_TB\rightarrow A\oplus B$ by our assumption that $\pi_A(\tau_A)=\pi_B(\tau_B)$.  Since $\pi_A$, $\pi_B$ are surjective, Lemma \ref{lem:InjSurj} implies that each Gysin map is injective, hence so is their direct sum.  This shows exactness on the left.  Exactness on the right is due to the natural projection being surjective.  To see exactness in the middle, it suffices to observe that the image of the direct sum of Gysin maps $\times\tau_A\oplus\times\tau_B\colon T(k-d)\rightarrow A\times_TB$ is exactly the principal ideal in $A\times_TB$ generated by the pair $(\tau_A,\tau_B)$, but this follows from the description of the Gysin map in Lemma \ref{lem:GysinMult}.
	
	Next, suppose that $C$ is an Artinian algebra with a map $\phi\colon A\times_TB\rightarrow C$ which fits into the $A\#_TB$ slot in the short exact sequence in \eqref{eq:SESCS}.  Then $\phi$ passes to an isomorphism $C\cong A\times_TB/\langle (\tau_A,\tau_B)\rangle= A\#_T B$ which is the desired isomorphism.
\end{proof}

\begin{lemma}[See also {\cite[Theorem 2.8]{AAM}}]
	\label{lem:Gorenstein}
Let $\left(A,\int_A\right)$, $\left(B,\int_B\right)$, and $\left(T,\int_T\right)$ be oriented AG algebras with socle degrees $d$, $d$, and $k$, respectively, and let  $\pi_A\colon A\rightarrow T$ and  $\pi_B\colon B\rightarrow T$ be surjective homomorphisms with Thom classes $\tau_A\in A_{d-k}$ and $\tau_B\in B_{d-k}$.  	Then the connected sum $A\#_TB$ is a (not necessarily standard) graded Artinian Gorenstein $\F$-algebra.
\end{lemma}

\begin{proof}
	It suffices to show that the socle is one dimensional.  Fix a homogeneous element  $c=\overline{(a,b)}\in\soc(A\#_TB)$.  Note that for each $x\in\ker(\pi_A)\subset A_+$, we have $(x,0)\in\left(A\times_TB\right)_+$, hence $\overline{(x,0)}\in\left(A\#_TB\right)_+$.  Thus we must have $\overline{(x,0)}\cdot\overline{(a,b)}=0$ which implies that $x\cdot a=\tau_A\cdot t$ for some $t\in T$, and also $0\cdot b=\tau_B\cdot t$ for the same $t\in T$.  Since $\pi_A$ is surjective, Lemma \ref{lem:InjSurj} implies that the multiplication map $\times \tau_A$, i.e. the Gysin map, is injective.  Hence we must have $t=0$ hence $x\cdot a=0$ for any $x\in\ker(\pi_A)$ and so  $a\in\left(0:\ker(\pi_A)\right)=\operatorname{Im}(\iota_A)$ by Remark \ref{rk:dualGysin}. Therefore we deduce that $a=\tau_A\cdot t_1$ for some $t_1\in T$.  A similar argument shows that $b=\tau_B\cdot t_2$ for some $t_2\in T$.  Note that $t_1=t_2$ if and only if $\overline{(a,b)}=0$ in $A\#_TB$.  Therefore we can replace the representative pair $(a,b)=\left(\tau_A\cdot t_1,\tau_B\cdot t_2\right)$ by $(a,b)-\left(\tau_A\cdot t_1,\tau_B\cdot t_1\right)$, and hence we may write $c=\overline{\left(0,\tau_B(t_2-t_1)\right)}$.  Then note that for any $y\in B_+$, and for any $x\in\pi_A^{-1}\left(\pi_B(y)\right)$ (which must exist by surjectivity of $\pi_B$), we have $\overline{(x,y)}\in\left(A\#_TB\right)_+$, and thus $\overline{(x,y)}\cdot c=\overline{(0,0)}$.  From this we deduce that $y\cdot\tau_B(t_2-t_1)=0$, hence it follows that $\tau_B(t_2-t_1)\in\soc(B)=B_d$.  This argument shows that 
	$$\soc\left(A\#_TB\right)\subset\left(A\#_TB\right)_d=\frac{\left(A\times_TB\right)_d}{\left\langle \left(\tau_A,\tau_B\right)\right\rangle_d}.$$ 
	Since $\left(A\times_TB\right)_d$ is two dimensional and $\left\langle\left(\tau_A,\tau_B\right)\right\rangle_d$ is one dimensional, we conclude that $\soc\left(A\#_TB\right)_d$ must be at most, hence exactly, one dimensional.    
\end{proof}

\begin{lemma}[See also {\cite[Theorem 3.3]{AAM}}]
	\label{lem:HFCS}
	The Hilbert series of the connected sum satisfies
	$$H(A\#_TB,t)=H(A,t)+H(B,t)-(1+t^{d-k})H(T,t).$$
	Equivalently, the Hilbert functions satisfy $H(A\#_TB)=H(A)+H(B)-H(T)-H(T)[d-k]$.
\end{lemma}
\begin{proof}
The result follows from Sequence \eqref{eq:SESCS} and Lemma \ref{lem:HFFP}.
\end{proof}

\begin{remark}
	We will see in Section \ref{sect:MDG} that the characterizing sequences \eqref{eq:SESFP} and \eqref{eq:SESCS} can be interpreted respectively as a Mayer-Vietoris sequence 
	\begin{equation}
	\label{eq:MVSeq}
	\xymatrix{0\ar[r] & Q/I_1\cap I_2\ar[r] & Q/I_1\oplus Q/I_2\ar[r]^-{\pi_1-\pi_2} & Q/I_1+I_2\ar[r] & 0}
	\end{equation}
	and a multiplication sequence:
	\begin{equation}
	\label{eq:mbtSeq}
	\xymatrix{0\ar[r] & Q/(I:\tau)(k-d)\ar[r]^-{\times\tau} & Q/I\ar[r] & Q/I+(\tau)\ar[r] & 0}
	\end{equation}
	where $Q$ is a polynomial ring, $I_1,I_2,I\subset Q$ homogeneous ideals, and $\tau\in Q_{d-k}$ a homogeneous polynomial. 
\end{remark}

\begin{example}[Fibered product and Connected Sum]\label{ex:FPEx}
	Let $A=\F[x,y]/(x^2,y^4)$ and $B=\F[u,v]/(u^3,v^3)$ each with the standard grading $\deg(x)=\deg(y)=\deg(u)=\deg(v)=1$.  Let $T=\F[z]/(z^2)$, and define maps $\pi_A:A\rightarrow T$, $\pi_A(x)=z, \ \pi_A(y)=0$ and $\pi_B\colon B\rightarrow T$, $\pi_B(u)=z, \ \pi_B(v)=0$.  Then the fibered product $A\times_TB$ is generated as an algebra by elements $z_1=(y,0)$, $z_2=(x,u)$, and $z_3=(0,v)$, all having degree one.  One can check that it has the following presentation:
	\begin{equation}
	\label{eq:FPEx}
	A\times_TB=\frac{\F[z_1,z_2,z_3]}{\left\langle z_1^4,z_2^3,z_3^3, z_1z_3,z_1z_2^2\right\rangle}.
	\end{equation}
	The Hilbert function of the fibered product is 
	\begin{align*}
	H(A\times_TB)= & (1,3,5,4,2)\\
	= & (1,2,2,2,1)+(1,2,3,2,1)-(1,1,0,0,0)\\
	= & H(A)+H(B)-H(T).
	\end{align*}
	
	Fix orientations on $A$, $B$, and $T$ by $\int_A\colon xy^3\mapsto 1$, $\int_B\colon u^2v^2\mapsto 1$, and $\int_T\colon z\mapsto 1$, respectively.  Then the Thom classes for $\pi_A\colon A\rightarrow T$ and $\pi_B\colon B\rightarrow T$ are, respectively, $\tau_A=y^3$, $\tau_B=uv^2$.  Note that $\pi_A(\tau_A)=0=\pi_B(\tau_B)$, hence $(\tau_A,\tau_B)\in A\times_TB$, and in terms of Presentation \eqref{eq:FPEx} we have $(\tau_A,\tau_B)=z_1^3+z_2z_3^2$.  Therefore we see that 
	\begin{equation}
	\label{eq:CSEx}
	A\#_TB=\frac{\F[z_1,z_2,z_3]}{\left\langle z_1^4,z_2^3,z_3^3, z_1z_3,z_1z_2^2,z_1^3+z_2z_3^2\right\rangle}.
	\end{equation}
	The Hilbert function of the connected sum is 
	\begin{align*}
	H(A\#_TB)= & (1,3,5,3,1)\\
	= & (1,2,2,2,1)+(1,2,3,2,1)-(1,1,0,0,0)-(0,0,0,1,1)\\
	= & H(A)+H(B)-H(T)-H(T)[3]
	\end{align*}
\end{example}

\begin{proposition}[See also {\cite[Proposition 2.4(b)]{ACLY}}]\label{prop:ConnSumF}
	Let $R,R'$ be graded polynomial rings with homogeneous maximal ideals $\mathfrak{m}_R$ and $\mathfrak{m}_{R'}$, respectively.  Let $\left(A=R/I,\int_A\right)$ and $\left(B=R'/I',\int_B\right)$ be oriented AG algebras each with socle degree $d$, and let $\pi_A\colon A\rightarrow\F$ and $\pi_B\colon B\rightarrow \F$ be the natural projection maps with Thom classes $\tau_A\in A_d$ and $\tau_B\in B_d$.  Then the fibered product $A\times_\F B$ has a presentation
	$$A\times_\F B\cong \frac{R\otimes_\F R'}{\mathfrak{m}_R\otimes\mathfrak{m}_{R'}+I\otimes R'+R\otimes I'}.$$
	and the connected sum $A\#_\F B$  has a presentation
	$$A\#_\F B\cong \frac{R\otimes_\F R'}{\mathfrak{m}_R\otimes\mathfrak{m}_{R'}+I\otimes R'+R\otimes I'+\left(\tau_A\otimes 1+1\otimes\tau_B\right)}.$$
	In particular, if $A$ and $B$ are standard graded then so are $A\times_\F B$ and $A\#_\F B$.

\end{proposition}
\begin{proof}

We first obtain presentations for $A$ and $B$ as quotients of $R\otimes_\F R'$
and subsequently use these to present the desired connected sum and fibered product. Indeed, we have 
\begin{align*}
A &\cong A\otimes_\F \F\cong R/I \otimes_\F R'/\m_{R'}\cong R\otimes_\F R'/I\otimes R'+R\otimes\m_{R'}\\
B &\cong \F\otimes_\F B \cong R/\m_R \otimes_\F R'/I'\cong R\otimes_\F R'/\m_R\otimes R'+R\otimes I'\\
F&\cong  \F\otimes_\F \F\cong R/\m_R\otimes_\F R'/\m_{R'}\cong R\otimes_\F R'/\m_R\otimes R'+R\otimes \m_{R'}.
\end{align*}

Notice that there is an equality of ideals
$$(I\otimes R'+R\otimes\m_{R'})+(\m_R\otimes R'+R\otimes I' )=\m_R\otimes R'+R\otimes \m_{R'}$$
which leads to the following short exact sequence
\begin{equation*}
0\to R\otimes_\F R'/(I\otimes R'+R\otimes\m_{R'})\cap(\m_R\otimes R'+R\otimes I' )\to A\oplus B \to \F \to 0.
\end{equation*}
Comparing the above to the  short exact sequence given in Equation (\ref{eq:SESFP})
\begin{equation*}
0\to A\times_\F B \to A\oplus B \stackrel{\pi_A-\pi_B}{\longrightarrow} \F\to 0
\end{equation*}
it follows that
\begin{equation}
\label{eq:pres1}
A\times_\F B \cong R\otimes_\F R'/(I\otimes R'+R\otimes\m_{R'})\cap(\m_R\otimes R'+R\otimes I').
\end{equation}
It remains to show that this simplifies to the claimed expression. To see this, notice that the containment
$$\mathfrak{m}_R\otimes\mathfrak{m}_{R'}+I\otimes R'+R\otimes I'\subseteq (I\otimes R'+R\otimes\m_{R'})\cap(\m_R\otimes R'+R\otimes I' )$$
is clear as $I\subseteq \m_R$ and $I'\subseteq\m_{R'}$. To show the reverse containment it is sufficient to note the following identities
\begin{align*}
&\frac{(I\otimes R'+R\otimes\m_{R'})\cap(\m_R\otimes R'+R\otimes I' )}{\mathfrak{m}_R\otimes\mathfrak{m}_{R'}}\\
 & =(I\otimes \F+ \F\otimes \m_{R'})\cap(\m_R\otimes \F+ \F \otimes I' )\\
& =I\otimes \F+\F \otimes I' \\
&=\frac{\mathfrak{m}_R\otimes\mathfrak{m}_{R'}+I\otimes R'+R\otimes I'}{\mathfrak{m}_R\otimes\mathfrak{m}_{R'}}.
\end{align*}
This gives the desired presentation for $A\times_\F B$. As for $A\#_\F B$, the given presentation follows from Definition \ref{def:ConnectedSum}. 

If $R$ and $R'$ are generated in degree 1 thens so is $R\otimes_\F R'$. This yields the assertion about $A\times_\F B, A\#_\F B$ being standard graded whenever $A$ and $B$ are.
\end{proof}

\section{Macaulay Dual Generators.}
\label{sect:MDG}
Let $Q=\F[x_1,\ldots,x_n]$ be a (not necessarily standard) graded polynomial ring and let $R=\F[X_1,\ldots,X_n]$ be the \emph{divided power algebra}, regarded as a $Q$ module with the contraction action $x_i\circ X_j^k=\begin{cases}X_j^{k-1}\delta_{ij} & \text{if} \ k>0\\ 0 & \text{otherwise}\\ \end{cases}$ where $\delta_{ij}$ is the Kronecker delta. Although we will regard $R$ and $Q$ as graded objects, the reader should keep in mind that the $Q$-module action on $R$ is not a \emph{graded} action since positive degree terms in $R$ \emph{lower} degrees of terms in $Q$, but we trust this will not cause too much confusion.

Then for each degree $i\geq 0$ the action of $Q$ on $R$ defines a non-degenerate $\F$-bilinear pairing 
\begin{equation}
\label{eq:MDPairing}
\xymatrixrowsep{.5pc}\xymatrix{Q_i\times R_i\ar[r] &  \F\\ 
	(f,F)\ar@{|->}[r] & f\circ F\\}
\end{equation}
This implies that for each $i\geq 0$ we have an isomorphism of $\F$-vector spaces $R_i\cong \Hom_\F(Q_i,\F)$ given by $F\mapsto\left\{f\mapsto f\circ F\right\}$.

The following fact is a classical result of Macaulay, cf. \cite[Theorem 21.6]{Ei}.  See also \cite[Proposition 2.67]{HMMNWW}.

\begin{fact}
	\label{fact:MDElts}
	An Artinian algebra $A=Q/I$ is level with socle degree $d$, and type $m$ if and only if there exists $m$ linearly independent homogeneous forms of degree $d$, $G_1,\ldots,G_m\in R_d$ such that $I=\Ann(G_1,\ldots,G_m)$, meaning that $I$ is the annihilator of the $Q$ submodule of $R$ generated by $G_1,\ldots,G_m$.
	
	In particular, if $A$ is Gorenstein with socle degree $d$, then $I=\Ann_Q(F_A)$ for some homogeneous polynomial $F_A\in R_d$.  Moreover this polynomial is unique up to a scalar multiple.   
\end{fact}

\begin{definition}
	\label{def:MDElts}
	The polynomials $G_1,\ldots,G_m$ or $F_A$ in Fact \ref{fact:MDElts} are referred to as a set of \emph{Macaulay dual generators} of the graded Artinian level algebra $A$.
\end{definition}

\begin{remarks}
	\label{rmk:MDOrient}
	\begin{enumerate}[(a)] 
	\item
	Note that an orientation on an Artinian level algebra $A=Q/I$ of socle degree $d$ and type $m$ is determined by fixing a set of Macaulay dual generators $G_1,\ldots,G_m\in R_d$:  
	$$\int_A a+\Ann(G_1,\ldots,G_m)\coloneqq (a\circ G_1(0),\ldots,a\circ G_m(0))^t\in\F^m, \ \ \forall \ a\in Q$$
	where $a\circ G_i(0)$ means evaluate the polynomial $a\circ G_i\in R$ at $X_1=0,\ldots,X_n=0$.
	Conversely, given an oriented AL algebra $\left(A,\int_A\right)$, and a presentation $A=Q/I$, there is a choice of Macaulay dual generators $G_1,\ldots,G_m$ such that $I=\Ann(G_1,\ldots,G_m)$ and $\int_A a+I=(a\circ G_1(0),\ldots,a\circ G_m(0))^t\in\F^m$ for all $a\in Q$.
	
	From now on, when we speak of Macaulay dual generators of an \emph{oriented} AL algebra $\left(A,\int_A\right)$ of type $m$ and socle degree $d$, we shall mean a set of $m$ homogeneous $d$-forms $G_1,\ldots,G_m\in R_d$ such that $A=Q/\Ann(G_1,\ldots,G_m)$ and
	$$\int_Aa+\Ann(G_1,\ldots,G_m)=(a\circ G_1(0),\ldots,a\circ G_m(0))^t, \ \ \forall a\in Q.$$  
	\item  Note that given a homomorphism  $\pi: A\to T$ of AG algebras having dual generators $F, H$ of degrees  $d$ and $k$, respectively, the Thom class of Definition~\ref{def:ThomClass} is the element $\tau$ of $A_{d-k}$ such that $\tau\circ F=H$. We sometimes regard $\tau$ as an element in $Q_{d-k}$, which is unique up to $\Ann (F)$.  Recall from Lemma \ref{lem:GysinMult} that the Gysin A-module map $\iota_A: T\to A$ is $\iota(t)=t\cdot \tau\in A$.
	\end{enumerate}
\end{remarks}

The following is essentially a direct consequence of Lemma \ref{lem:LevelThom}.
\begin{lemma}
	\label{cor:MDLevelGor}
	Suppose that $\left(L,\int_L\right)$ is an oriented AL algebra of socle degree $d$ and type $m$, suppose that $\left(K,\int_K\right)$ is an oriented AG algebra of socle degree $k\leq d$, and suppose that $\pi\colon L\rightarrow K$ is a surjective algebra map between them.  If $L$ has a presentation $L=Q/I$ with Macaulay dual generators $G_1,\ldots,G_m\in R_d$, then $K$ has a presentation $K=Q/J$ with Macaulay dual generator $F$ given by 
	\begin{equation}
	\label{eq:MDLG}
	F=\sum_{i=1}^ma_i\cdot\left(\tau\circ G_i\right)
	\end{equation}
	for some constants $a_1,\ldots,a_m\in\F$ and for some homogeneous polynomial $\tau\in Q_{d-k}$.  Moreover with these presentations, the generalized Thom class for the map $\pi\colon L\rightarrow K$ is the pair $(\psi\coloneqq (a_1,\ldots,a_m),\tau+\Ann(G_1,\ldots,G_m))$.
	
\end{lemma}

\begin{proof}
	Let $\left(L,\int_L\right)$, $\left(K,\int_K\right)$ and $\pi\colon L\rightarrow K$ be as above, and suppose that $L=Q/I$ has Macaulay dual generators $G_1,\ldots,G_m\in R_d$.  Then since $\pi\colon L\rightarrow K$ is surjective, $K$ also has a presentation of the form $K=Q/J$ for some $J\supseteq I$.  Since $K$ is Artinian Gorenstein of socle degree $k$, there is some homogeneous form $F\in R_k$ for which $J=\Ann(F)$, by Fact~\ref{fact:MDElts}.  To complete the proof we need to find constants $a_1,\ldots,a_m\in\F$ and a homogeneous polynomial $\tau\in Q_{d-k}$ such that Equation \eqref{eq:MDLG} holds.  Let $\left(\psi,\tau_0\right)$ be the generalized Thom class for the map $\pi\colon L\rightarrow K$, where $\psi\colon \Hom_\F(\F^k,\F)$ and $\tau_0\in L_{d-k}$.  Choose a basis for $L_d$ so that $\psi=(a_1,\ldots,a_m)$ and choose $\tau\in Q_{d-k}$ such that $\tau_0=\tau+I$.  Then for $y\in Q$, set $y_0=y+I\in L$, and compute: 
	\begin{align*}
	\int_K \pi(y_0)= & \psi\left(\int_L\tau_0\cdot y_0)\right)\\
	= & (a_1,\ldots,a_m)\cdot \left((y\cdot\tau)\circ G_1(0),\ldots,(y\cdot\tau)\circ G_m(0)\right)^t\\
	= & \sum_{i=1}^ma_i(y\cdot\tau)\circ G_i(0)\\
	= & y\circ\left(\sum_{i=1}^ma_i\cdot \tau\circ G_i\right)(0).  
	\end{align*}
	Set $F=\sum_{i=1}^ma_i\cdot \tau\circ G_i\in R_k$.  It remains to show that $J=\Ann(F)$.
	Suppose that $x\in \Ann(F)$, and set $x_0=x+I\in L$.  Then for any other $z\in Q$ with $z_0=z+I$, we have
	$$\int_K\pi(z_0)\cdot\pi(x_0)=0$$
	which implies that $\pi(x_0)=x+J=0+J$, since $K$ is Gorenstein.  Therefore $x\in J$.  Conversely, suppose that $x\in J$.  Then again for any other $z\in Q$, we must have 
	$xz\in\Ann(F)$, and taking $z=1$ shows that $x\in\Ann(F)$.  Hence $J=\Ann(F)$ and the proof is complete.	
\end{proof}

\begin{example}[Macaulay dual generators]
\label{ex:FP3}
Let 
$$A=\F[x]/(x^4), \ B=\F[u,v]/(u^3,v^2), \ T=\F[z]/(z^2),$$  of Hilbert functions $H(A)=(1,1,1,1)$ and $H(B)=(1,2,2,1)$.
Define maps $\pi_A\colon A\rightarrow T$, $\pi_A(x)=z$ and $\pi_B\colon B\rightarrow T$, $\pi_B(u)=z$, $\pi_B(v)=0$.  Then the fibered product has the presentation -- here $z_1,z_2$ have degree one, and $z_3$ has degree two --
$$A\times_TB=\frac{\F[z_1,z_2,z_3]}{\left(z_1^4,z_2^2,z_3^2,z_1z_3,z_1^2z_2-z_2z_3\right)}, \ \ \text{where} \begin{cases} z_1= & (x,u)\\
z_2= & (0,v)\\
z_3= & (0,u^2)\\
\end{cases}$$
A set of Macaulay dual generators for $A\times_TB$ are given by  the homogeneous degree
three forms
$$G_1=Z_1^3, \ \text{and} \ G_2=Z_1^2Z_2+Z_2Z_3,$$  
and $A\times_TB$ has Hilbert function $(1,2,3,2)=H(A)+H(B)-H(T)$.

Set $Q=\F[z_1,z_2,z_3]$ with variable degrees $1,1,2$. The natural projection maps given by $\phi_A\colon A\times_T B\rightarrow A$ and $\pi_B\colon A\times_TB\rightarrow B$ give $A$ and $B$ the new presentations
$$A\cong \frac{\F[z_1,z_2,z_3]}{(z_1^4,z_2,z_3)}, \ \text{and} \ B\cong \frac{\F[z_1,z_2,z_3]}{(z_1^4,z_2^2,z_3^2,z_1z_3,z_1^2z_2-z_2z_3,z_1^2-z_3,z_1^3)}= \frac{\F[z_1,z_2,z_3]}{(z_1^3,z_2^2,z_1^2-z_3)}$$
which have Macaulay dual generators, respectively,
$$F_A=G_1=Z_1^3 \ \text{and} \ F_B=G_2=Z_1^2Z_2+Z_2Z_3.$$
Also composing projection maps, say $\pi_A\circ \phi_A\colon A\times_TB\rightarrow T$ we get a presentation for $T$ as 
$$T=\frac{\F[z_1,z_2,z_3]}{(z_1^2,z_2,z_3)}=\frac{\F[z_1,z_2,z_3]}{\Ann(Z_1)}.$$
With these orientations, we see that the Thom classes for the maps $\pi_A\colon A\rightarrow T$ and $\pi_B\colon B\rightarrow T$ are given by
$$\tau_A=z_1^2+\Ann(Z_1^3)\in A_1, \ \ \text{and} \ \ \tau_B=z_1z_2+\Ann(Z_1^2Z_2+Z_2Z_3),$$
and we have $\pi_A(\tau_A)=0=\pi_B(\tau_B)$.  Set $\tau=(z_1^2-z_3)+z_1z_2$ $\in Q_2$ so that we have 
\begin{align*}\tau_A&=\tau+\Ann(F_A), \ \ \text{and} \ \ \tau_B=\tau+\Ann(F_B);\\
\tau\circ G_1&=\tau\circ G_2=Z_1, \text{ and }\Ann(\tau\circ G_1)=\Ann (G_1)+\Ann( G_2)=(z_1^2,z_2,z_3).
\end{align*}
Then $\tau+\Ann(G_1,G_2)=(\tau_A,\tau_B)\in A\times_TB$. We then have a presentation for the connected sum $C=A\#_TB=A\times_T B/(\tau)$,
whence
\begin{align*}
A\#_TB \cong & \frac{\F[z_1,z_2,z_3]}{\left(z_1^4,z_2^2,z_3^2,z_1z_3,z_1^2z_2-z_2z_3,(z_1^2-z_3)+z_1z_2\right)}\cong \frac{\F[z_1,z_2]}{(z_1^3+z_1^2z_2,z_2^2)}
\\
= & \frac{\F[z_1,z_2,z_3]}{\Ann(Z_1^3-(Z_1^2Z_2+Z_2Z_3))}=  \frac{\F[z_1,z_2,z_3]}{\Ann(F_A-F_B)}
\end{align*}
has Hilbert function $H(C)=(1,2,2,1)=H(A)+H(B)-H(T)-H(T)[1]$ as in Lemma \ref{lem:HFCS}.  It is interesting to note that the connected sum $A\#_TB$ has a standard grading whereas the fibered product $A\times_TB$ does not.  Also interesting is that $A\#_TB$ is a complete intersection, an anomaly that does not occur for connected sums over a field $T=\F$, \cite[Theorem 8.3]{AAM}.
\end{example}

\subsection{Characterization of the connected sum.}
We generalize Example \ref{ex:FP3}, to characterize the connected sum via Macaulay duality. This is our first main result, Theorem \ref{thm:A} from the Introduction, which we restate here for the convenience of the reader. 

\begin{theorem}\label{theorem1}

	Let $Q=\F[x_1,\ldots,x_n]$ be a (possibly non-standard) graded polynomial ring, and let $R=\F[X_1,\ldots,X_n]$ be its dual ring (a divided power algebra).
	Let $F,G\in R_d$ be two linearly independent homogeneous forms of degree $d$, and suppose that there exists $\tau\in Q_{d-k}$ (for some $k<d$) satisfying
	\begin{enumerate}[(a)]
		\item $\tau\circ F=\tau\circ G\neq 0$, and 
		\item $\Ann(\tau\circ F=\tau\circ G)=\Ann(F)+\Ann(G)$.
	\end{enumerate} 
	Define the oriented AG algebras 
	$$A=\frac{Q}{\Ann{F}}, \ B=\frac{Q}{\Ann(G)}, \ T=\frac{Q}{\Ann(\tau\circ F=\tau\circ G)},$$
	and let $\pi_A\colon A\rightarrow T$ and $\pi_B\colon B\rightarrow T$ be the natural projection maps.  Then the Thom classes of $\pi_A$ and $\pi_B$ are given by $\tau_A=\tau+\Ann(F)$ and $\tau_B=\tau+\Ann(G)$, and we have algebra isomorphisms
	$$A\times_TB\cong \frac{Q}{\Ann(F)\cap\Ann(G)}, \ \ A\#_TB\cong \frac{Q}{\Ann(F-G)}.$$
	And, conversely, every connected sum $A\#_T B$ of graded AG algebras of the same socle degree over a graded AG algebra $T$ arises in this way.
\end{theorem}   

\begin{proof}
	For the forward direction, let $F,G\in R_d$, $\tau\in Q_{d-k}$, $A$, $B$, $T$, $\pi_A\colon A\rightarrow T$ and $\pi_B\colon B\rightarrow T$ be as in the statement of the Theorem.  For the fibered product we look at the Mayer-Vietoris sequence \eqref{eq:MVSeq} with $I_1=\Ann(F)$ and $I_2=\Ann(G)$:
	\begin{equation}
	\label{eq:MVProof}
	\xymatrixcolsep{3pc}\xymatrix{0\ar[r] & \frac{Q}{\Ann(F)\cap\Ann(G)}\ar[r] & \frac{Q}{\Ann(F)}\oplus \frac{Q}{\Ann(G)}\ar[r]^-{\pi_A-\pi_B} & \frac{Q}{\Ann(F)+\Ann(G)}\ar[r] & 0.}
	\end{equation} 
	Note that by Condition (b), the rightmost term in \eqref{eq:MVProof} is $T$, and thus by Lemma \ref{lem:FPid}, we must have an algebra isomorphism 
	$$A\times_TB=\frac{Q}{\Ann(F)\cap\Ann(G)}=\frac{Q}{\Ann(F,G)}.$$
	For the Thom classes, note that for any $y\in Q_k$ we have
	$$(y\cdot \tau)\circ F=y\circ(\tau\circ F), \ \text{and} \ (y\cdot\tau)\circ G=y\circ(\tau\circ G)$$
	which shows that $\tau_A=\tau+\Ann(F)$ and $\tau_B=\tau+\Ann(G)$.
	Hence the total Thom class $(\tau_A,\tau_B)\in A\times_TB$ is identified with $\tau+\Ann(F,G)\in Q/\Ann(F,G)$.  
	
	For the connected sum, consider the multiplication by $\tau$ sequence \eqref{eq:mbtSeq}, where $I=\Ann(F,G)$:
	\begin{equation}
	\label{eq:mbtProof}
	\xymatrixcolsep{3pc}\xymatrix{0\ar[r] & \frac{Q}{\left(\Ann(F,G):\tau\right)}(k-d)\ar[r]^-{\times\tau} & \frac{Q}{\Ann(F,G)}\ar[r] & \frac{Q}{\Ann(F,G)+(\tau)}\ar[r] & 0.}
	\end{equation} 
	Looking at the leftmost term, note that we have 
	\begin{align*}
	(\Ann(F,G):\tau)= & \left(\Ann(F)\cap\Ann(G):\tau\right)\\
	= & \left(\Ann(F):\tau\right)\cap\left(\Ann(G):\tau\right)\\
	= & \Ann(\tau\circ F)\cap\Ann(\tau\circ G)\\
	= & \Ann(\tau\circ F=\tau\circ G).
	\end{align*}
	Hence the leftmost term is $T$, and we have already seen that the middle term is $A\times_TB$, hence by Lemma \ref{lem:CSid} the rightmost term must be isomorphic to the connected sum $A\#_TB$, i.e.
	$$A\#_TB\cong \frac{Q}{\Ann(F,G)+(\tau)}.$$
	It remains to see that $\Ann(F,G)+(\tau)=\Ann(F-G)$.  But since $A\#_TB$ is Gorenstein, Lemma~\ref{cor:MDLevelGor} implies that there are constants $a,b\in\F$ for which $\Ann(F,G)+(\tau)=\Ann(aF-bG)$.  On the other hand, since $\tau\in\Ann(aF-bG)$, condition (a) guarantees that $a=b=1$, and the proof of the forward implication is complete.
	
	For the converse, suppose that $\left(A,\int_A\right)$, $\left(B,\int_B\right)$ are oriented AG algebras of socle degree $d$, and $\left(T,\int_T\right)$ is an oriented AG algebra with socle degree $k<d$, such that $\pi_A\colon A\rightarrow T$ and $\pi_B\colon B\rightarrow T$ are surjective algebra maps with Thom classes $\tau_A\in A_{d-k}$ and $\tau_B\in B_{d-k}$, and consider the associated fibered product $A\times_TB$ and connected sum $A\#_TB$.  We follow the method of Example \ref{ex:FP3}.  Since $A\times_TB$ is an AL algebra, it has a presentation 
	$$A\times_TB\cong \frac{Q}{\Ann(H_1,H_2)}$$
	for some polynomial ring $Q$, with dual ring $R$, and some two linearly independent homogeneous $d$-forms $H_1,H_2\in R_d$.  Then the projection maps $\phi_A\colon A\times_TB\rightarrow A$ and $\phi_B\colon A\times_TB\rightarrow B$ give $A$ and $B$ presentations of the form $A=Q/I_A$ and $B=Q/I_B$.  By Lemma \ref{cor:MDLevelGor}, $A$ and $B$ have respective Macaulay dual generators $F=a_1H_1+a_2H_2$ and $G=b_1H_1+b_2H_2$ such that for all $a,b\in Q$
	$$\int_A(a+I_A)=a\circ F, \ \ \text{and} \ \ \int_B(b+I_B)=b\circ G.$$
	Note that $F$ and $G$ must be linearly independent.  Indeed because $A\times_TB$ contains elements of the form $(a,0)$ with $a\neq 0\in A$ (e.g. take $a\in A_d$ a socle generator, so that $\pi_A(a)=0$ since $k<d$), there exists $q\in Q$, $q+\Ann(H_1,H_2)\cong (a,0)$ for $\phi_A(q+\Ann(H_1,H_2))=q+\Ann(F)\cong a\neq 0$ but $\phi_B(q+\Ann(H_1,H_2))=q+\Ann(G)\cong 0$.  This implies that $F$ and $G$ are not scalar multiples of each other (since we found $q\in\Ann(F)\setminus\Ann(G)$), and hence they must be linearly independent.  Therefore they form a basis for the $\F$-span of $H_1,H_2$, and we can write
	$$A\times_TB=\frac{Q}{\Ann(H_1,H_2)}=\frac{Q}{\Ann(F,G)}.$$
		
	For $T$, note that the surjective map $\pi_A\colon A\cong \frac{Q}{\Ann(F)}\rightarrow T$ gives a presentation for $T$, $T=Q/J$, which is the same as the presentation given by the map $\pi_B\colon B\cong Q/\Ann(G)\rightarrow T$ by the commutativity of the diagram
	$$\xymatrix{A\times_TB\ar[r]\ar[d] & A\ar[d]\\
	B\ar[r] & T.}$$ 
	Therefore $J\supseteq \Ann(F),\Ann(G)$, and hence $J\supseteq\Ann(F)+\Ann(G)$.  Let $\iota\colon Q/\Ann(F)+\Ann(G)\rightarrow Q/J\cong T$ be the natural projection map.  Then comparing the Mayer-Vietoris sequence with the fibered product sequence 
	$$\xymatrixcolsep{3pc}\xymatrix{0\ar[r] & \frac{Q}{\Ann(F)\cap\Ann(G)}\ar[r]\ar[d] & \frac{Q}{\Ann(F)}\oplus\frac{Q}{\Ann(G)}\ar[r]\ar[d] & \frac{Q}{\Ann(F)+\Ann(G)}\ar[r]\ar[d]^-\iota & 0\\
		0\ar[r] & A\times_TB\ar[r] & A\oplus B\ar[r]_-{\pi_A-\pi_B} & T\ar[r] & 0}$$
	we see that since the left two vertical arrows are isomorphisms, so is $\iota$.  Hence $J=\Ann(F)+\Ann(G)$.  On the other hand, let $\tau_A\cong\tau_1+\Ann(F)$ and $\tau_B\cong \tau_2+\Ann(G)$ be the Thom classes for $\pi_A$ and $\pi_B$.  Since $\pi_A(\tau_A)=\pi_B(\tau_B)$, there is a total Thom class $(\tau_A,\tau_B)\in A\times_TB$, and hence there is a $\tau\in Q_{d-k}$ such that $(\tau_A,\tau_B)\cong \tau+\Ann(F,G)$; in particular we have $\tau_A=\tau+\Ann(F)$ and $\tau_B=\tau+\Ann(G)$.  We will show that $\tau$ satisfies conditions (a) and (b).
	First, by definition of Thom class we have for every $q\in Q$ 
	\begin{align*}
	\int_Aq\cdot\tau+\Ann(F)\coloneqq & \int_T\pi_A(q+\Ann(F))=\int_Tq+J\\
	= & (q\cdot \tau)\circ F=q\circ\left(\tau\circ F\right).
	\end{align*}
	Similarly we have 
	\begin{align*}
	\int_Bq\cdot\tau+\Ann(G)\coloneqq & \int_T\pi_B(q+\Ann(G))=\int_Tq+J\\
	= & (q\cdot \tau)\circ G=q\circ\left(\tau\circ G\right).
	\end{align*}
	Since $q\circ(\tau\circ F)=q\circ(\tau\circ G)$ for all $q\in Q$, we must have 
	\begin{enumerate}[(a)]
		\item $\tau\circ F=\tau\circ G\neq 0$.
	\end{enumerate}
	Also since $\left(T,\int_T\right)$ is oriented AG algebra with socle degree $k$ and presentation $T=Q/J$, there must be a homogeneous $k$-form $H\in R_k$ such that $J=\Ann(H)$ and $\int_Tq+\Ann(H)=q\circ H$ for all $q\in Q$.  But from the equations above, we must have $H=\tau\circ F=\tau\circ G$, and hence 
	\begin{enumerate}[(b)]
		\item $\Ann(\tau\circ F=H=\tau\circ G)=\Ann(F)+\Ann(G)=J$.
	\end{enumerate}
	Finally from the natural projection $\Phi\colon A\times_TB\cong Q/\Ann(F,G)\rightarrow A\#_TB$, we deduce from Lemma \ref{cor:MDLevelGor} that the connected sum $A\#_TB$ has a presentation of the form
	$$A\#_TB\cong\frac{Q}{\Ann(aF-bG)},$$
	and since $\tau\in\Ann(aF-bG)$, that $a=b=1$.
\end{proof}

\begin{example}\label{ex:cs}
	Let $Q=\F[x,y,z]$, $F_1=XY$, and $F_2=YZ$, and set 
	\begin{align*}
	A=Q/\Ann(XY)= & \frac{\F[x,y,z]}{(x^2,y^2,z)}\\
	B=Q/\Ann(YZ)= & \frac{\F[x,y,z]}{(x,y^2,z^2)}\\
	D=Q/\Ann(XY,YZ)= & \frac{\F[x,y,z]}{(x^2,y^2,z^2,xz)}\\
	C=Q/\Ann(XY-YZ)= & \frac{\F[x,y,z]}{(x^2,y^2,z^2,xz,x+z)}.
	\end{align*}
	Note that $\tau=x+z$ satisfies
	\begin{enumerate}[(a)]
		\item $\tau\circ F=Y=\tau\circ G\neq 0$, and 
		\item $\Ann(Y)=(x,y^2,z)=\Ann(F)+\Ann(G)$
	\end{enumerate}
	Set 
	$$T=Q/\Ann(Y)=\frac{\F[x,y,z]}{(x,y^2,z)}$$
	and let $\pi_A\colon A\rightarrow T$ and $\pi_B\colon B\rightarrow T$ be the natural maps.  Then the Thom classes are $\tau_A=x+\Ann(XY)$ and $\tau_B=z+\Ann(YZ)$.  The fibered product is given by 
	\begin{align*}
	A\times_TB= & \F[(x,0),(y,y),(0,z)]\\
	\cong & \frac{\F[x,y,z]}{(x^2,y^2,z^2,xz)}=D
	\end{align*}
	the total Thom class is $\tau=x+z$, and the connected sum is 
	\begin{align*}
	A\#_TB= & \frac{A\times_TB}{(\tau)}\cong \frac{\F[x,y,z]}{(x^2,y^2,z^2,xz,x+z)}=C.
	\end{align*}
	\end{example}
	
\begin{remark}
	\label{rmk:converse}
	Theorem \ref{theorem1} shows that the Macaulay dual generator of a connected sum $A\#_TB$ is \emph{always} the difference of the Macaulay dual generators of $A$ and $B$ over the \emph{correct} polynomial ring $\hat{Q}$.  Moreover, as indicated in the proof, this 
	\emph{correct} polynomial ring can always be determined from a Macaulay dual presentation of the fibered product $A\times_TB$.  The next Example \ref{ex:notcs} shows the importance of finding the correct ring $\hat{Q}$.
\end{remark}

\begin{example}\label{ex:notcs}
	Let $Q=\F[x,y]$, $F=X^2$, and $G=XY$, and set 
	\begin{align*}
	A=Q/\Ann(X^2)= & \frac{\F[x,y]}{(x^3,y)}\\
	B=Q/\Ann(XY)= & \frac{\F[x,y]}{(x^2,y^2)}\\ 
	D=Q/\Ann(X^2,XY)= & \frac{\F[x,y]}{(x^3,y^2,x^2y)}, \ \text{and}\\
	C=Q/\Ann(X^2-XY)= & \frac{\F[x,y]}{(x^3,y^2,x^2y,x^2+xy)}.
	\end{align*}
	Note that $\tau=x^2+xy\in Q_2$ satisfies 
	\begin{enumerate}[(a)]
		\item $\tau\circ F=1=\tau\circ G\neq 0$
	\end{enumerate} but it does not satisfy condition (b), i.e.
	$$\Ann(1)=(x,y)\neq(x^2,y)=\Ann(F)+\Ann(G).$$
	Here if we set  
	$$T=Q/\Ann(1)=\frac{\F[x,y]}{(x,y)}=\F,$$
	then a simple Hilbert function computation shows that $C\neq A\#_TB$:
	
	\begin{align*}
	H(A\#_TB)= & H(A)+H(B)-H(T)-H(T)[d-k=2]\\
	= & (1,1,1)+(1,2,1)-(1,0,0)-(0,0,1)\\
	= & (1,3,1)\\
	\neq & (1,2,1)=H(C).
	\end{align*}
	On the other hand, note that given $A$, $B$, and $T$ as above, and the natural projection maps $\pi_A\colon A\rightarrow T$ and $\pi_B\colon B\rightarrow T$, the associated fibered product is given by
	$$A\times_TB= \F[(x,0),(0,x),(0,y)]\cong \frac{\F[t_1,t_2,t_3]}{(t_1^3,t_2^2,t_3^2,t_1t_2,t_1t_3)}=\hat{Q}/\Ann(T_1^2,T_2T_3)$$
	where $\hat{Q}=\F[t_1,t_2,t_3]$. The projection maps $\phi_A\colon A\times_TB\rightarrow A$ and $\phi_B\colon A\times_TB\rightarrow B$ have kernels $K_A=\left(t_2,t_3\right)$ and $K_B=(t_1)$ respectively.  Thus $A$ and $B$ have presentations
	$$A=\frac{\F[t_1,t_2,t_3]}{(t_1^3,t_2,t_3)}=\hat{Q}/\Ann(T_1^2), \ \ B=\frac{\F[t_1,t_2,t_3]}{(t_1,t_2^2,t_3^2)}=\hat{Q}/\Ann(T_2T_3).$$
	In this case the Thom classes are just the socle generators $\tau_A=t_1^2$ and $\tau_B=t_2t_3$, and hence the total Thom class is $\tau=t_1^2+t_2t_3$.  Therefore, as predicted by Theorem \ref{theorem1}, the connected sum has presentation
	$$A\#_TB=\frac{F[t_1,t_2,t_3]}{(t_1^3,t_2^2,t_3^2,t_1t_2,t_1t_3,t_1^2+t_2t_3)}=\hat{Q}/\Ann(T_1^2-T_2T_3).$$
\end{example}

Theorem \ref{theorem1} also gives the well-known characterization of connected sums over $T=\F$, cf. \cite[Proposition II.2.4]{MS}. See also \cite[\S4]{SS}, \cite{BBKT} and \cite[\S3]{ACLY} for further work on the decomposability of AG algebras as connected sums over a field.
\begin{corollary}[{\cite[Proposition II.2.4]{MS}}]
	\label{cor:CSoverF}
	The oriented AL algebra $D=Q/\Ann(F,G)$ and the oriented AG algebra $C=Q/\Ann(F-G)$ are, respectively, the fibered product and the connected sum of $A=Q/\Ann(F)$ and $B=Q/\Ann(G)$ over $T=\F$, with the natural projection maps $\pi_A\colon A\rightarrow T$ and $\pi_B\colon B\rightarrow T$ if and only if $F$ and $G$ can be expressed (after possible change of coordinates) in disjoint sets of variables, i.e. there is a basis of linear forms $y_1,\ldots,y_n\in Q_1$ such that 
	$$y_i\circ F\neq 0 \Rightarrow y_i\circ G=0.$$
\end{corollary}
\begin{proof}
	Assume that $D$ and $C$ are, respectively, the fibered product and connected sum of $A$ and $B$ over $T=\F$, as above.  Then by Theorem \ref{theorem1} there is $\tau\in Q_d$ such that 
	\begin{enumerate}[(a)]
		\item $\tau\circ F=1=\tau\circ G$, and 
		\item $\Ann(1)=\mathfrak{m}_Q=\Ann(F)+\Ann(G)$.
	\end{enumerate}
	Choose any basis of linear forms for $Q_1\subset \mathfrak{m}_Q$, say $x_1,\ldots,x_n$.  For each $i$, condition (b) implies that we may write
	$x_i=y_i+z_i$ where $y_i\in\Ann(F)$ and $z_i\in\Ann(G)$.  Then the set of linear forms $\left\{y_1,\ldots,y_n,z_1,\ldots,z_n\right\}$ spans $Q_1$, and hence we may choose a basis for $Q_1$ from this spanning set, say 	$$\left\{\alpha_1,\ldots,\alpha_n\right\}, \ \text{where} \ \alpha_i\in\Ann(F) \ \text{or} \ \alpha_i\in\Ann(G).$$
	Then $\alpha_i\circ F\neq 0 \Rightarrow \alpha_i\circ G=0$, as desired.
	
	Conversely assume that $F,G\in R_d$ are two linearly independent $d$-forms, and suppose that $y_1,\ldots,y_n\in Q_1$ is a basis of linear forms such that 
	$$y_i\circ F\neq 0 \Rightarrow y_i\circ G=0.$$
	Let $f,g\in Q_d$ be polynomials for which $f\circ F=1$ and $g\circ G=1$ but $f\circ G=0$ and $g\circ F=0$ (such a choice is always possible as long as $F$ and $G$ are linearly independent).  Set $\tau=f+g$.  Then we have
	\begin{enumerate}[(a)]
		\item $\tau\circ F=1=\tau\circ G$, and
		\item $\Ann(1)=(y_1,\ldots,y_n)=\Ann(F)+\Ann(G)$.
	\end{enumerate}
	To see condition (b), note that $\Ann(F)+\Ann(G)$ contains all linear forms $y_i$ by our assumption on $F$ and $G$.
	Thus from Theorem \ref{theorem1} we deduce that $D=Q/\Ann(F,G)$ is the fibered product $A\times_TB$ and that $C=Q/\Ann(F-G)$ is the connected sum $A\#_TB$ for $T=Q/\Ann(\tau\circ F=1=\tau\circ G)=\F$. 
\end{proof}

\begin{definition}
\label{def:indec}
We say that a graded AG algebra $C$ is {\em decomposable} over a graded AG algebra $T$ if there exist AG algebras $A$ and $B$ with the same socle degree $d$ such that $C\cong A\#_TB$.  Otherwise we say that $C$ is {\em indecomposable} over $T$. We say that a graded AG algebra $C$ is {\em totally indecomposable} if it is indecomposable over every graded AG algebra $T$.
\end{definition}

 Decomposability for graded algebras as connected sums over a field was proposed as question in \cite[Problem 1.1.2]{MS}. Algebras defined by projective bundle ideals are in general indecomposable over $\F$ \cite[Corollary 4.3]{SS} and a criterion for indecomposability over a field is given in \cite{BBKT}. A large class of local AG algebras are shown to be indecomposable over $\F$ in \cite[Theorems 3.6,3.9]{ACLY}. The problem has been studied in special cases over $T$ in \cite[\S 8]{AAM}.  We address this also in our Corollary \ref{cor:BBKT} and Example \ref{ex:totallyindec} generalizing to indecomposability over $T$ a result of \cite{BBKT}.

\begin{lemma}[See also {\cite[Lemma 2.2]{SS3}}]
	\label{lem:Quadratic}
	Let $\F$ be a field of characteristic not equal two, and let $C=Q/\Ann(F)$ be an oriented AG algebra of socle degree two (i.e. $\deg(F)=2$). Then either $C\cong \F[z]/\Ann(Z^2)$, or $C$ is a connected sum of such algebras, i.e. $C\cong A_1 \#\cdots \# A_n$ where $A_i=\F[z_i]/\Ann(Z_i^2)$.
\end{lemma}
\begin{proof}
	Let $Q=\F[x_1,\ldots,x_n]$ and $F\in Q_2$ some quadratic form.  Since every quadratic form can be diagonalized by some change of coordinates, say $X_i\mapsto Z_i$, we  may write
	$$F=a_1Z_1^2+\cdots+a_nZ_n^2, \ \ \text{some $a_i\in\F$}.$$
	Hence by Corollary \ref{cor:CSoverF}, $C$ is a connected sum of the desired structure.
\end{proof}

\begin{remark}[Decomposability into a connected sum depends on the field]\label{decomposerem} Take a degree-$j$ form $q\in R=\F [X,Y]$ and let $C=Q/\Ann(q)$.  Assume $j\ge 3$, then by \cite{BBKT} $C$ is a connected sum over $T=\F$ if $I=\Ann F=(h,g)\subset Q$ where $h\in I_2$ factors over $\F$ as $h=\ell_1\cdot\ell_2$, and $g\in Q_j$.\footnote{The authors of \cite{BBKT} work over an algebraically closed field; that these statements extend to arbitrary fields is shown in \cite[Proposition 1.36]{IK}.}  Thus, to obtain $q$ not a direct sum over $\mathbb Q$, but a direct sum over a quadratic extension $\F$, take $q\in \mathbb Q[X,Y]$, such that $h\in I_2$ factors over $\F$ but not over $\mathbb Q$. We thank Z.~Teitler for the example $q=(X+iY)^4+(X-iY)^4=2\left(X^4-6X^2Y^2+Y^4\right)$; then $h=x^2+y^2$ can be factored over $\F=\mathbb Q[i]$ but not over $\mathbb Q$, so $C$ is not a connected sum over $\Q$. Or consider
 $$q=(\sqrt{3}X+Y)^4+(\sqrt{3}X-Y)^4=18(X^4+2X^2Y^2+Y^4),$$
 for which $h=x^2-3y^2=(x-\sqrt{3}y)(x+\sqrt{3}y)$; then $C$ 
is a connected sum over $\F=\mathbb Q(\sqrt{3})$ but not over $\mathbb Q$.\footnote{ For the Macaulay duality we here let $Q=\F[x,y]$ act on $R$ as differential operators (not contraction), so we use the usual powers in the polynomial ring $S$; if we used the contraction action of $R$ on $S$, we would need to replace the usual powers $L^j$ by divided powers.} \par
H. Huang, H. Lu, Y. Ye, and C. Zhang in \cite[Examples 5.3,5.4]{HLYZ} give 
forms in $\Q[X,Y,Z]$ and $\Q[X,Y,Z,W] $ whose decomposability as a connected sum depends on the field $\F$. 

Later, we propose a similar example, over a $T$ not equal to a field (Example 
\ref{CconnectFex}).
\end{remark}

\begin{example}[Indecomposable over $\F$, but decomposable over $T$]
	\label{ex:Alexandra}
	Define the AG algebras (with their standard gradings)
	\begin{align*}
	A=\F[x,y,z]/\Ann(X^2Y)= & \F[x,y,z]/(x^3,y^2,z),\\ 
	B=\F[x,y,z]/\Ann(Y^2Z)= & \F[x,y,z]/(x,y^3,z^2),\\
	C=\F[x,y,z]/\Ann(X^2Y-Y^2Z)= & \F[x,y,z]/(z^2,xz,y^3,xy^2,x^2+yz)\\
	D=\F[x,y,z]/\Ann(X^2Y,Y^2Z)= & \F[x,y,z]/(z^2,xy,y^3,xy^2)\\	
	T=\F[x,y,z]/\Ann (Y)= & \F[x,y,z]/(x,y^2,z)
	\end{align*} 
	and maps $\pi_A\colon A\rightarrow T$ and $\pi_B\colon B\rightarrow T$ by $\pi_A(y)=y=\pi_B(y)$ and $\pi_A(x)=\pi_A(z)=\pi_B(x)=\pi_B(z)=0$.  Then $F_1=X^2Y$, $F_2=Y^2Z$, and $\tau=x^2+yz$, and $\Ann(F_1-F_2)=\tau+\Ann(F_1,F_2)$.  Moreover $\Ann(F_1)+\Ann(F_2)=\Ann(\tau\circ F_1=\tau\circ F_2)=\Ann(Y)$.  Thus, conditions (a) and (b) from Theorem \ref{theorem1} are satisfied, hence $C$ is a connected sum $C=A\#_TB$.  On the other hand $C$ is indecomposable over $T=\F$ according to \cite[Example 1.4]{BBKT}.
\end{example}

	W. Buczy\'nska et al. showed that if an AG algebra $C=Q/\Ann(F)$ of socle degree $d$ is decomposable over $\F$ then the ideal $\Ann(F)$ must contain a minimal generator of degree $d=\deg(F)$ \cite[Theorem 1.1]{BBKT}.  One can derive the following corollary, which generalizes their result, from Theorem \ref{theorem1}.  

\begin{corollary}
	\label{cor:BBKT}
	If $C=Q/\Ann(F)$ is an AG algebra of socle degree $d$ that is decomposable over another AG algebra $T=Q'/\Ann(F')$ of socle degree $k$, then the ideal $\Ann(F)$ must contain a minimal generator of degree $d-k=\deg(F)-\deg(F')$. 
\end{corollary}
\begin{proof}
	Assume that $C=Q/\Ann(F)=A\#_TB$.  Then we can find $F_1$ and $F_2$ such that $F=F_1-F_2$ where $A=Q/\Ann(F_1)$ and $B=Q/\Ann(F_2)$, and by Theorem \ref{theorem1}, we must have $\Ann(F)=\Ann(F_1-F_2)=\tau+\Ann(F_1,F_2)$ for some non-zero $\tau\in Q_{d-k}$.  Thus in particular $\tau\in\Ann(F)$ is a minimal generator of degree $\deg(\tau)=d-k$.   
\end{proof}

\begin{example}[Indecomposable algebras]\label{ex:totallyindec}
\begin{enumerate}[(i).]
\item
 Any graded AG $\F$-algebra of embedding dimension one and socle degree at least two is totally indecomposable. Indeed, such an algebra has the form $C=\F[x]/(x^n)$ for some $n\geq 3$ by the structure theorem for modules over a PID. Since the Macaulay dual generator for $A$ is $F_C=X^{n-1}$, if $C$ were decomposable over some graded AG algebra $T$, Corollary \ref{cor:BBKT} would give that $n-1=(n-1)-k$, where $k$ is the socle degree of $T$. Therefore $k=0$ and $T=\F$. But if $C=A\#_\F B$ with $A$ and $B$ having socle degree at least three, then Lemma \ref{lem:HFCS} gives that $\dim_\F( C_1)\geq 2$, a contradiction.
 
 % - the next item is wrong
% \item Any standard graded AG $\F$ algebra $C$ of Hilbert function $H=(1,3,6,3,1)$ (or, more generally, any codimension three standard graded AG extremal Gorenstein $C$ of socle degree at least four)\footnote{An extremal AG algebra \cite{Sch} is one having the maximum Hilbert function $H(C)_i=\max\{ \dim Q_i, Q_{d-i}\}$ possible, for an even socle-degree $d$.}  is totally indecomposable. To see this for $H=(1,3,6,3,1)$, an AG $C=\F[x,y,z]/I$, note that the defining ideal $I$ has 7 generators, all of degree 3; by Corollary \ref{cor:BBKT},
% $d-k=4-k=3$ so $T$ has socle degree $1$. A check of the possible component Hilbert functions shows $H(A)=H(B)=(1,2,3,2,1)$,  but each ideal is generated in degree 3, and the ideal defining $T=\F[x]/(x^2)$ is generated in degree two, so condition (b) of Theorem \ref{theorem1} cannot be satisfied.
% 
% \item The non-standard graded $C=A\#_{F}B$ has $H(C)=(1,3,6,3,1)$ where $A$ is the standard-graded $A= F[y,z,w]/ \Ann(W^4+Y^3W+Z^2W^2)$ with $H(A)=(1,3,5,3,1)$, and $B=F[x]/(x^3)$ with $\deg (x)=2$ and $ H( B )= (1,0,1,0,1)$.
 
\item In contrast to the case over $T=\F$ where complete intersections are always indecomposable \cite[Theorem 8.3]{AAM}, over other $T$ this is no longer true.  For example 
let $Q=\F[x,y,z]$, $A=Q/\Ann(XYZ), B=Q/\Ann(Z^3), \tau=xy+z^2, T=Q/\Ann(Z)$ then, by Theorem \ref{theorem1}, $C=A\#_T B\cong Q/\Ann(Z^3-XYZ)=Q/(x^2,y^2,z^2+xy)$, which is a complete intersection.\par 
Now consider $C\#_T B$, $Q=\F[x,y,z]$, $\tau=z^2$, so $\tau\circ (Z^3-XYZ)=\tau \circ Z^3=Z,$ and
$(x,y,z^2)=\Ann(Z)=(z^2-xy,x^2,y^2)+(x,y,z^4)=\Ann(Z^3-ZXY)+\Ann(Z^3)$, so we have by Theorem \ref{theorem1} $C\#_T B=\Ann (Z^3-XYZ-Z^3)=\Ann (-XYZ)=\Ann (XYZ)=A$, so $A=C\#_TB$. \par
Over a field $T=\F$ each AG algebra of socle degree at least 3 has a unique decomposition into 
indecomposable summands \cite[Proposition 3.1]{SS3}. The above pair of decompositions $C=A\#_T B$ and $A=C\#_T B$ illustrate that the AG algebras of Hilbert function $H=(1,3,3,1)$ over $T=\F[z]/(z^2)$ form a class of decomposable algebras: we cannot write $C$ as the connected sum of indecomposable summands!

%In contrast to the case over $T=\F$ where complete intersections are always indecomposable \cite[Theorem 8.3]{AAM}, over higher $T$ this is no longer true.  For example 
%are that are of Hilbert function $H=(1,3,3,1)$ over $T=\F[z]/(z^2)$ form a class of decomposable algebras  but there are no indecomposable summands!
%let $Q=\F[x,y,z]$, $A=Q/\Ann(XYZ), B=Q/\Ann(Z^3), \tau=xy+z^2, T=Q/\Ann(Z)$ then, by Theorem \ref{theorem1}, $C=A\#_T B\cong Q/\Ann(Z^3-XYZ)=Q/(x^2,y^2,z^2+xy)$, which is a complete intersection. 
%But, likewise, we may write $A$ as $C\#_T B$.
 \end{enumerate}
\end{example}
\subsection{The monomial case.}
We will now describe the Macaulay dual of a connected sum $F-G$ in the special case where $F$ and $G$ are monomials.  First some notation.

Let $R=\F[X_1,\ldots,X_n]$ and $Q=\F[x_1,\ldots,x_n]$ be as above, and let $m=x_1^{a_1}\cdots x_n^{a_n}\in Q$ and $M=X_1^{b_1}\cdots X_n^{b_n}\in R$ be monomials.  Define their \emph{dual monomials} by 
\begin{align*}
m^*= & X_1^{a_1}\cdots X_n^{a_n}\in R\\
M^*= & x_1^{b_1}\cdots x_n^{b_n}\in Q
\end{align*} 
Note that 
\begin{equation}
\label{eq:Prop1}
m\circ M\neq 0 \ \Leftrightarrow \ \text{$m^*$ divides $M$ in which case} \ m^*\cdot \left(m\circ M\right)=M.
\end{equation}
Note that if $m$ and $M$ have the same degree, i.e. $a_1+\cdots+a_n=b_1+\cdots+b_n$ then 
$$m\circ M=\begin{cases} 1 & \text{if $m=M^*$}\\ 0 & \text{if $m\neq M^*$}.\\ \end{cases}$$
Define the support of a monomial to be those variables which divide it, i.e.
$$\Supp(m \ (\text{or} \ M))\coloneqq \left\{x_i \ (\text{or} \ X_i)\left|a_i \ (\text{or} \ b_i) \ >0\right.\right\}$$
First, an easy Lemma:
\begin{lemma}
	\label{lem:Monomial}
	The ideal $\Ann(F)$ is generated by monomials if and only if $F$ is a monomial.
\end{lemma}
\begin{proof}
	If $F=X_1^{a_1}\cdots X_n^{a_n}$ is a monomial, then $\Ann(F)=(x_1^{a_1+1},\ldots,x_n^{a_n+1})$ is a monomial ideal.  Conversely, suppose that $F$ has a monomial expansion with non-zero weights on at least two distinct monomials say $M_1$ and $M_2$ of degree $d$, i.e. 
	$$F=c_1M_1+c_2M_2+\left(\substack{\text{linear combination of other}\\ \text{monomials of degree $d$}}\right).$$  
	Then $c_2M_1^*-c_1M_2^*\in\Ann(F)$.  But if $\Ann(F)$ were a monomial ideal, then $M_1^*$ and $M_2^*$ would also be in $\Ann(F)$, but they are not, e.g. $M_1^*\circ F=c_1\neq 0$.  Therefore $\Ann(F)$ is not a monomial ideal.
\end{proof}

\begin{lemma}
	\label{lem:taubinom}
	Suppose that $F,G\in R_d$ are distinct monomials of degree $d$, and suppose that $\tau\in Q_{d-k}$ satisfies the conditions of Theorem \ref{thm:A}, i.e. 
	\begin{enumerate}[(a)]
		\item $\tau\circ F=\tau\circ G\neq 0$, and
		\item $\Ann(\tau\circ F=\tau\circ G)=\Ann(F)+\Ann(G)$.
	\end{enumerate}
	Then there exist monomials $m_F\in \Ann(G)_{d-k}$ and $m_G\in\Ann(F)_{d-k}$ such that 
	$$\tau\equiv m_F+m_G \ \text{mod} \Ann(F,G).$$
\end{lemma}

\begin{proof}
	Let $\tau\in Q_{d-k}$ be as above, and let 
	$$\tau=\sum_{i=1}^mc_im_i$$
	be its monomial expansion with coefficients $c_i$ and distinct monomials $m_i\in Q_{d-k}$.  Note that since $\Ann(F)$ and $\Ann(G)$ are monomial ideals, so is their sum $\Ann(F)+\Ann(G)$, and hence also by our assumption $\Ann(\tau\circ F=\tau\circ G)$.  Therefore $M_0\coloneqq \tau\circ F=\tau\circ G$ must be a monomial.  On the other hand, we have 
	\begin{equation}
	\label{eq:MEtau}
	\sum_{i=1}^mc_im_i\circ F=M_0=\sum_{i=1}^mc_im_i\circ G.
	\end{equation}
	Note that since the monomials $m_i$ are distinct, so are the monomials $m_i^*$, and hence by \eqref{eq:Prop1} so are the monomials $m_i\circ F$, as well as the monomials $m_i\circ G$.  In particular, the LHS and RHS of \eqref{eq:MEtau} are monomial expansions of the single monomial $M_0$.  Hence looking at this monomial expansion on the LHS of \eqref{eq:MEtau} we conclude that there is a unique index $i_F\in\left\{1,\ldots,m\right\}$ such that 
	$$c_im_i\circ F=\begin{cases} M_0 & \text{if $i=i_F$}\\ 0 & \text{if $i\neq i_F$}\\ \end{cases}$$
	and similarly for $G$ and after rescaling $\tau$, we may assume that $c_{i_F}=1=c_{i_G}$.  Then if we set $m_F\coloneqq m_{i_F}$ and $m_G\coloneqq m_{i_G}$, we see that 
	$$\tau-(m_F+m_G)\in\Ann(F)\cap\Ann(G)=\Ann(F,G).$$
	It remains to show that $m_F\in\Ann(G)$ and that $m_G\in\Ann(F)$.  To see this, it suffices to see that the indicies $i_F$ and $i_G$ are distinct.  But they must be, for if they were the same index, we would have by \eqref{eq:Prop1}
	$$F=m_{i_F}^*\cdot(m_{i_F}\circ F)=m_{i_F}^*\cdot M_0=m_{i_G}^*\cdot M_0=m_{i_G}^*\cdot (m_{i_G}\circ G)=G$$
	contradicting our assumption that $F$ and $G$ are distinct.  
\end{proof}

\begin{proposition}
	\label{prop:TFAE}
	Let $F,G\in R_d$ be two distinct monomials of the same degree $d$.  The following are equivalent:
	\begin{enumerate}[(i).]
		\item There exists $\tau\in Q_{d-k}$ satisfying 
		\begin{enumerate}[(a)]
			\item $\tau\circ F=\tau\circ G$, and
			\item $\Ann(\tau\circ F=\tau\circ G)=\Ann(F)+\Ann(G)$.
		\end{enumerate}
		\item If $M_0=\operatorname{gcd}(F,G)$ then $F\nmid M_0^2$ and $G\nmid M_0^2$.
	\end{enumerate}
	These imply that the algebra $C=Q/\Ann (F-G)$ is a connected sum $A\#_TB$ with $A={Q}/{\Ann{F}}, \ B={Q}/{\Ann(G)}$ over $T=Q/{\Ann(\tau\circ F=\tau\circ G)}.$

\end{proposition}

\begin{proof}
	Assume that (1) holds.  Then by Lemma \ref{lem:taubinom}, there exist monomials $m_F\in\Ann(G)$ and $m_G\in\Ann(F)$ such that $\tau\equiv m_F+m_G$ mod $\Ann(F,G)$.  Hence we may assume without loss of generality that $\tau=m_F+m_G$.  Set $M_0=\tau\circ F=m_F\circ F$ which by condition (a) is equal to $M_0=m_G\circ G=\tau\circ G$.  Set $M_F=m_F^*\in R_{d-k}$ and $M_G=m_G^*\in R_{d-k}$.  Then by \eqref{eq:Prop1} we have 
	$$F=M_F\cdot M_0, \ \ \text{and} \ \ G=M_G\cdot M_0.$$
	To see that $M_0=\operatorname{gcd}(F,G)$, set $m_0=M_0^*$ and note that for every variable $x_i$, we have $x_i\cdot m_0\in\Ann(M_0)=\Ann(F)+\Ann(G)$, by condition (b).  But since $x_i\cdot m_0$ is a monomial, and since $\Ann(F)$ and $\Ann(G)$ are monomial ideals, it follows that $x_i\cdot m_0\in\Ann(F)$ or $\Ann(G)$.  Also note that $(x_i\cdot m_0)\circ F=x_i\circ M_F$ and $(x_i\cdot m_0)\circ G=x_i\circ M_G$.  In particular it follows that $M_F$ and $M_G$ must be relatively prime, and therefore $M_0=\operatorname{gcd}(F,G)$.  Since $m_F\in\Ann(G)$ we see that $m_F\circ G=0$, and hence again by \eqref{eq:Prop1}, $M_F\nmid G=M_G\cdot M_0$, which in turn implies that $F=M_F\cdot M_0\nmid M_0^2$.  A similar argument shows that since $m_G\in \Ann(F)$ that $G\nmid M_0^2$ as well.

Next assume that (2) holds, i.e. let $M_0=\operatorname{gcd}(F,G)$ and assume that $F\nmid M_0^2$ and $G\nmid M_0^2$.  Set $m_F=M_F^*$ and $m_G=M_G^*$ and set $\tau\coloneqq m_F+m_G$.  Note first that $m_F\in \Ann(G)$.  Indeed if $m_F\circ G\neq 0$, then by \eqref{eq:Prop1} $M_F$ would divide $G=M_G\cdot M_0$.  But $M_F$ and $M_G$ are relatively prime since $M_0=\operatorname{gcd}(F,G)$, and hence we must have $M_G\mid M_0$ which implies that $G\mid M_0^2$, contradicting our assumption that (2) holds.  Therefore $m_F\in\Ann(G)$, and similarly $m_G\in \Ann(F)$.  Hence we must have 
	$$\tau\circ F=m_F\circ F=M_0=m_G\circ G=\tau\circ G$$
	which is condition (a).  To see that condition (b) holds, it suffices to see that for every monomial $m\in\Ann(M_0)=\Ann(\tau\circ F=\tau\circ G)$, that either $m\in\Ann(F)$ or $m\in\Ann(G)$.  To this end, fix $m\in\Ann(M_0)$, and suppose that $m\notin\Ann(F)$.  Then $m\circ F\neq 0$ and hence by \eqref{eq:Prop1}, the monomial $M=m^*$ must divide $F=M_F\cdot M_0$.  Also, since $m\circ M_0=0$, $M$ must not divide $M_0$.  But this implies that $M$ must not divide $G$ either, since otherwise $M$ would be a common divisor of $F$ and $G$, and would necessarily divide their greatest common divisor, $M_0$ which it does not.  It follows from \eqref{eq:Prop1} that $m\in \Ann(G)$, as desired.  This shows that , which shows that condition (b) holds and thus (1) holds.  
	The last statement is from Theorem \ref{theorem1}.
\end{proof}

\begin{example}
	\label{ex:easy}
	Let $Q=\F[x,y,z]$, and set 
	$$F=X^3Y, \ \ G=XZ^3.$$
	Then $M_0=X=\operatorname{gcd}(F,G)$ and clearly $F\nmid M_0^2$ and $G\nmid M_0^2$ so that condition (2.) from Proposition \ref{prop:TFAE} is satisfied.
Set 
	\begin{alignat*}{3}
	A &= & \frac{Q}{\Ann(X^3Y)} &= & \frac{\F[x,y,z]}{(x^4,y^2,z)}\\
	B &= & \frac{Q}{\Ann(XZ^3)} &= & \frac{\F[x,y,z]}{(x^2,y,z^4)}\\
	T &= & \frac{Q}{\Ann(X)} &= & \frac{\F[x,y,z]}{(x^2,y,z)}
	\end{alignat*}	
	with the natural surjections $\pi_A,\pi_B\colon A,B\rightarrow T$.  Then $\tau_A=x^2y$ and $\tau_B=z^3$, and the total Thom class is $\tau=x^2y+z^3$.  The fibered product is 
	\begin{align*}
	A\times_TB\cong & \F[(x,x),(y,0),(0,z)]\cong \frac{\F[x,y,z]}{(x^4,y^2,z^4,yz,x^2z)}
	= & \frac{Q}{\Ann(X^3Y,XZ^3)}
	\end{align*}
	and the connected sum is given by 
	\begin{align*}
	A\#_TB= & \frac{A\times_TB}{(\tau=x^2y+z^3)}\cong \frac{\F[x,y,z]}{(x^4,y^2,z^4,yz,x^2z,x^2y+z^3)}
	= & \frac{Q}{\Ann(X^3Y-XZ^3)}
	\end{align*}
	as implied by Proposition \ref{prop:TFAE}.
\end{example}
In the next example, we see that a change of $A,B$ to $A^\prime=R/\Ann F^\prime, B^\prime=R/\Ann G^\prime$, may be needed to write
	 $C=R/\Ann(F-G)$ as a connected sum.
\begin{example}
	\label{ex:Tony}
	Let $Q=\F[x,y]$, and set 
	$$F=X^2Y, \ \ G=XY^2.$$
	Then $M_0=\operatorname{gcd}(F,G)=XY$, but in this case we have both $F\mid M_0^2$ and $G\mid M_0^2$.  
	We can verify that $F-G=X^2Y-XY^2=XY(X-Y)$ is not a connected sum of $F$ and $G$:
	Set 
	\begin{alignat*}{3}
	A &= & \frac{Q}{\Ann(X^2Y)} & = & \frac{\F[x,y]}{(x^3,y^2)}\\
	B &= & \frac{Q}{\Ann(XY^2)} & = & \frac{\F[x,y]}{(x^2,y^3)}\\
	T &= & \frac{Q}{\Ann(XY)} & = & \frac{\F[x,y]}{(x^2,y^2)}\\
	\end{alignat*}
	with the natural projections $\pi_A,\pi_B\colon A,B\rightarrow T$ and Thom classes $\tau_A=x$ and $\tau_B=y$.  The fibered product is 
	$$A\times_TB=\F[(x,x),(y,y)]=\frac{\F[x,y]}{(x^3,y^3,x^2y^2)}=\frac{Q}{\Ann(X^2Y,XY^2)}.$$
	But now we see that $\pi_A(\tau_A)=x\neq y=\pi_B(\tau_B)$, so we do not have a total Thom class!  Therefore it follows that for this particular choice of $A$, $B$, and $T$ (and projection maps) there can be no connected sum!  Note however that in this case, if we set
	$$C=\frac{\F[x,y]}{\Ann(X^2Y-XY^2)}=\frac{\F[x,y]}{(x^3,y^3,x^2+y^2+xy,x^2y+xy^2)}$$
	then we do get the expected Hilbert function $H(C)=(1,2,2,1)$, i.e.  
	\begin{align*}
	H(C)= & H(A)+H(B)-H(T)-H(T)[1]\\
	= & (1,2,2,1)+(1,2,2,1)-(1,1,1,0)-(0,1,1,1)=(1,2,2,1).
	\end{align*}
	On the other hand, we claim that $C$ \emph{is} a connected sum for some \emph{different choices} $F^\prime, G^\prime$.  We will need them to satisfy
	$$F'-G'=X^2Y-XY^2=X(XY-Y^2)=\frac{1}{4}\left[X(X^2-(X-2Y)^2)\right].$$
	Hence after the change of coordinates on $R$
	$$\begin{cases} X\mapsto & X\\
					Y\mapsto & U=(X-2Y)\\
					\end{cases}$$
	we can take ${F'}=X^3, {G'}=XU^2$, so that $M_0'=\operatorname{gcd}(F',G')=X$.  With these choices we have $F'\nmid \left(M_0'\right)^2$ and $G'\nmid \left(M_0'\right)^2$ hence condition (2.) of Proposition \ref{prop:TFAE} is satisfied.
	Setting 
	\begin{alignat*}{3}
	A' &= & \frac{Q}{\Ann(X^3)} &= & \frac{\F[x,u]}{(x^4,u)}\\
	B' &= & \frac{Q}{\Ann(XU^2)} &= & \frac{\F[x,u]}{(x^2,u^3)}\\
	T &'= & \frac{Q}{\Ann(X)} &= & \frac{\F[x,u]}{(x^2,y)}
	\end{alignat*}
	(with natural projections) we see that 
	$$C\cong A'\#_{T'}B',$$
	and 
	\begin{align*}
	H(C)= & H(A)+H(B)-H(T)-H(T)[2]\\
	= & (1,1,1,1)+(1,2,2,1)-(1,1,0,0)-(0,0,1,1)=(1,2,2,1).
	\end{align*}	
	Note that in this case, the total Thom class is $\tau'=x^2+u^2\cong x^2+y^2+xy$.			
\end{example}

\section{Lefschetz Properties.}
\label{sect:LP}

In this section, we study the Lefschetz properties for fibered products and connected sums of AG algebras.  The weak and strong Lefschetz properties  are algebraic properties of Artinian graded algebras modeled after the property of cohomology rings of complex projective varieties given by the Hard Lefschetz Theorem. This theorem yields that multiplication by powers of a class of a hyperplane is a vector space isomorphism between a graded component of the cohomology ring and the component of complementary degree. In algebraic terms, this can be translated to the following:

\begin{definition}\label{def:Lefschetz}
Let $A$ be any graded Artinian algebra.
A linear form $\ell\in A_1$ is said to be a \emph{weak Lefschetz element} if the multiplication maps $\times \ell\colon A_i\rightarrow A_{i+1}$ have maximal rank for each degree $i$.  It is a \emph{strong Lefschetz element} if the multiplication maps $\times\ell^k\colon A_i\rightarrow A_{i+k}$ have maximal rank for every degree $i$ and every exponent $k$.

The algebra $A$ is said to have the \emph{weak Lefschetz property} or WLP (respectively the \emph{strong Lefschetz property} or SLP) if it has a weak (respectively strong) Lefschetz element. 	
\end{definition}
	  \begin{remark}
  	\label{rem:SLP}
 
  		If $A$ is a graded Artinian algebra with socle degree $d$ with a symmetric Hilbert function (for example, if $A$ is an AG algebra), then an element $\ell\in A_1$ is strong Lefschetz if and only if the multiplication maps  $\times\ell^{d-2i}\colon A_i\rightarrow  A_{d-i}$ are bijections for $0\leq i\leq \left\lfloor\frac{d}{2}\right\rfloor$.

 A graded Artinian algebra $A$ with socle degree $d$
 is said to have the \emph{strong Lefschetz property in the narrow sense} if the multiplication maps  $\times\ell^{d-2i}\colon A_i\rightarrow  A_{d-i}$ are bijections for $0\leq i\leq \left\lfloor\frac{d}{2}\right\rfloor$ (see \cite[p.143]{HMMNWW}). Note that having this property implies that the algebra has a symmetric Hilbert function.
\end{remark}

It is well known that the strong Lefschetz properties behave differently over fields of characteristic zero and fields of positive characteristic, e.g. \cite{Cook},\cite{CN}.  To avoid these complications, we shall often  assume that our ground field $\F$ has characteristic zero, or is an infinite field of characteristic greater than $d$, the socle degree of $A$.

Let $Q$ be a (not necessarily standard-)graded polynomial ring and write $Q=\F[x_1,\ldots,x_n]$ so that $x_1,\ldots,x_r$ $(r\leq n)$ are the algebra generators of degree one; in particular, $Q$ has the standard grading if and only if $r=n$.  Then for any polynomial $F=F(X_1,\ldots,X_n)\in R=\F[X_1,\ldots,X_n]$ (its divided power algebra described in section \ref{sect:MDG}), we may regard $F$ as a polynomial in the coordinates of the $n$-dimensional vector space $\mathfrak{m}_R/\mathfrak{m}_R^2=\operatorname{span}_{\F}\left\{x_1,\ldots,x_n\right\}$, the span of $\{x_1,\ldots,x_n\}$.  We shall write $F_{1}\in\F[X_1,\ldots,X_r]$ to mean the restriction of $F$ to the subspace of linear forms in $Q$, i.e. 
$$F_{1}(X_1,\ldots,X_r)=F(X_1,\ldots,X_r,0,\ldots,0).$$  

\begin{lemma}
	\label{lem:0Lef}
	If $A=Q/\Ann(F)$ is an AG algebra of socle degree $d$ over a field $\F$ of characteristic zero or over an infinite field of characteristic larger than $d$ and if $A_{1}$ is not identically zero, then there is a linear form $\ell\in A_1$ such that $\ell^d\neq 0$.  In particular, if $A$ has the standard grading, then there is always a linear form $\ell_A\in A_1$ such that $\ell_A^d\neq 0$. 
\end{lemma} 
\begin{proof}
	We have the following general formula 
	\begin{equation}
	\label{eq:Der}
\left(C_1x_1+\cdots+C_rx_r\right)^d\circ F_1(X_1,\ldots,X_r)=F(C_1,\ldots,C_r).
	\end{equation}
	Since $\rm{char}(F)=0$ or $\rm{char}(F)>d$, we have $\left(C_1x_1+\cdots+C_rx_r\right)^d\neq 0$. Thus, assuming that $F(X_1,\ldots,X_r)$ is not identically zero, there must be some linear form $\ell_A=C_1x_1+\cdots+C_rx_r\in Q_1$ for which $\ell_A^d\notin\Ann(F)$.
\end{proof}

\begin{remark}
	\begin{enumerate}[(a).]
		\item In general, Lemma \ref{lem:0Lef} can fail in small positive characteristic, even in the standard graded case, cf. \cite{Cook}. An example is $A=\Z_2 [x,y]/(x^2,y^2)$.
		\item Formula \eqref{eq:Der} was employed in \cite{MW} (see \cite[Theorem 3.76]{HMMNWW}) for standard graded algebras to derive a connection between Lefschetz properties and Hessians.
		\item Note that Lemma \ref{lem:0Lef} implies that the set 
		$\mathcal{D}_A=\left\{\ell_A\in A_1\left|\ell_A^d\neq 0\right.\right\}$
		is a nonempty, Zariski open set in $\Spec(\Sym(A_1))$.
	\end{enumerate}
\end{remark}

\begin{definition}\label{Leflocusdef}
In what follows we shall write $\L_A\subset A_1$ to denote the (possibly empty) Zariski open set in $\P(A_1)$ of strong Lefschetz elements for $A$. The set $\L_A$ is called the {\em Lefschetz locus} for the Artinian algebra $A$. 
\end{definition}

Note that Definition \ref{def:Lefschetz} can be reformulated to say that an Artinian algebra has the SLP or WLP respectively if and only if the corresponding Lefschetz locus is non-empty.

\subsection{Strong Lefschetz property of a connected sum.}

The following two results show that the classes of $\F$-algebras having the SLP and WLP are closed under taking fibered products and connected sums over the field $\F$. Their subclasses consisting of standard graded algebras are also closed under taking factors in either the fibered product or the connected sum.

\begin{proposition}
	\label{prop:SLPFP}
If $A$ and $B$ are AG algebras of the same socle degree that each have the SLP, then the fibered product $D=A\times_\F B$ over a field $\F$  also has the SLP.  If $A$ and $B$ have the standard grading, then the converse holds as well. 	
\end{proposition}
	\begin{proof}
		Assume that $A$ and $B$ have the SLP, let $\ell_A\in \L_A$ and $\ell_B\in\L_B$ be Lefschetz elements.  For $0< i\leq i+k\leq d$, $D_i=A_i\oplus B_i$, and $D_{i+k}=A_{i+k}\oplus B_{i+k}$. Furthermore  the multiplication maps $\ell_A^k:A_i\to A_{i+k}$ and $\ell_B^k:B_i\to B_{i+k}$ are either both injective or both surjective, hence the multiplication map 
		$$\times(\ell_A,\ell_B)^{k}\colon D_i\rightarrow D_{i+k}$$
is also injective or surjective respectively.  For $i=0$, we have $D_0=\F$. Fix $k>0$ and consider 
		$$\times(\ell_A,\ell_B)^{k}\colon D_{0}\rightarrow D_k.$$
		Since $A, B$ have the SLP, $\ell_A^k\neq 0$ and $\ell_B^k\neq 0$, hence also  $c\ell_A^k\neq 0$ and $c\ell_B^k\neq 0$ for any $c\in\F$. Since the map displayed above takes $c\mapsto(c\ell_A,c\ell_B)\neq 0$ for all $c\in \F$, we see that it is injective.
		
		Conversely, assume that $A$ and $B$ have the standard grading and suppose that $D$ has the SLP, and let $\mathcal{L}_D\subset D_1\cong A_1\times B_1$ be the set of strong Lefschetz elements.  Define the subset $\mathcal{D}_A\subseteq A_1$ by 
		$$\mathcal{D}_A=\left\{\ell_A\in A_1\left|\ell_A^d\neq 0\right.\right\}.$$
		and similarly for $\mathcal{D}_B\subseteq B_1$.  By Lemma \ref{lem:0Lef}, these are non-empty open sets, hence so is their product $\mathcal{D}_A\times\mathcal{D}_B\subset A_1\times B_1$.  This implies that the intersection $\mathcal{L}_D\cap\left(\mathcal{D}_A\times\mathcal{D}_B\right)$ is non-empty, hence there exists a strong Lefschetz element, $\ell=(\ell_A,\ell_B)\in D_1$ such that $\ell_A^d\neq 0\neq \ell_B^d$.  We will show that $\ell_A\in A_1$ is a strong Lefschetz element for $A$, and the argument for $\ell_B\in B_1$ is similar.  Fix $0\leq i\leq \frac{d}{2}$, and consider $\times\ell^{d-2i}\colon A_i\rightarrow A_{d-i}$.  Suppose that $\alpha\in A_i$ is in the kernel.  If $i>0$ then $(\alpha,0)\in D_i$ and this element is in the kernel of $\times\ell\colon D_i\rightarrow D_{d-i}$ which implies that $\alpha=0$.  If $i=0$, and $\alpha\neq 0$, then we can take $\alpha=1$, and deduce that $\ell_A^d=0$.  But this contradicts our choice of $\ell\in\mathcal{D}_A\times\mathcal{D}_B$.  This shows that $A$ and $B$ both have the SLP. 
\end{proof}

\begin{proposition}
	\label{prop:SLPCS}
	If $A$ and $B$ both have the SLP, then the connected sum $C=A\#_\F B$ over a field $\F$ also has the SLP.  If $A$ and $B$ have the standard grading, then the converse holds as well.
\end{proposition}	
	\begin{proof}
		Fix $\ell_A\in\L_A$ and $\ell_B\in\L_B$, and consider the multiplication map
		$$\times (\ell_A,\ell_B)^{d-2i}\colon C_i\rightarrow C_{d-i}.$$
		If $0<i<d$, then $C_i=(A\oplus B)_i$, hence the multiplication map is an isomorphism.  It remains to see the isomorphism for the case $i=0$:  but $(\ell_A,\ell_B)^d=(\ell_A^d,\ell_B^d)\equiv 0$ in $C$ if and only if $\ell_A=a\tau_A$ and $\ell_B=a\tau_B$ for the same $a\in\F$.  But if this occurs, we can replace $\ell_B$ by $b\ell_B$ for some $b\in \F$ such that $b^d\neq 1$.  Then the element $(\ell_A,b\ell_B)$ will be strong Lefschetz for $C$. 
	
		Conversely, assume that $A$ and $B$ have the standard grading and suppose that $C$ has the SLP, and let $\mathcal{L}_C\subset C_1$ be the set of strong Lefschetz elements.  As before, define the subsets $\mathcal{D}_A\subset A_1$ and $\mathcal{D}_B\subset B_1$ as the linear forms whose $d^{th}$ power is non-vanishing.  Then Lemma~\ref{lem:0Lef} implies that $\mathcal{D}_A$, $\mathcal{D}_B$ are non-empty open sets in $A_1$, $B_1$, respectively.  Therefore the non-empty Zariski open set $\L_C\subset C_1\cong A_1\times B_1$ must intersect the non-empty open set $\mathcal{D}_A\times\mathcal{D}_B\subset A_1\times B_1$, hence there exists $\ell=(\ell_A,\ell_B)\in\L_C$ such that $\ell_A^d\neq 0$ and $\ell_B^d\neq 0$.  We claim that $\ell_A\in A_1$ is strong Lefschetz for $A$ and similarly for $B$.  Indeed, fix $0\leq i\leq \frac{d}{2}$ and consider the multiplication map
		$$\times\ell_A^{d-2i}\colon A_i\rightarrow A_{d-i}.$$
		Suppose that $\alpha\in A_i$ is in the kernel.  If $i>0$, then $(\alpha,0)\in C_i$ is also in the kernel of $\times\ell^{d-2i}\colon C_i\rightarrow C_{d-i}$, which implies that $\alpha=0$.  If $i=0$, and $\alpha\neq 0$, then we can assume that $\alpha=1$, and we must have $\ell_A^d=0$, which was ruled out by our choice of $\ell$.  Thus, both $A$ and $B$ must also have the SLP. 
\end{proof}

\begin{remark}
\label{rem:Hessian}
If $T=\F$ and $A$ and $B$ are standard graded, then the connected sum $A\#_TB$ is also standard graded, by Proposition \ref{prop:ConnSumF}.  Hence in this case, one can also use Theorem~\ref{theorem1} regarding the Macaulay dual generator of a connected sum and the theory of Hessians to establish Proposition \ref{prop:SLPCS}.  See \cite[Theorem 3.76 and Proposition 3.77(ii)]{HMMNWW}.
\end{remark}

\begin{remark}
\label{rem:LefschetzLocusProd}
From the proof of Proposition \ref{prop:SLPFP}, we have the inclusion
 of Lefschetz loci $\L_A\times\L_B\subseteq\L_{A\times_\F B}$.  The following example shows that in general $\L_A\times\L_B\not\subseteq\L_{A\#_\F B}$.
\end{remark}

\begin{example}
	Let $A= \F[x]/(x^2)$, $B= \F[y]/(y^2)$, and $T=\F$
	with usual orientations and projection maps $\pi_A\colon A\rightarrow T$ and $\pi_B\colon B\rightarrow T$, so that the Thom classes are the socle generators $\tau_A=x$ and $\tau_B=y$.  Then the fibered product is given by  

	$$D=\frac{\F[z_1,z_2]}{(z_1^2,z_2^2,z_1z_2)}, \ \begin{cases} z_1= & (x,0)\\ z_2= & (0,y)\\ \end{cases}$$
	which has the strong Lefschetz element $\ell=(x,y)\cong z_1+z_2$.  But we have $(\tau_A,\tau_B)\cong z_1+z_2$, hence the connected sum is given by 
	$$C=\frac{\F[z_1,z_2]}{(z_1^2,z_2^2,z_1z_2,z_1+z_2)},$$ 
	which means that $\ell=(x,y)\cong z_1+z_2$ is zero, hence is not strong Lefschetz in $C_1$.  Note however that $\ell'=(x,2y)\in C_1$ is a strong Lefschetz element. 
\end{example}

The next example shows that we should not expect the converse of either Proposition \ref{prop:SLPFP} or Proposition \ref{prop:SLPCS} to hold in the non-standard graded case even when $T=\F$.
\begin{example}
	\label{ex:NonStd}
	Let $A= \F[x]/(x^2)$ with $\deg(x)=2$ and orientation $\int_A\colon x\mapsto 1$ (hence socle degree $d=2$), and let $B=\F[y]/(y^3)$ with the standard grading and the standard orientation $\int_B\colon y^2\mapsto 1$.  Finally set $T=\F$ with the usual projection maps $\pi_A\colon A\rightarrow T$ and $\pi_B\colon B\rightarrow T$, so that the Thom classes are $\tau_A=x$ and $\tau_B=y^2$.  Then the fibered product $D$ and connected sum $C$ are given by  
	\begin{align*}
	D= & \frac{\F[z_1,z_2]}{(z_1^2,z_2^3,z_1z_2)}, & \begin{cases} z_1= & (x,0)\\ z_2= & (0,y)\\ \end{cases}\\
	C= & \frac{\F[z_1,z_2]}{(z_1^2,z_2^3,z_1z_2,z_1+z_2^2)}, & \tau=(\tau_A,\tau_B)=z_1+z_2^2 . 
	\end{align*}
Then note that both $D$ and $C$ satisfy the SLP (with strong Lefschetz element $\ell=z_2$), but $A$ does not have SLP (it has no linear elements).
\end{example}

We end the discussion on the SLP by providing another class of algebras which are connected sums and satisfy this property. In geometry, the connected sum of an $n$-dimensional manifold with $\P^n$ is diffeomorphic as an oriented manifold to the blow up of $M$ at a point \cite[p.~101]{Huybrechts} and blowing up a smooth projective variety along a smooth projective subvariety preserves projectivity \cite[Proposition 7.16]{Ha}.  This is one reason one might expect $A\# _{\F} B$ to have SLP  when $A$ has SLP and $B=\F[x] / (x^n) $. More generally, we consider below the class of rings of the form $A\# _{T}T[x] / (x^n) $  and we show that they satisfy the strong Lefschetz property.

\begin{theorem}
\label{thm:blowup}
Let $A, T$ be AG algebras with socle degrees $d, k$ respectively and let $\pi_A\colon A\rightarrow T$ be a surjective ring homomorphism such that its Thom class $\tau_A$ satisfies $\pi_A(\tau_A)=0$.  Let $x$ be an indeterminate of degree one, set $B=T[x]/(x^{d-k+1})$, and define $\pi_B\colon B\rightarrow T$ to be the natural projection map satisfying $\pi_B(t)=t$ and $\pi_B(x)=0$.  In this setup, if $A$ and $T$ both satisfy the SLP,  then the fibered product $A\times_T B$ also satisfies the SLP.  Moreover if the field $\F$ is algebraically closed, then the connected sum $A\#_T B$  also satisfies the SLP.
\end{theorem}

\begin{proof}
Let $D=A\times_TB$ be the fibered product.  We have a graded $A$-module decomposition
	\begin{equation}
	\label{eq:Ddecomp}
	D\cong A\oplus \underbrace{Tx\oplus\cdots\oplus Tx^{d-k}}_J
	\end{equation} 
	where we identify $(a,b)\cong (a,\pi_A(a))\oplus (0,\sum_{i=1}^{d-k}t_ix^i)$, with $b=\pi_A(a)+\sum_{i=1}^{d-k}t_ix^i\in B=T[x]/(x^{d-k+1})$.  Since $J\subset D$ is an ideal in $D$, it follows that for any $\ell=(\ell_A,\ell_B)\in D$, the multiplication map $\times\ell\colon D\rightarrow D$ can be represented by a block matrix of the form 
	\begin{equation}
	\label{eq:Dblock}
	\left(\begin{array}{cc}
	\times\ell_A & 0\\ * & \times\ell_B|_J\\ \end{array}\right).
	\end{equation} 
	In particular, we see that the rank of the multiplication map $\times\ell\colon D\rightarrow D$ is the sum of the ranks of the multiplication maps $\times\ell_A\colon A\rightarrow A$ and $\times\ell_B\colon J\rightarrow J$.  Since $A$ and $T$ have SLP, and since $\pi_A\colon A\rightarrow T$ is surjective, we may choose a strong Lefschetz element $\ell_A\in A_1$ for $A$ such that $\pi_A(\ell_A)=\ell_T\in T^1$ is strong Lefschetz for $T$.  Then by a standard argument $\ell_B=\ell_T+t\in B_1$ is also strong Lefschetz for $B$.  Set $\ell=(\ell_A,\ell_B)$ as above.  We want to show that for each integer $m$ and for each degree $i$, the multiplication maps $\times\ell^m\colon D_i\rightarrow D_{i+m}$ have maximal rank.  According to the above discussion, it suffices to see that each of the multiplication maps $\times\ell_A^m\colon A_i\rightarrow A_{i+m}$ and $\times \ell_B^m\colon J_i\rightarrow J_{i+m}$ have maximal rank.   Since $\ell_A$ is strong Lefschetz for $A$, the map $\times\ell_A^m\colon A_i\rightarrow A_{i+m}$ has maximal rank.  Since $\ell_B$ is strong Lefschetz for $B$, the map $\times\ell_B^m\colon B_i\rightarrow B_{i+m}$ has maximal rank.  Also, since $B\cong T\oplus J$ (as $T$-modules) and $J$ is an ideal in $B$, it follows that $\times\ell_B^m$ decomposes as in \eqref{eq:Dblock} as 
	$$\left(\begin{array}{cc} \times\ell_T^m & 0\\ * & \times \ell_B^m|_J\\ \end{array}\right).$$
	Finally since $\ell_T$ is strong Lefschetz for $T$ we know $\times\ell_T^m$ has maximal rank too, and thus so does $\times\ell_B^m\colon J_i\rightarrow J_{i+m}$.  This shows that $D$ has SLP.

	Note that the Thom class for $\pi_B\colon B\rightarrow T$ is $\tau_B=x^{d-k}$, and by the assumptions above we have $\pi_A(\tau_A)=0=\pi_B(\tau_B)$.  Therefore $(\tau_A,\tau_B)=(\tau_A,x^{d-k})=(\tau_A,0)+(0,x^{d-k})\in D=A\times_TB$ and we can form the connected sum; let $C=A\#_TB=A\times_TB/(\tau_A,\tau_B)$ be that connected sum.  We have, similar to \eqref{eq:Ddecomp}, an $A$-module decomposition 
\begin{equation}
\label{eq:Cdecomp}
C\cong A\oplus \underbrace{Tx\oplus\cdots\oplus Tx^{d-k-1}}_I
\end{equation}
	Note, there is a missing factor from $D$ because in $C$ the element $(0,tx^{d-k})\in Tx^{d-k}$ is identified with $(a\tau_A,0)\in A$, where $a\in A$ is any element with $\pi_A(a)=t$.
	Also, the summand $I\subset C$ is not an ideal, and therefore we cannot decompose multiplication maps as in \eqref{eq:Dblock}.  On the other hand, we can find a flat family that deforms $C$ to an $A$-module whose multiplication maps do decompose as in \eqref{eq:Dblock}.

	Specifically, we define the parameter algebra 
	\begin{equation}
	\label{eq:1par}
	\mathcal{R}=\frac{\left(A\times_TB\right)[z]}{\left(z(\tau_A,0)+(0,\tau_B=x^{d-k})\right)}.
	\end{equation}
	Clearly $\mathcal{R}$ is a finite ring extension over $\F[z]$.  It is also clear that for each $c\in\F$, the fiber 
	\begin{equation}
	\label{eq:fiber}
	\mathcal{R}_c\coloneqq \frac{\mathcal{R}[z]}{(z-c)\cdot\mathcal{R}}\cong\frac{A\times_TB}{(c\cdot(\tau_A,0)+(0,\tau_B))}
	\end{equation}
	decomposes as in \eqref{eq:Cdecomp}, hence all the fibers have the same length. Since the field $\F$ is algebraically closed, this implies that $\mathcal{R}$ is flat over $\F[z]$.\footnote{This criterion for flatness is true in much more generality (see \cite[Exercise II.5.8]{Ha}); in this simple case it is easy to see using the decomposition of $\mathcal{R}$ as a module over $\F[z]$ (a PID).}
	We now show\vskip 0.2cm\par\noindent
	{\bf Claim}. For every $0\neq c\in\F$ there is an $\F$-algebra isomorphism \begin{equation}
	\label{eq:isom}
	\psi_c\colon\mathcal{R}_c\rightarrow\mathcal{R}_1=A\#_TB.
	\end{equation}
	To see this claim note first that for any graded Artinian $\F$-algebra $A$ and for each $t\in\F$ there is an $\F$-algebra homomorphism $\phi_{A,t}\colon A\rightarrow A$ defined by $\phi_{A,t}(a)=t^{\deg(a)}\cdot a$.  Also, given any graded Artinian $\F$-algebra $T$, and setting $B=T[x]/(x^{d-k+1})$ there is an $\F$-algebra homomorphism $\phi_{B,t}\colon B\rightarrow B$ defined by 
	$$\phi_{B,t}(b=b_0+b_1x+\cdots+b_{d-k}x^{d-k})=t^{\deg(b_0)}\cdot b_0+\cdots+t^{\deg(b_{d-k})}\cdot b_{d-k}x^{d-k}.$$
	(Note that in order to be an $\F$-algebra map, we need $x^{d-k+1}\equiv 0$ in $B$).
	Then for each $t\in\F$ the product map gives an $\F$-algebra map on the fibered product
	$$(\phi_{A,t},\phi_{B,t})\colon A\times_TB\rightarrow A\times_TB.$$
	Moreover since $\phi_{A,t},\phi_{B,t}$ are isomorphisms for $t\neq 0$, so is their product map $(\phi_{A,t},\phi_{B,t})$, which passes to an $\F$-algebra isomorphism on the fibers
	$$(\phi_{A,t},\phi_{B,t})\colon \frac{A\times_TB}{(c\cdot(\tau_A,0)+(0,x^{d-k}))}\rightarrow\frac{A\times_TB}{(c\ t^{\deg(\tau_A)}\cdot(\tau_A,0)+(0,x^{d-k}))}$$
	Then taking $t=c^{-1/(d-k)}\in\F$ (which exists since $\F$ is algebraically closed) yields the claim.\footnote{In deformation theory, $\mathcal{R}_c$ is called a \emph{jump deformation} of $\mathcal{R}_0$; see \cite[\S 6]{Fox}.}\vskip 0.2cm

	Here the special fiber is 
	$$\mathcal{R}/z\cdot \mathcal{R}\cong \frac{A\times_TB}{(0,x^{d-k})}=C'$$
	and in terms of the decomposition of \eqref{eq:Cdecomp} the summand $I\subset C'$ is an ideal.  Therefore, taking $\ell=(\ell_A,\ell_B=\ell_T+x)\in C'_1$ as before, the multiplication maps $\times\ell^m\colon C'_i\rightarrow C'_{i+m}$ decompose as in \eqref{eq:Dblock}.  Hence the special fiber $\mathcal{R}/z\cdot\mathcal{R}$ has SLP.  Since SLP is a determinantal condition, and since $\mathcal{R}$ is flat, it follows that there exists some $0\neq c\in\F$ for which the fiber $\mathcal{R}/(z-c)\mathcal{R}$ also has SLP, and from the isomorphism in \eqref{eq:isom} we conclude that the connected sum $A\#_TB=\mathcal{R}_1\cong \mathcal{R}_c$ also has SLP. 
\end{proof}
In a forthcoming paper, the authors together with P. Macias Marques and J. Watanabe will investigate further this interesting subclass of connected sums called cohomological blow ups \cite{I-W}.
%In a sequel the authors with P. Macias Marques and J. Watanabe consider a subclass of connected sums, the blow up graded algebras \cite{I-W}.

We next give an application of Theorem \ref{thm:blowup} to questions pertaining to the strong Lefschetz property of codimension three AG algebras. Despite known structure theorems for this class of algebras, little is known about their Lefschetz properties, even more so regarding the SLP.  The most prominent result is that  in characteristic zero all standard graded complete intersection Artinian algebras of codimension three have the WLP \cite{HMNW} and it is an open problem whether in fact all codimension three standard graded AG algebras have the WLP in characteristic zero. For the case of codimension three AG algebras that are not necessarily complete intersections, it is known that for each possible Hilbert function an AG example exists having the WLP \cite{Harima}.\footnote{T. Harima in \cite{Harima} establishes more generally the existence of a weak Lefschetz AG example for any symmetric Hilbert function of socle degree $j$ satisfying the SI condition: the first difference $\Delta H_{\le j/2}$ is an O-sequence -- that is, $\Delta H_{\le j/2}$ occurs as the Hilbert function of some Artinian 
algebra.} The issue of the weak Lefschetz property of all AG codimension three algebras is reduced to determining the property for compressed AG algebras (maximum Hilbert function given the socle degree) in \cite{BMMNZ}; they also show that all codimension three AG algebras of socle degree five or less are strong Lefschetz. However, the existence of an AG algebra of a given non-CI Hilbert function that satisfies the strong Lefschetz property is open in general, even in codimension three. As an application of Theorem \ref{thm:blowup} we next give an infinite family of non-CI Hilbert functions of codimension three AG algebras, for each of which there is an AG algebra having SLP.  We also show a closure property for the set $\mathcal {SL}(d)$ of Artinian Gorenstein Hilbert functions admitting an algebra that satisfies the SLP (Corollary \ref{closurecor}).

\begin{corollary}\label{heightthreecor}[Codimension three AG Hilbert functions having SLP]
Let $\F$ be an algebraically closed field. For each choice of positive triple $(a,d,k) $ such that $a\leq d-a$ and $ 2k<d$,  the class of codimension three AG $\F$-algebras having one of the Hilbert functions displayed below contains at least one member which satisfies the SLP:\\
Case 1: assume $k\leq a\leq b=d-a$, and set\small
\begin{equation}\label{firstHFeqn}
H=(\underbrace{1,3,5, \dots, 2k+1}_{\Delta H= 2}, \underbrace{2k+2, \dots, k+a+1}_{\Delta H=1}, \underbrace{k+a+2, \dots, k+a+2}_{ b-a+1 \text{ times }},\underbrace{\dots 2k+2, 2k+1, \dots, 5, 3, 1}_{\text{symmetric to the first part}} ).
\end{equation}
(If $k<a$ then the second subsequence is $2k+3, 2k+4,\ldots, k+a+1$, if there is room.)\vskip 0.2cm\par\noindent
Case 2: assume $a\leq k< b=d-a$ and set
\begin{equation}\label{secondHFeqn}
H=(\underbrace{1,3,5, \dots, k+a+1}_{\Delta H =2}, \underbrace{k+a+2, \dots, k+a+2}_{ d-2k+1 \text{ times }}, \underbrace{k+a+1, \dots, 5, 3, 1}_{\Delta H= -2}).
\end{equation}
In particular, the graded Artinian Gorenstein algebra 
\begin{equation}\label{CSLeq}
C=\F [s,x,y]/(xs,s^{a+1},x^{d-k+1},y^{b+1},s^ay^{b-k}-x^{d-k})
\end{equation} 
with Macaulay dual generator $f=S^aY^b-X^{d-k}Y^k$, where $a,b,k$ satisfies one of the conditions Case~1 or Case 2, has the respective Hilbert function, and is strong Lefschetz. 
\end{corollary}
\normalsize
\begin{proof}
Fix integers $a,d,k$ as in the statement of the claim, let $b=d-a$ and note that $2b \geq a+b=d>2k$ yields $b>k$.
Consider the rings 
\begin{equation} A=\F[s,y]/(s^{a+1},y^{b+1}), T=\F[y]/(y^{k+1}), B=T[x]/(x^{d-k+1})=\F[x,y](x^{d-k+1}, y^{k+1})
\end{equation} 
with Hilbert functions
\begin{eqnarray*}
H(A) &=(1, 2, 3 \dots, a, \underbrace{a+1, \dots, a+1}_{b-a+1}, a, \dots, 3, 2, 1) \\
H(B) &=(1, 2, 3 \dots, k, \underbrace{k+1, \dots, k+1}_{d-2k+1}, k, \dots, 3, 2, 1) \\
H(T) &=(\underbrace{1, 1, \dots, 1}_{k})
\end{eqnarray*}
Set $\pi_A:A\to T$ to map $s\mapsto 0, y\mapsto y$ and let $\pi_B:B\to T$ map $x\mapsto 0, y\mapsto y$. The connected sum $C=A \#_TB$ satisfies the SLP by Theorem \ref{thm:blowup}; by Lemma \ref{lem:HFCS}  $H(C)$ is equal to one of the displayed Hilbert functions \eqref{firstHFeqn} or \eqref{secondHFeqn} depending on whether $k\leq a$ or $a\leq k$. The Macaulay dual generator for $A$ is $S^aY^b$, for $B$ is $X^{d-k}Y^k$. By Theorem \ref{theorem1} the Macaulay dual generator for $C$ is $f=S^aY^b-X^{d-k}Y^k$.
\end{proof}
We denote by $\mathcal {SL}(d)$ the family of symmetric Gorenstein sequences $H$ having socle degree $d$, such that there is a standard-graded Artinian Gorenstein algebra $A$ of Hilbert function $H$ having the strong Lefschetz property.
Recall that the subscheme PGor$(H)\subset \mathbb P(R_d)$, the projectivization of $R_d$,  parametrizes the graded Artinian Gorenstein quotients of $R$ having Hilbert function $H$, via the Macaulay dual generator (Fact \ref{fact:MDElts}).\footnote{For a discussion of the subtleties of parametrizing PGor$(H)$ see \cite[Theorem 44 ff]{Kl}.}
 \begin{corollary}\label{closurecor}[Closure]  Let $2k<d$, assume that $\F$ is algebraically closed of characteristic $\cha \F=0$ or $\cha \F>d$,  and let  \par
$W(k,d)=(0_0,1,2,\ldots ,k, \underbrace{k+1, \dots, k+1}_{d-2k-1}, k, \ldots, 3, 2, 1_{d-1},0_d) $. Then\par
A. The set $\mathcal {SL}(d)$ includes all CI Hilbert functions of socle degree $d$ and codimension two.\par
B. the set $\mathcal {SL}(d)$ is invariant under the addition of any such sequence $W(k,d)$.  \par
C. Suppose $H\in \mathcal {SL}(d)$. Then there is an irreducible component of PGor$(H)$ whose generic element has SLP. In particular, if $H$ is of codimension three ($H(1)=3$) with $H\in \mathcal {SL}(d)$, then the generic element of $PGor(H)$ is SL.\par
\end{corollary}
\begin{proof} The statement (A) is well-known (\cite[Theorem 2.9]{Ia}).  The statement (C) follows from (B) as the strong Lefschetz property is an open condition, since the maps involved are semicontinuous for  a fixed Hilbert function (\cite[p. 102]{HMNW}, \cite[\S 2.7]{IMM}); also, the irreducibility of PGor$(H)$ in codimension three is well known \cite{Di}. To show (B), let $A$ be a codimension $c$ strong Lefschetz Artinian Gorenstein algebra of Hilbert function $H\in \mathcal {SL}(d)$, and let $\ell\in A_1$ be a strong Lefschetz element for $A$. Take a new basis of $A_1$, $\ell_1=x_1, x_2,\ldots, x_c$, let $T=\F[y]/(y^{k+1}),  B=T[x]/(x^{d-k+1})=\F[x,y](x^{d-k+1}, y^{k+1})$. Define a homomorphism $\pi_A: A\to T$ by $\pi_A (x_1)=y$ and $\pi_A(x_i)=0$ for $i\neq 1$; and $\pi_B (x)=0, \pi_B(y)=y$. By Theorem \ref{theorem1} $C=A\#_T B$ is strong Lefschetz, and by Lemma  \ref{lem:HFCS} it has the Hilbert function $H(A)+W(k,d)$. This completes the proof.
\end{proof}
For example, by Corollaries \ref{heightthreecor} and \ref{closurecor} the sequence $H=(1,3,5,7,7,5,3,1)+W(3,7)=(1,3,5,7,7,5,4,1)+(0,1,2,3,3,2,1,0)=(1,4,7,10,10,7,4,1)$ is in $\mathcal {SL}(7)$.
	
\vskip 0.3cm

In the proof of Theorem \ref{thm:blowup}, we used the $\F$-algebraically closed condition to 
\begin{enumerate}[(a).]
	\item get a flatness criterion:  If $\mathcal{R}$ is a finite ring extension over $\F[z]$ with equidimensional fibers then $\mathcal{R}$ is flat over $\F[z]$, and
	
	\item to show that the fibers of the flat extension $\mathcal{R}_c$ all had the same isomorphism type for all $0\neq c\in\F$.
\end{enumerate}

The following examples show the troubles that can occur over a non-algebraically closed field (e.g. $\Q$).   

\begin{example}
	\label{ex:bad}
	Let $A=\Q[x,y]/(x^2,y^2)$, and set 
	$$\mathcal{R}=\frac{\Q[x,y][z]}{((z^2+1)x^2,y^2)}.$$
	Then for every $c\in\Q$ we have 
	$$\mathcal{R}/(z-c)\mathcal{R}\cong \frac{\Q[x,y]}{((c^2+1)x^2,y^2)}\cong A$$
	hence $\mathcal{R}$ has all fibers of the same length, but is not $\Q[z]$-free.
\end{example}

The next remarkable example is taken from \cite[Example 3.1]{AAM}.
\begin{example}
	Let $A=\Q[x]/(x^3)$, $B=\Q[y]/(y^3)$, and $T=\Q$ with $\pi_A\colon A\rightarrow T$ and $\pi_B\colon B\rightarrow T$ the natural projections and Thom classes $\tau_A=x^2$ and $\tau_B=-y^2$ (give $B$ the opposite orientation).  Then the fibered product is $A\times_TB$ generated as a $\Q$-algebra by $X=(x,0)$ and $Y=(0,y)$.  Therefore we see that the fibered product is 
	$$A\times_TB=\frac{\Q[X,Y]}{(X^3,Y^3,XY)}$$
	and the connected sum is 
	$$A\#_TB=\frac{\Q[X,Y]}{(X^2-Y^2,XY)}.$$
	In \cite{AAM}, it is observed that there is no $\Q$-algebra isomomorphism
	\begin{equation}
	\label{eq:Qiso}
	\psi_p\colon\mathcal{R}_p=\frac{\Q[X,Y]}{(X^2-pY^2,XY)}\rightarrow\frac{\Q[X,Y]}{(X^2-Y^2,XY)}=\mathcal{R}_1
	\end{equation}
	for $p\in\Z\subset\Q$ a prime not congruent to $3$ mod $4$.  We slightly extend this:
	
	{\em If $p\in\Q$ is any rational number that does not have a square root in $\Q$ (e.g. if $p$ is prime) then there can be no $\Q$-algebra isomorphism as in \eqref{eq:Qiso}.}
	
	Indeed, assume that there is a $\Q$-algebra homomorphism as in \eqref{eq:Qiso}.  Then for some $a,b,c,d\in\Q$ we must have $\psi_p(X)=aX+bY$ and $\psi_p(Y)=cX+dY$, 
	and in order for $\psi_p$ to be a $\Q$-algebra homomorphism these parameters must satisfy the following equations 
	\begin{align}
	\label{eq:1}
	(a^2+b^2)-p(c^2+d^2)= & 0\\
	\label{eq:2}
	ac+bd= & 0.
	\end{align}
	
	Moreover, in order for $\psi_p$ to be a $\Q$-algebra isomorphism, it is necessary that the determinant $ad-bc$ is non-zero.  There are two cases to consider:
	
	\textbf{Case 1:}  One of the variables, say $a=0$.  In this case, \eqref{eq:2} implies that $b$ or $d$ equals zero and \eqref{eq:1} implies it must be $d=0$, in which case we have 
	$b^2=pc^2$.
	This implies $p$ is a perfect square.
	
	\textbf{Case 2:} None of the variables are zero.  In this case \eqref{eq:2} gives
	$\frac{a}{b}=-\frac{d}{c}$
	and \eqref{eq:1} gives
	\begin{align*}
	b^2\cdot \left(\left(\frac{a}{b}\right)^2+1\right)=  pc^2\left(1+\left(\frac{d}{c}\right)^2\right)
	\Rightarrow 
	b^2=  p\cdot c^2
	\end{align*} 
	which again implies $p$ must be a perfect square.  
	
	Of course, if we replace $\Q$ with an algebraically closed field, this issue disappears.
\end{example}
\subsection{Weak Lefschetz property of connected sums.}

The aim of the next results is to connect the WLP for fibered products over $\F$ to the WLP for fibered products over arbitrary graded Artinian algebras $T$. 
	
\begin{lemma}\label{lem:injFiberProd}
For connected, graded $\F$-algebras $A,B$ and $T$ and $\F$-algebra homomorphisms $\pi_A:A\to T$ and  $\pi_B:B\to T$ there is an inclusion of graded $F$-algebras $A\times_TB \hookrightarrow A\times_\F B$.
\end{lemma}	
\begin{proof}
In order to prove this we first discuss the compatibility of the maps involved. We say that the fibered products $A\times_\F B$ and $A\times_T B$ or the connected sums $A\#_\F B$ and $A\#_T B$ are compatible if they arise from maps $\pi_A,\pi_B$ and $\pi'_A,\pi'_B$, respectively, which make the following diagram commute, where 
 the map $T\to\F$ is defined to be the canonical projection $T\to T/T_{>0}=\F$

\begin{equation*}
	\label{eq:compatible}
	\xymatrix{& A\ar[d]^-{\pi'_A} \ar[ddrr]^-{\pi_A}\\  B\ar[r]^{\pi'_B} \ar[drrr]_-{\pi_B} & T \ar[drr]\\ & & & \F}
	\end{equation*}

Note that, because the maps above are graded, the restriction of each of the maps $\pi_A,\pi_B$ and $\pi'_A,\pi'_B$ to the positive degree component of its domain followed by projection to $\F$ is the zero map. This observation yields the commutativity of the above diagram,

The containment $A\times_TB \hookrightarrow A\times_\F B$ now follows from the definition of the fibered product after noticing that commutativity of the above diagram implies that if  $\pi'_A(a)=\pi'_B(b)$ then $\pi_A(a)=\pi_B(b)$.
\end{proof}

In order to be able to transfer the Lefschetz properties from $A\times_\F B$ to $A\times_T B$ we must understand whether the Lefschetz elements of $A\times_\F B$ remain available in $A\times_T B$. Since in general $(A\times_T B)_1\subsetneq (A\times_\F B)_1$, this is a delicate, but manageable, task.

\begin{lemma}
\label{lem:survival}
Assume that $A$ and $B$ are AG algebras having either the WLP or the SLP and their respective Lefschetz loci are $\L_A$ and $\L_B$. Assume further that $T$ is an AG algebra such that the maps $\pi_A:A\to T$ and $\pi_B:B\to T$ used to define the fibered product $A\times_T B$  are surjections. Then there is a common element $(\ell, \ell^\prime)$ of $\L_A\times\L_B$ and $(A\times_T B)_1$, upon identifying the latter set with a subset of $(A\times_\F B)_1=A_1\times B_1$. That is, 
\begin{equation}\label{commonelementeq}
\exists\  (\ell,\ell')\in(\L_A\times\L_B)\cap(A\times_T B)_1.
\end{equation}
\end{lemma}
\begin{proof}
Denote by $K_A$ and $K_B$ the kernels of the $\F$-linear homomorphisms $\pi_A|_{A_1}$ and $\pi_B|_{B_1}$ and note that the isomorphisms $A_1/K_A\cong T_1 \cong B_1/K_B$ induce 
 corresponding isomorphisms 
 $$A_1\cong K_A \times A_1/(K_A)_1\cong K_A\times T_1 $$
 $$B_1\cong K_B \times B_1/(K_B)_1\cong K_B \times T_1.$$
Since the projection maps $\P(A_1)\to\P(T_1)$ and res\-pectively $\P(B_1)\to\P(T_1)$ are open (see \cite[\href{http://stacks.math.columbia.edu/tag/037G}{Lemma 037G}]{stacks}) it follows that  $\pi_A(\L_A), \pi_B(\L_B)$ are Zariski open sets of $\P(T_1)$. Since these sets must have a nonempty intersection, there exist $\ell\in\L_A$ and $\ell'\in\L_B$ such that $\pi_A(\ell)= \pi_B(\ell')$ yielding the desired element $(\ell,\ell')$ satisfying Equation \eqref{commonelementeq}.
\end{proof}

For the proof of Theorem \ref{theorem2} below, we will need well-known results that  combine to make checking the WLP a task of checking a few specific maps.\begin{lemma}\cite[Lemma 2.1 and Corollary 2.2]{CN}
\label{lem:WLPcheck}
Let $S$ be a standard graded ring and let $M$ be a graded $S$-module such that the degrees of its minimal generators are at most $v$. Let $\ell\in S$ be a linear form. If the map $\times\ell\colon M_{v-1}\rightarrow M_v$ is surjective, then the map  $\times\ell\colon M_{j-1}\rightarrow M_j$ is surjective for all $j \geq v$.

Let $M$ be an Artinian graded $S$-module such that the degrees of its non-trivial socle elements are at least $u-1$. Let $\ell\in S$ be a linear form. If the map $\times \ell\colon M_{u-1}\rightarrow M_u$ is injective, then  the map $\times \ell\colon M_{j-1}\rightarrow M_j$ is injective for all $j \leq u$.
\end{lemma}

We point out the necessity of the hypothesis that $S$ has standard grading for the proof of the previous lemma. Indeed, the proof relies on the fact that for $\ell\in S_1$, having $j\geq v$ and $(M/\ell)_j=0$ implies that $(M/\ell)_{j+1}=0$ since then the multiplication map $S_1\otimes (M/\ell)_j\rightarrow (M/\ell)_{j+1}$ is surjective for $j\geq v$.

With an eye towards applying Lemma \ref{lem:WLPcheck} for fibered products and connected sums, we establish bounds on the degrees of generators of these rings as modules over the standard graded polynomial algebra generated by their degree 1 components.

\begin{lemma}
\label{lem:DegreeGens}
Let $A$, $B$ be standard graded AL algebras having socle degree $d$ and let $T$ be a graded AG algebra of socle degree $k$ endowed with surjective $\F$-algebra homomorphisms $\pi_A:A\to T$ and $\pi_B:B\to T$.  Then there exists a standard graded polynomial algebra $S$ such that $A\times_TB$ and $A\#_TB$ are $S$-modules generated in degree at most $k+1$ and with socle concentrated in degree $d$.
\end{lemma}
\begin{proof}
Set $K_A=\ker(\pi_A), K_B=\ker(\pi_B)$, let $V$ be a basis for $T_1\cup (K_A)_1\cup (K_B)_1$ and let $Q=\F[V]$.
 We claim that there are ideals $I_A, I_B$ and $I_T$ such that $A\cong Q/I_A, B\cong Q/I_B$ and $T\cong Q/I_T$ and $\pi_A, \pi_B$ are the canonical projection maps between the respective cyclic modules. Indeed, the restrictions of the maps $\pi_A$ and $\pi_B$ to $A_1$ and $B_1$ respectively induce vector space isomorphisms $A_1\cong T_1\oplus (K_A)_1$ and $B_1\cong T_1\oplus (K_B)_1$ which give rise to surjective $\F$-algebra maps $q_A:Q\to A, q_B:Q\to B$ and $q_T:Q\to T$, since $A, B$ and $T$ are standard graded. Setting $I_A, I_B, I_T$ to be the kernels of these maps gives the claimed isomorphisms $\overline{q_A}:Q/I_A\to A$ etc. Furthermore, the induced map $\pi_A\circ \overline{q_A}$ is a canonical projection by construction and similarly for $\pi_A\circ \overline{q_B}$. 

Notice that 
$$(A\times_TB)_1=\{(\overline{q_A}(t),\overline{q_B}(t))\mid t\in T_1\} \oplus \left((K_A)_1\times 0 \right) \oplus \left(0 \times (K_B)_1 \right), $$
let $U$ be an $\F$-basis for this vector space and set $S=\F[U]$ to be a polynomial ring with the standard grading. Let $G$ be an $\F$-basis for $(A\times_T B)_{\leq k+1}$.
We claim that $G$ generates both $A\times_TB$ and  $A\#_TB$ as $S$-modules. It is sufficient to establish this for the former module, since the latter is its quotient. Notice that $(A\times_TB)_{\leq k+1}\subseteq SG$ by definition of $G$ and $(A\times_TB)_i=(A_i\times 0) \oplus (0\times B_i)$ for $i\geq k+1$ and $A_i\times 0=S_{i-k-1}(A_{k+1}\times 0)$ and $0\times B_i=S_{i-k-1}(0\times B_{k+1})$ since the projection of $S_1$ onto the first coordinate is $A_1$ while the projection onto the second coordinate is $B_1$ and both $A$ and $B$ are standard graded. This shows that $(A\times_T B)_{\geq k+1}\subset SG$ as well, yielding the desired statement about the generator degrees of  $A\times_T B$ as an $S$-module.

The proof of Lemma \ref{lem:Level} establishes the claim about the socle degree.
\end{proof}

We are now in a position to prove our second main result, stated as Theorem \ref{thm:B} in the introduction, which we restate here for convenience.
\begin{theorem}\label{theorem2}
Let $A$ and $B$ be standard graded AG algebras of socle degree $d$ satisfying the SLP, and let $T$ be a graded AG algebra of socle degree $k$, with $k<\lfloor \frac{d-1}{2} \rfloor$,
endowed with surjective $\F$-algebra homomorphisms $\pi_A:A\to T$ and $\pi_B:B\to T$. Then the resulting fibered product  $A\times_TB$ and the connected sum $A\#_T B$ both satisfy the WLP.
\end{theorem}
\begin{proof}
We aim to apply  Lemma \ref{lem:WLPcheck} for the ``middle degrees", i.e.,
$$
\begin{cases}
 u=\frac{d}{2} \text{ and } v=\frac{d}{2}+1 & \text{ if } d \text{ is even}\\
u=v=\lfloor\frac{d}{2}\rfloor & \text{ if } d \text{ is odd}.
\end{cases}
$$
 We proceed in several steps. First we show that  $A\times_TB$ and $A\#_T B$ satisfy the hypotheses of the Lemma \ref{lem:WLPcheck}. Next we show that, due to our hypothesis on the socle degree of $T$, the multiplication maps by a linear form on  $A\times_TB$, $A\#_T B$ and $A\times_\F B$ coincide in degrees $u-1$ to $u$ and $v-1$ to $v$. Finally, using the WLP for $A\times_\F B$ we conclude  $A\times_TB$ and $A\#_T B$ have the WLP.

{\em Step 1:} It follows from Lemma \ref{lem:DegreeGens} that there is a standard graded polynomial ring $S$ generated by $(A\times_TB)_1$ such that $A\times_TB$ and $A\#_T B$ are $S$ modules generated in degrees at most $k+1$ and with socle in degree $d$. Because of the assumption that $k<\lfloor \frac{d-1}{2} \rfloor$, the conditions $k+1\leq v$ and $d\geq u$ are satisfied.

{\em Step 2:} We show that $A\times_T B, A\times_F B$ and $A\#_T B$ have certain graded components in common, in particular $(A\times_T B)_i\cong(A\#_T B)_i$ whenever $i\in\{u-1,u,v-1,v\}$. 
Indeed, from Equations \eqref{eq:SESFP} and \eqref{eq:SESCS} and our assumption on $k$ we  have that
 $$(A\#_T B)_i \cong (A\times_T B)_i\cong (A\#_\F B)_i\cong (A\times_\F B)_i\cong A_i\times B_i$$
 for $i\in\{u-1,u,v-1,v\}$ for the above defined values of $u$ and $v$.

{\em Step 3:} By Lemma \ref{lem:WLPcheck}, in order to establish that WLP holds both for $A\times_TB$ and $A\#_T B$,  it suffices to check that there exists $L\in(D=A\times_T B)_1$ satisfying the following properties :
\begin{itemize}
\item $D_{u-1}\stackrel{\times L}{\longrightarrow} D_u$ is injective 
\item $D_{v-1}\stackrel{\times L}{\longrightarrow} D_v$ is surjective.

\end{itemize}
Take $L$ to be a Lefschetz element for $A\times_\F B$ which is also in $(A\times_T B)_1$. It is possible to choose such an element $L$ because of Remark \ref{rem:LefschetzLocusProd} and Lemma~\ref{lem:survival}. The inclusion $D=A\times_T B\hookrightarrow A\times_\F B=D'$ of Lemma \ref{lem:injFiberProd} induces an isomorphisms between the graded components of $A\#_T B,A\times_T B$ and $A\times_\F B$ in degree $i\in\{u-1,u,v-1,v\}$, as established in Step 2. Since $D'$ satisfies the SLP by Theorem \ref{prop:SLPFP}, the map $D'_{u-1}\stackrel{\times L}{\longrightarrow} D'_u$ is injective and the map  $D'_{v-1}\stackrel{\times L}{\longrightarrow} D'_v$ is surjective. Using the isomorphisms $D_i\cong D'_i$ for $i\in\{u-1,u,v-1,v\}$ it follows that multiplication by $L$ induces maximal rank maps on $D$ in the desired degrees.
 \end{proof}

We next discuss the necessity of the hypothesis on the socle degrees of $A,B$ and $T$ in the previous theorem, as well as whether this theorem can be extended to cover the SLP. We give a family of algebras showing that the SLP does not hold in general for connected sums, even those satisfying the hypotheses of the previous theorem. Our next result also shows that in the absence of the hypothesis on the socle degrees of $A,B$ and $T$ given in Theorem \ref{theorem2}, the connected sums may or may not satisfy the WLP.

\begin{proposition}
	\label{ex:NonSLP}
Let $A=\F[x]/(x^m)$, $B=\F[y]/(y^m)$, $T=\F[z]/(z^t)$ (where $m>t>1$), and let $\pi_A\colon A\rightarrow T$ and $\pi_B\colon B\rightarrow T$ be the algebra maps defined by $\pi_A(x)=z$ and  $\pi_B(y)=z$.  Let $\int_A\colon x^{m-1}\mapsto 1$, $\int_B\colon y^{m-1}\mapsto 1$, $\int_T\colon z^{t-1}\mapsto 1$ be orientations for $A$, $B$, and $T$, respectively so that the Thom classes for $\pi_A$ and $\pi_B$ are $\tau_A=x^{m-t}$ and $\tau_B=y^{m-t}$.  

Then $A\#_T B$ has the WLP if and only if $t\neq \frac{m}{2}$, but $A\#_T B$ never has the SLP. 
\end{proposition}

\begin{proof}
Under our hypothesis the fibered product $A\times_T B$ has a presentation given by
$$A\times_T B=\frac{\F[z_1,z_2]}{\left(z_1^m,z_2^{\left\lceil\frac{m}{t}\right\rceil},z_1^{m-t}z_2, z_1^tz_2-z_2^2\right)} \text{ where } \begin{cases} z_1= & (x,y)\\
z_2= & (x^t,0)\\
\end{cases}.$$

Note that $(\tau_A,\tau_B)=z_1^{m-t}$, hence we get the following  presentation for the connected sum:

\begin{equation}\label{Ceqn}C=A\#_T B=\frac{\F[z_1,z_2]}{\left(z_1^{m-t}, z_1^tz_2-z_2^2\right)}.
\end{equation}

Additionally we observe that the Macaulay dual generators for the fibered product are given by 
	$$\begin{cases}
		H_1= & Z_1^{m-1}\\
		H_2= & Z_1^{m-1-t}Z_2+Z_1^{m-1-2t}Z_2^2+\cdots+Z_1^{m-1-\left(\left\lceil\frac{m}{t}\right\rceil-1\right)t} Z_2^{\left\lceil\frac{m}{t}\right\rceil-1}\\
	\end{cases}$$
With this information, we can compute the ``correct'' Macaulay dual generators for $A$, $B$, and $T$ by looking at the projections from $A\times_TB$:
\begin{align*}
A\cong & \frac{\F[z_1,z_2]}{\Ann(F=H_1+H_2)}\cong\frac{\F[z_1,z_2]}{(z_1^m,z_1^t-z_2)}\\
B\cong & \frac{\F[z_1,z_2]}{\Ann(G=H_1)}\cong \frac{\F[z_1,z_2]}{(z_1^m,z_2)}\\
T\cong & \frac{\F[z_1,z_2]}{\Ann(\tau\circ F=Z_1^{t-1}=\tau\circ G)}\cong\frac{\F[z_1,z_2]}{(z_1^t,z_2)}.
\end{align*}
The Macaulay dual generator of the connected sum is $H_2=F-G$, i.e.
$$C=\frac{\F[z_1,z_2]}{\Ann\left(Z_1^{m-1-t}Z_2+\cdots+Z_1^{m-1-\left(\left\lceil\frac{m}{t}\right\rceil-1\right)t} Z_2^{\left\lceil\frac{m}{t}\right\rceil-1}\right)}\cong \frac{\F[z_1,z_2]}{(z_1^{m-t},z_1^tz_2-z_2^2)}.$$
Since $C$ is an Artinian Gorenstein algebra of embedding dimension two, it must be a complete intersection (of socle degree $m-1$).  From Lemma~\ref{lem:HFCS}, the Hilbert function for $C$ is
\begin{equation}\label{HF2eqn}
H(C)=\begin{cases}
(\underbrace{1,\ldots, 1}_{t}, \underbrace{2, \ldots,2}_{m-2t},\underbrace{1,\ldots, 1}_{t}) &\text{ if } t<\frac{m}{2}\\
(\underbrace{1,\ldots, 1}_{m}) &\text{ if } t=\frac{m}{2}\\
(\underbrace{1,\ldots, 1}_{m-t}, \underbrace{0, \ldots, 0}_{2t-m},\underbrace{1,\ldots, 1}_{m-t}) &\text{ if } t>\frac{m}{2}\\
\end{cases} .
\end{equation}
Clearly the only candidate for a WL or SL element in $C=A\#_T B$ is a multiple of $L=z_1\in C_1$. Notice that the map  $\times L^{m-1}\colon C_0\rightarrow C_{m-1}$ is the zero map since $L^{m-1}=0$, thus $C$ never has the SLP. To study the WLP, we first note that a $\F$-basis for $C$ consists of the monomials 
$$\mathcal{B}=\{z_1^i, z_1^iz_2 \ | \ 1\leq i\leq m-t-1\}.$$ 

Examining the multiplication map $\times L\colon C_{m-t-1}\rightarrow C_{m-t}$, which takes $z_1^{m-t-1}$  to $0$ shows that $C$ does not have the WLP in the case $t=\frac{m}{2}$. In all the other cases the multiplication maps $\times L\colon C_{i}\rightarrow C_{i+1}$ have maximal rank; this can be easily seen by expressing these multiplication maps in terms of the given basis. 
\end{proof}

Equation \eqref{HF2eqn} implies that when both $2\le t$ and $t\not=m/2$ then $A,B,T$ are standard graded, but $C$ is not standard graded.  If $t=1$ then $C$ in \eqref{Ceqn} is a standard graded connected sum over $T=\F$ and it has both WLP and SLP provided the characteristic of $\F$ is sufficiently large.
Some topological implications of Proposition \ref{ex:NonSLP} are discussed in Remark \ref{rem:CSNP}.

\begin{definition}
\label{def:Jordantype}
The {\em Jordan type} of an Artinian graded algebra $A$ with $A_1\neq 0$ is the multiset of sizes of blocks in the Jordan matrix representing the (nilpotent) action  of multiplication by a general linear form on $A$. The Jordan type is a partition of the length (vector space dimension) of $A$ .
\end{definition}

In the last part of this section we focus on the interplay between the property of an algebra of being decomposable as a connected sum and its Jordan type. The importance of the Jordan type in the theory of the Lefschetz properties for Artinian graded algebras is given by the following remark. 

\begin{remark} Let $\ell$ be a general enough non-unit of $A$ and consider the multiplication map $m_\ell: A \to A$; denote by $P_\ell$ its Jordan type. Recall from \cite[Proposition~2.10]{IMM}  (or \cite[Proposition~3.64]{HMMNWW} for $A$ standard graded) that for a (not necessarily standard) graded Artinian algebra $A$, the SLP is equivalent to the fact that $P_\ell$ is the conjugate of the partition of the length of $A$ given by the Hilbert function of $A$.\par
For standard graded algebras $A$ having symmetric Hilbert function the WLP is equivalent to the number of parts of $P_\ell$ being equal to the largest value of the Hilbert function of $A$ (\cite[Proposition~3.5]{HMMNWW}). 

Proposition \ref{ex:NonSLP} shows that this WLP criterion is no longer true for non standard graded algebras, even in the case of complete intersections. 
For the connected sums discussed in Proposition~\ref{ex:NonSLP}, the Jordan matrix of the action of multiplication by $z_1$ on $A\#_TB$ expressed in terms of the basis $\mathcal{B}$ has two blocks of size $m-t$, corresponding to the two subsets $\mathcal{B}_1=\{z_1^i \ | \ 0\leq i\leq m-t-1\}$ and $\mathcal{B}_2=\{z_1^iz_2 \ | \ 0\leq i\leq m-t-1\}$ of the basis. This Jordan type remains the same independently of whether  $A\#_TB$ has the WLP (when $d>m/2$ ) or not.
Furthermore, when $d>m/2$ then $A\#_TB$ has the WLP, the number of Jordan blocks is two, and the largest value of the Hilbert function $H(A\#_TB)$ is one.
\end{remark}

In the following proposition we establish a partial converse to Proposition \ref{ex:NonSLP}.

\begin{proposition}
\label{prop:2blocks}
Let $C$ be a graded AG $\F$-algebra with $\F$-algebraically closed and $C_1\neq 0$. Assume that the Jordan type of $C$ consists of two equal parts of size $a>1$, which correspond to basis elements of degrees $0$ to $a-1$ and $t$ to $t+a-1$ respectively for some positive integer $t>0$.  Then one of the following possibilities is true:
\begin{enumerate}[(i).]
\item $C\cong \F[u,v]/(u^a, v^2)$ with $\deg(u)=1$, $\deg(v)=t$.
\item $C\cong  \F[u,v]/(u^a, v^2-u^tv)$  with $\deg(u)=1$, $\deg(v)=t$. In this case, $C$ is a connected sum $C \cong A\#_TB$ for $A=\F[x]/(x^m), B=\F[y]/(y^m),T=\F[z]/(z^t)$, $\pi_A(x)=z$, $\pi_B(y)=z$ and $a=m-t$.
\end{enumerate}
In all cases, $C$ does not have the strong Lefschetz property, and if $t=1$ then the field characteristic necessarily satisfies $0<\cha(\F)\leq a$.
\end{proposition}

\begin{proof}
Let $u\in C_1$ be a general linear form. The hypothesis yields that a basis for $C$ is given by $\mathcal{B}=\{1,u,\ldots, u^{a-1}, v,uv,\ldots, u^{a-1}v\}$ for some form $v\in C_t$ where $t$ is a positive integer. Note that $C$ cannot have the strong Lefschetz property since evidently the socle degree is $t+a-1\geq a$ but $u^a=0$ for a general linear form.  Also recall that standard graded AG algebras in two variables have SLP if $\cha(\F)=0$ or $\cha(\F)>a$ \cite[Theorem 2.9]{Ia}.

{\em Case 1:} $a \leq t$. In this case the socle degree of $C$ is $t+a-1<2t$, thus $C_{2t}=0$ and hence $v^2=0$. Therefore $C=\F[u,v]/(u^a,v^2)=\F[u,v]/(u^a,v^2-u^tv)$ where the second equality above follows from the assumption $a\leq t$.  Note in this case $C$ necessarily has a non-standard grading since $t\geq a>1$.   

{\em Case 2:} $t<a\leq 2t$. The hypothesis and the numerical constraints on $a$ and $t$ yield that $C_{2t}$ is 1-dimensional and it contains the monomials $v^{2}$ and $u^tv\neq 0$. Thus in $C$ there is a relation of the form $v^2-\alpha u^tv=0$ where $\alpha\in F$. If $\alpha=0$ then we are in case (i) since there is an obvious surjection $k[u,v]/(u^a, v^2)\to C$  and the two rings have the same length. If $\alpha\neq 0$ then we are in case (ii). Indeed, set $u'=\sqrt[t]{\alpha}u$ and notice that there is a surjection $\F[u',v]/((u')^a, v^2-(u')^tv)\to C$, which by comparing vector space dimensions must be an isomorphism. %It remains to observe that by equation \eqref{Ceqn} the ring $\F[u',v]/((u')^a, v^2-(u')^tv)$ is a connected sum of the form described in (2).
In this case $C$ may or may not have a standard grading.

{\em Case 3:} $a>2t$. The hypothesis and the numerical constraints on $a$ and $t$ yield that $C_{2t}$ is 2-dimensional and it contains the monomials $v^{2},u^tv, u^{2t}$, where $u^{2t}, u^tv$ are linearly independent. Thus in $C$ there is a relation of the form $v^2+\alpha u^tv+\beta u^{2t}=0$ where $\alpha,\beta\in \F$. Since $\F$ is algebraically closed this relation can be factored as $(v-\delta u^t)(v-\epsilon u^t)=0$. If $\delta=\epsilon$, then setting $v'=v-\delta u^t$ gives  $C=\F[u,v']/(u^a,(v')^2)$. Otherwise, setting $v'=v-\delta u^t$, $u'=\sqrt[t]{\epsilon-\delta} \cdot u$ gives $C=\F[u',v']/((u')^a,(v')^2-u^tv')$.  In this case $C$ may or may not have a standard grading.

Note that the two rings $ \F[u,v]/(u^a, v^2)$ and $ \F[u,v]/(u^a, v^2-u^tv)$ are not isomorphic if $a>2t$ since the defining ideal of the former contains the square of a degree $t$ form, whereas the degree $2t$ component of the latter does not.  If $t<a\leq 2t$, the rings $\F[u,v]/(u^a, v^2)$ and $\F[u,v]/(u^a, v^2-u^tv)$ are isomorphic as one can see by noticing that in the latter $0=v^2-u^tv=v^2-u^tv+\frac{1}{4}u^{2t}=(v')^2$ where $v'=v-\frac{1}{2} u^t$.  In particular, the numerical condition $a\leq 2t$ forces $C$ to be a connected sum as in \eqref{Ceqn}.
%Conversely, if $\deg(v)=t>1$,  $u$ is the only general form up to scalar multiple. For each of the algebras listed in the statement, the Jordan type of multiplication by $u$ has two equal parts by construction. However, if $\deg(v)=t=1$ then the Jordan type of multiplication by a general linear form, e.g. $u+v$, on $\F[u,v]/(u^a, v^2)$ if $\cha(\F)$ does not divide $a$ as well as on $ \F[u,v]/(u^a, v^2-u^tv)$ regardless of $\cha(\F)$ is $(a,a-2)$. This eliminates these possibilities from our classification and completes the proof of the proposition.
\end{proof}

\begin{remark}
	\label{cor:2blocks}
	Proposition \ref{prop:2blocks} shows that Jordan type can impose severe restrictions on the structure of an AG algebra.  In a similar vein, one can also show that if an AG algebra $C$ has $C_1\neq 0$ and generic Jordan type $(a,b)$ with two unequal parts that differ by two, then $C$ has SLP if and only if $C$ is standard graded.  Indeed if $C$ has SLP with socle degree $s$, then a general linear form $u\in C_1$ satisfies $u^s\neq 0$, and from its Jordan basis $\{1,u,\ldots,u^{a},v,uv,\ldots,u^{b}v\}$ we conclude that $a=s$ and hence $\deg(v)=t< a-b=2$.  Conversely if $C$ has standard grading then $\deg(v)=t=1$ and $C=\F[u,v]/I$ has SLP by \cite[Proposition~3.64]{HMMNWW}. 
	
	%The WLP for connected sums of the type relevant for case (ii) of Proposition~\ref{prop:2blocks} is analyzed in Proposition \ref{ex:NonSLP}. The same proof extends to study the WLP for the AG algebra $\F[u,v]/(u^a,v^2)$ with $\deg(u)=1$, $\deg(v)=t>1$. In either case, the conclusion is that, $C$ has the WLP if and only if $t\neq \frac{a+1}{2}$. 
\end{remark}

\begin{remark}\label{WLCIrem} 
An important open problem in the study of the algebraic Lefschetz properties is establishing whether {\em all} Artinian standard graded complete intersection (CI) algebras satisfy the WLP -- or, even more ambitious, might satisfy the SLP -- or finding a counterexample. In codimension three the WLP for such CI algebras when $\cha \F=0$ is shown in \cite{HMNW}. 
Proposition~\ref{prop:2blocks} singles out a certain class of AG algebras of codimension two, which are automatically complete intersections by the Hilbert-Burch theorem. Standard graded Artinian algebras of codimension two all satisfy the SLP when $\cha \F=0$ or is greater than the socle degree $d$. Proposition~\ref{prop:2blocks} shows that there are many non-standard graded CI algebras, even in codimension two, that fail to posses the SLP.
\end{remark}

We now consider the hypothesis that $\F$ be algebraically closed in Proposition \ref{prop:2blocks}.
 We had already discussed the simpler analogue for the case $T=\F$ in Remark \ref{decomposerem}. 
\begin{example}[Dependence of $C$ being a connected sum on the field $\F$]\label{CconnectFex}  Consider the Artinian algebra $C=\mathbb Q[u,v]/(u^5, v^2+u^4-u^2v)$, where $\deg u=1, \deg v=2$. Evidently, $C$ is a complete intersection, the Jordan type of multiplication by $u$ is $(5,5)$ with two strings, the classes in $C$ of $1,u,u^2,u^3,u^4; v,uv,u^2,u^3v,u^4v$: thus, $C$ satisfies the hypotheses of Proposition \ref{prop:2blocks} with $a=5, t=2$, except the hypothesis of the closure of $\F$. Over $\mathbb Q$ the degree $4$ form $f_4=v^2+u^4-u^2v=(v-u^2/2)^2+3u^4/4$ is irreducible, but over $\F=\mathbb Q(\omega), \omega=\sqrt{-3}$ it can be factored as $$f_4=\left(v-u^2/2+(\omega/2)u^2\right)\cdot \left(v-u^2/2-(\omega/2)u^2\right).$$ 
	In the latter case, after a suitable change of variables, $C$ has the form (2) of Proposition~\ref{prop:2blocks}, so is a connected sum of $A=\F[x]/(x^7)=B$ over $T=\F[z]/(z^2)$ as in \eqref{Ceqn}. However, over $\Q$, one can show that there is no such isomorphism, i.e. 
	$$C=\frac{\Q[u,v]}{(u^5,v^2+u^4-u^2v)}\not\cong \frac{\Q[z_1,z_2]}{(z_1^5,z_2^2-z_1^2z_2)}=A\#_TB$$
	and hence over $\Q$, $C$ is not a connected sum as in \eqref{Ceqn}.
\par	We are not aware whether the $\mathbb Q$-algebra $C$ is decomposable as a connected sum over some other $T$ after a change of coordinates.
\end{example}

 \appendix
 \section{Cohomology of the Connected Sums of Manifolds.}
 \subsection{Topological Construction.}
 Suppose that $M_1$ and $M_2$ are two smooth, connected, compact, orientable $2d$-dimensional manifolds, and suppose that $N$ is a smooth compact connected orientable $2k$-dimensional manifold.  Furthermore, suppose that $\iota_1,\iota_2\colon N\rightarrow M_1,M_2$ are smooth embeddings with images $N_1,N_2\subset M_1,M_2$.  Let $\pi_i\colon\mathcal{E}_i\rightarrow N_i$, $i=1,2$ denote the (orientable) normal bundles of $N_i$ in $M_i$.  By the tubular neighborhood theorem (e.g. \cite[Theorem 11.1]{MS1}), there is an open neighborhood $U_i\subset M_i$ containing $N_i$ which is diffeomorphic to the total space of the normal bundle $\mathcal{E}_i$, under which the submanifold $N_i$ is identified with the zero section of $\mathcal{E}_i$.  We will \textbf{assume} that we have an isomorphism of normal bundles $\phi\colon \mathcal{E}_2\rightarrow \mathcal{E}_1$, and hence a diffeomorphism of tubular neighborhoods $\phi\colon U_2\rightarrow U_1$ which restricts to a diffeomorphism $\phi|_{N_2}\colon N_2\rightarrow N_1$, making the diagram commute 
 $$\xymatrix{& N\ar[dl]_-{\iota_2}\ar[dr]^-{\iota_1} & \\
 	N_2\ar[rr]_-{\phi|_{N_2}}& & N_1.\\}$$
 Define the \emph{(topological) fibered product} $M_1\times_N M_2$ as the adjunction space obtained by gluing $M_1$ and $M_2$ along $U_1$ and $U_2$ via $\phi$: 
 \begin{equation}
 \label{eq:FPspace}
 M_1\times_NM_2=\frac{M_1\sqcup M_2}{\phi(x)\sim x}.
 \end{equation}  
 
 Note that $M_1\times_N M_2$ is not Hausdorff.  To see this, fix a point $z\in\partial U_2$ in the boundary of $U_2$, and let $\left\{z_n\right\}_{n=0}^\infty\subset U_2$ be a sequence that converges to $z$.  Then $\left\{\phi(z_n)\right\}_{n=0}^\infty\subset U_1$ is a sequence of points of $U_1$ converging to some point $w\in\partial U_1$ on the boundary of $U_1$.  Then $z,w\in M_1\times_NM_2$ are two distinct points in the fibered product that cannot be separated by disjoint open sets.
 
 We can remedy this non-Hausdorff issue in the following way:  Let $\mathcal{E}_0$ denote the total space of the vector bundle minus the zero section, so that $\mathcal{E}_0\cong U\setminus N$.  Assuming that we have defined some metric on $\mathcal{E}=\mathcal{E}_2$,  we may define an ``orientation reversing'' bundle isomorphism
 $$\alpha\colon\mathcal{E}_0\rightarrow\mathcal{E}_0, \ \alpha((x,v))=\left(x,\frac{v}{|v|^2}\right).$$
 Identifying $U_2\setminus N_2\cong \mathcal{E}_0\cong U_1\setminus N_1$, we get an orientation reversing diffeomorphism of deleted tubular neighborhoods
 $$\psi=\alpha\circ\phi\colon U_2\setminus N_2\rightarrow U_1\setminus N_1.$$
 The \emph{(topological) connected sum} $M_1\#_NM_2$ is defined to be the adjunction space.
 \begin{equation}
 \label{eq:CSspace}
 M_1\#_NM_2=\frac{\left(M_1\setminus N_1\right)\sqcup\left(M_2\setminus N_2\right)}{\psi(x)\sim x}.
 \end{equation}
 Geometrically, we are gluing the punctured neighborhoods $V_1= U_1\setminus N_1$ and $V_2= U_2\setminus N_2$ by identifying points close to the ``zero boundary`` (i.e. $N_1$) of $V_1$ with points close to the ``infinity boundary'' (i.e. $\partial U_2$) of $V_2$.  The connected sum is a Hausdorff topological space which can be endowed with a smooth structure compatible with smooth structures on $M_1$ and $M_2$, meaning that the natural inclusions $M_i\setminus N_i\hookrightarrow M_1\#_NM_2$ are smooth (open) maps (embeddings) \cite[Chapter VI]{Kos}.

\subsection{Cohomology Computation.}
For a topological space $X$, let $H^i(X)=H^i(X,\Q)$ denote the $i^{th}$ singular cohomology group of $X$ with coefficients in $\Q$. We write $H^*(X)=\bigoplus_{i\in\Z}H^i(X)$ for the direct sum of all cohomology groups, which has the structure of a graded-commutative ring via the cup product, i.e. for $\alpha\in H^i(X)$ and $\beta\in H^j(X)$, we have $\alpha\cdot\beta=(-1)^{i+j}\beta\cdot\alpha\in H^{i+j}(X)$.  In particular, the even degree part $H^{2*}(X)=\bigoplus_{i\in\Z}H^{2i}(X)$ forms a commutative ring.  In fact, if $X$ is a smooth compact orientable manifold of even dimension $2d$, then its even degree cohomology $A_*=H^{2*}(X)$ is a graded AG algebra with socle degree $d$ (halving degrees).  We write $\Cb(X)$ for the complex of cochain groups (over $\Q$) of $X$, so that the cohomology of $X$ is the cohomology of that complex, i.e. $H^*(X)=H^*(\Cb(X))$.    
\begin{theorem}
	\label{thm:TopCS}
	Let $M_1$, $M_2$, and $N$ be as above, with fibered product $M_1\times_NM_2$ as in \eqref{eq:FPspace} and connected sum $M_1\#_NM_2$ as in \eqref{eq:CSspace}. 
	Assume further that 
	\begin{enumerate}[(i).]
		\item the odd degree cohomology groups of $M_1$ and $M_2$ vanish, i.e.
		$$H^{2q-1}(M_i)=0, \ i=1,2, \ \forall q\in\Z.$$
		\item the smooth embeddings $\iota_i\colon N\hookrightarrow M_i$ induce surjections on cohomology 
		$$\iota_i^*\colon H^*(M_i)\rightarrow H^*(N), \ i=1,2.$$
	\end{enumerate}
	Fix orientations on $M_1$, $M_2$, and $N$, and set
	$$A_*\coloneqq H^{2*}(M_1), \ B_*\coloneqq H^{2*}(M_2), \ T_*\coloneqq H^{2*}(N),$$
	with projection maps $\pi_A\coloneqq \iota_1^*\colon A\rightarrow T$ and $\pi_B\coloneqq \iota_2^*\colon B\rightarrow T$.  Then 
	\begin{enumerate}[(a).]
		\item the Thom classes for the maps $\pi_A$ and $\pi_B$ are correspond to the Thom classes for the oriented normal bundles $\mathcal{E}_i\rightarrow N_i$
		\item the algebraic fibered product $A\times_TB$ is isomorphic as a ring to the (even degree) cohomology ring of the topological fibered product $M_1\times_NM_2$, i.e.
		$$\left(A\times_TB\right)_*\cong H^{2*}(M_1\times_NM_2).$$
		\item the algebraic connected sum $A\#_TB$ is isomorphic as a graded vector space to the (even degree) cohomology of the topological connected sum $M_1\#_NM_2$, i.e. 
		$$\left(A\#_TB\right)_*\cong H^{2*}(M_1\#_NM_2).$$
	\end{enumerate}
\end{theorem}

Before we prove Theorem \ref{thm:TopCS}, we review some basic facts from algebraic topology.  Our main reference is J.W. Milnor and J.D. Stasheff's book \cite{MS1}.

\subsection{Exact Sequences.}
There are two short exact sequences of cochain complexes of topological spaces that we need:  the exact sequence for pairs, and the Mayer-Vietoris sequence.

Recall that if $X$ is a topological space and $A\subset X$ is a subspace, there is a short exact sequence of cochain complexes, called the short exact sequence for the pair $(X,A)$: 

\begin{equation}
\label{eq:SESpair}
\xymatrix{0\ar[r] & \mathcal{C}^\bullet(X,A)\ar[r] & \mathcal{C}^\bullet(X)\ar[r] & \mathcal{C}^\bullet(A)\ar[r] & 0.}
\end{equation} 
Sequence \eqref{eq:SESpair} yields a long exact sequence of cohomology groups
\begin{equation}
\label{eq:LESpair}
\xymatrix{\cdots\ar[r] & H^j(X,A)\ar[r] & H^j(X)\ar[r] & H^j(A)\ar[r] & H^{j+1}(X,A)\ar[r] & \cdots}
\end{equation}

Recall that if $X$ is a topological space and $U,V\subset X$ are two open sets which cover $X$, i.e. $X=U\cap V$, then we get a short exact sequence of cochain complexes called the Mayer-Vietoris sequence:
\begin{equation}
\label{eq:SESMV}
\xymatrix{0\ar[r] & \mathcal{C}^\bullet(X)\ar[r] & \mathcal{C}^\bullet(U)\oplus\mathcal{C}^\bullet(V)\ar[r]^-{\rho_U^*-\rho_V^*} & \mathcal{C}^\bullet(U\cap V)\ar[r] & 0}
\end{equation}
where $\rho_U,\rho_V\colon U\cap V\hookrightarrow U,V$ are the natural inclusion maps.

Sequence \eqref{eq:SESMV} also yields a long exact sequence of cohomology groups
\begin{equation}
\label{eq:LESMV}
\xymatrix{\cdots\ar[r] & H^j(X)\ar[r] & H^j(U)\oplus H^j(V)\ar[r] & H^j(U\cap V)\ar[r] & \cdots}
\end{equation}

\subsection{Thom Classes and Gysin Maps.}
Suppose that $M$ is a smooth compact connected oriented $2d$-dimensional manifold $M$ and that $N\subseteq M$ is a smooth compact connected oriented $2k$-dimensional embedded submanifold with normal bundle $\pi\colon\mathcal{E}\rightarrow N$.  Let $\pi_0\colon \mathcal{E}_0\rightarrow N$ be the deleted normal bundle, i.e. the normal bundle with the zero section removed.  In this case, the normal bundle is orientable itself (cf. \cite[p. 66]{BT}), meaning that each fiber $\mathcal{E}_x\subset\mathcal{E}$ of the normal bundle has a preferred orientation class $\mu_x\in H^{2(d-k)}(\mathcal{E}_x,\left(\mathcal{E}_x\right)_0)$ such that each point $x\in N$ has some open neighborhood $V\subset N$ and a cohomology class $\mu_V\in H^{2(d-k)}(\pi^{-1}(V),\pi^{-1}(V)_0)$ which restricts to the preferred generator $\mu_y\in H^{2(d-k)}(\mathcal{E}_y,(\mathcal{E}_y)_0)$ for each $y\in V$, (cf. \cite[p. 96]{MS1}). 

The following fact is referred to by Milnor-Stasheff as the Thom isomorphism theorem \cite[Theorem 10.4]{MS1}:
\begin{fact}[Thom isomorphism theorem]
	\label{lem:Thom}
	Fix an orientation for the normal bundle $\pi\colon\mathcal{E}\rightarrow N$ with preferred orientation classes $\mu_x\in H^{2(d-k)}(\mathcal{E}_x\left(\mathcal{E}_x\right)_0)$ as above. Then there is a unique cohomology class $\tau\in H^{2(d-k)}(\mathcal{E},\mathcal{E}_0)$ with the property that its restriction to each fiber $\tau_x=\mu_x\in H^{2(d-k)}(\mathcal{E}_x,\left(\mathcal{E}_x\right)_0)$ is the preferred orientation class for that fiber.  Moreover there is a well defined map 
	\begin{equation}
	\label{eq:Thom}
	\xymatrixrowsep{.5pc}\xymatrix{H^j(\mathcal{E})\ar[r] & H^{j+2(d-k)}(\mathcal{E},\mathcal{E}_0)\\
		y\ar@{|->}[r] & y\cdot \tau}
	\end{equation}	
	that is an isomorphism for every $j\in\Z$.
\end{fact}
This class $\tau\in H^{2(d-k)}(\mathcal{E},\mathcal{E}_0)$ is called the \emph{Thom class} of the oriented normal bundle $\pi\colon \mathcal{E}\rightarrow N$, and the isomorphism $H^*(\mathcal{E})\rightarrow H^{*+2(d-k)}(\mathcal{E},\mathcal{E}_0)$ is called the \emph{Thom isomorphism}.  

Note that the bundle map $\pi\colon \mathcal{E}\rightarrow N$ is a retraction onto the zero section, so induces an isomorphism on cohomology rings $\pi^*\colon H^*(N)\overset{\cong}{\rightarrow} H^*(\mathcal{E})$.  Hence the Thom isomorphism can also be given as the map
\begin{equation}
\label{eq:ThomIsom}
H^*(N)(-2(d-k))\cong H^*(\mathcal{E})\rightarrow H^*(\mathcal{E},\mathcal{E}_0), \ y\mapsto y\cdot\tau
\end{equation}
Setting $M^*=M\setminus N$, there is a canonical isomorphism of cohomology rings
$$H^*(\mathcal{E},\mathcal{E}_0)\cong H^*(M,M^*).$$
This evidently follows from an excision argument, cf. \cite[Corollary 11.2]{MS1}.  Also note that the inclusion of pairs $(M,\emptyset)\hookrightarrow (M,M^*)$ induces a map on cohomology rings
\begin{equation}
\label{eq:excision}
H^*(\mathcal{E},\mathcal{E}_0)\cong H^*(M,M^*)\rightarrow H^*(M).
\end{equation}
We shall abuse notation slightly, and use the same letter (i.e. $\tau$) and same name (i.e. Thom class) for the image of the Thom class under this map.  The Thom class $\tau\in H^{2(d-k)}(M)$ for $N$ satisfies the following:
$$\langle x\cup\tau_N,\mu_M\rangle = \langle \iota^*(x),\mu_N\rangle$$
where $\mu_M\in H_{2d}(M)$ (homology) is the fundamental class for $M$, $\mu_N\in H_{2k}(N)$
is the fundamental class for $N$, $\iota^*\colon H^*(M)\rightarrow H^*(N)$ the restriction map induced by the embedding $\iota\colon N\hookrightarrow M$, and $\langle a,b\rangle$ is the natural pairing between cohomology and homology, cf. \cite[Problem 11-C p. 136]{MS1} (see also \cite[p. 67]{BT}).  In particular if we set $A_*=H^{2*}(M)$ and $T=H^{2*}(N)$ oriented AG algebras with orientations $\int_Aa=\langle a,\mu_M\rangle$ and $\int_Tt=\langle t,\mu_N\rangle$, $\iota\colon N\hookrightarrow M$ the embedding of $N$ into $M$, and $\iota^*\colon A_*\cong H^{2*}(M)\rightarrow H^{2*}(N)\cong T_*$ the restriction map, then the (algebraic) Thom class for the map $\iota^*$ coincides exactly with the Thom class of $\tau\in H^{2(d-k)}(M)$ above.

\begin{remark}
	\label{rem:PDGysin}
	The image of the Thom class $\tau\in H^{2(d-k)}(M)$ under Map \eqref{eq:excision} is typically referred to by topologists as the \emph{Poincar\'e dual class} of the submanifold $N\subset M$, cf. \cite[Definition p. 120]{MS1}.  See also \cite[Proposition 6.24]{BT}.  Furthermore, the composition of Maps \eqref{eq:ThomIsom} and \eqref{eq:excision} is referred to, e.g. \cite[p. 212]{Fulton},\footnote{W. Fulton gives a different but equivalent description of Map \eqref{eq:GysinTop} in terms of Poincar\'e duality.} as the \emph{Gysin map} for the inclusion $\iota\colon N\rightarrow M$, i.e.
	\begin{equation}
	\label{eq:GysinTop}
	\xymatrixrowsep{.5pc}\xymatrix{\iota_*\colon H^{2*-2(d-k)}(N)\cong H^{2*}(M,M^*)\ar[r]^-{\times\tau} & H^{2*}(M)\\
	y\ar@{|->}[r] & y\cdot \tau}
	\end{equation}
and is analogous to the ``algebraic Gysin map'' we defined in Definition \ref{def:GysinMap} in Section \ref{sect:Pre}.\footnote{In topology, one also encounters the related \emph{Gysin sequence}, which is obtained from the long exact sequence of cohomology groups corresponding to the pair of spaces $(\mathcal{E},\mathcal{E}_0)$ after applying the Thom isomorphism cf. \cite[p. 143]{MS1}.  See also \cite[Proposition 14.33]{BT}.  Gysin maps for Chow groups and their algebraic analogues have also appeared in the algebraic geometry literature, e.g. \cite{Fulton2}, \cite{Laksov}.}  Note that if the restriction map $\iota^*\colon H^{2*}(M)\rightarrow H^{2*}(N)$ is surjective, then the Gysin map $\iota_*$ is injective by Lemma \ref{lem:GysinMult}.
\end{remark}

\subsection{Proof of Theorem \ref{thm:TopCS}.}
With notation as above, set $M_i^*=M_i\setminus N_i$ and set $U_i^*=U_i\setminus N_i$.  Then we have the following grid of cochain complexes:
\begin{equation}
\label{eq:Grid}
\xymatrixcolsep{2.5pc}\xymatrix{ & 0 & 0 & 0 & \\
	0\ar[r] & \mathcal{C}^\bullet(M_1\#_NM_2)\ar[r]\ar[u] & \Cb(M_1^*)\oplus\Cb(M_2^*)\ar[r]^-{\iota_1^*-\iota_2^*\circ\psi^*}\ar[u] & \Cb(U_1^*)\ar[r] \ar[u] & 0\\
	0\ar[r] & \mathcal{C}^\bullet(M_1\times_NM_2)\ar[r] & \Cb(M_1)\oplus\Cb(M_2)\ar[r]^-{\iota_1^*-\iota_2^*\circ\phi^*}\ar[u] & \Cb(U_1)\ar[r] \ar[u] & 0\\
	0\ar[r] & \Cb(\mathcal{E},\mathcal{E}_0)\ar[r]\ar@{-->}[u]^-{\rho} & \Cb(M_1,M_1^*)\oplus\Cb(M_2,M_2^*)\ar[r]^-{\iota_1^*-\iota_2^*\circ\phi^*}\ar[u] & \Cb(U_1,U_1^*)\ar[r] \ar[u] & 0\\
	& 0\ar[u] & 0\ar[u] & 0\ar[u] & }	
\end{equation}
\normalsize
The rows of Grid \eqref{eq:Grid} are the short exact Mayer-Vietoris sequences (the first map on the bottom row is just the canonical isomorphism on each summand).  The two columns on the right are the short exact sequences for pairs.  The maps in the left column are less clear, and need to be specified.  It follows from the commutativity of the lower right square that there exists a (unique) well-defined (injective) map of complexes $\rho\colon \Cb(\mathcal{E},\mathcal{E}_0)\rightarrow \Cb(M_1\times_NM_2)$ making the lower left square of Grid \eqref{eq:Grid} commute.
Unfortunately, we cannot employ the same technique to get a map between the complexes $\Cb(M_1\times_NM_2)$ and $\Cb(M_1\#_NM_2)$ because the upper right square does not commute.  On the other hand we claim that the upper right square does commute on the level of cohomology!  To wit:
\vskip 0.2cm\noindent
{\bf Claim.}
	\label{lem:Same}
	The two maps of cochain complexes
	\vskip 0.2cm
	\begin{center}
	$\xymatrixcolsep{3pc}\xymatrix{&\\
		\Cb(M_1^*)\oplus\Cb(M_2^*)\ar@/^2pc/[r]^-{\iota_1^*-\iota_2^*\circ\phi^*} \ar@/_2pc/[r]_-{\iota_1^*-\iota_2^*\circ\psi^*}& \Cb(U_1^*)}$
		\end{center}	
		\par
	induce the same map on cohomology groups.\par\noindent
{\it Proof of Claim.}
	For the vector bundle $\mathcal{E}$, let $\mathcal{D}=\left\{(x,v)\in\mathcal{E} \ | \ |v|=1\right\}$ be the associated unit disk bundle (with respect to some choice of metric on $\mathcal{E}$).  Then the map from the punctured bundle to the disk bundle $\mathcal{E}_0\hookrightarrow\mathcal{D}$ defined by $(x,v)\mapsto \left(x,\frac{v}{|v|}\right)$ is a retraction, hence induces an isomorphism on cohomology rings $H^*(\mathcal{E}_0)\cong H^*(\mathcal{D})$.  But note that the orientation reversing bundle isomorphism $\alpha\colon U^*\cong \mathcal{E}_0\rightarrow\mathcal{E}_0\cong U^*$ defined by $\alpha(x,v)=\left(x,\frac{v}{|v|^2}\right)$ restricts to the identity map on $\mathcal{D}\subset\mathcal{E}_0$.  Thus, the induced map on cohomology $\alpha^*\colon H^*(U^*)\rightarrow H^*(U^*)$ must be the identity map, which implies that the composition $\psi=\alpha\circ\phi\colon U_2\setminus N_2\rightarrow U_1\setminus N_1$ induces the same map as $\phi\colon U_2\setminus N_2\rightarrow U_1\setminus N_1$ on cohomology.
\vskip 0.2cm
Therefore the partial almost-commuting grid of short exact sequences of complexes yields a partial commuting grid of cohomology groups \eqref{eq:LongGrid}

\begin{equation}
\label{eq:LongGrid}
\xymatrix{\vdots & \vdots & \vdots \\
	H^{2q+1}(\mathcal{E},\mathcal{E}_0)\ar[r]\ar[u] & H^{2q+1}(M_1,M_1^*)\oplus H^{2q+1}(M_2,M_2^*)\ar[r]\ar[u] & H^{2q+1}(U_1,U_1^*)\ar[u]\\
	H^{2q}(M_1\#_NM_2)\ar[r]\ar[u] & H^{2q}(M_1^*)\oplus H^{2q}(M_2^*)\ar[r]^-{\iota_1^*-\iota_2^*\circ\psi^*}\ar[u] & H^{2q}(U_1^*)\ar[u]& \\
	H^{2q}(M_1\times_NM_2)\ar[r]\ar@{-->}[u]^-{\theta} & H^{2q}(M_1)\oplus H^{2q}(M_2)\ar[r]^-{\iota_1^*-\iota_2^*\circ\phi^*}\ar[u] & H^{2q}(U_1)\ar[u] \\
	H^{2q}(\mathcal{E},\mathcal{E}_0)\ar[r]\ar[u]^-{\rho} & H^{2q}(M_1,M_1^*)\oplus H^{2q}(M_2,M_2^*)\ar[r]^-{\iota_1^*-\iota_2^*\circ\phi^*}\ar[u] & H^{2q}(U_1,U_1^*)\ar[u]\ar[lluuu]\\
	H^{2q-1}(M_1\#_NM_2)\ar[r]\ar@{-->}[u]^-{\theta} & H^{2q-1}(M_1^*)\oplus H^{2q-1}(M_2^*)\ar[r]\ar[u] & H^{2q-1}(U_1^*)\ar[u]\ar[lluuu]\\
	\vdots\ar[u] & \vdots\ar[u] & \vdots\ar[u]\\}
\end{equation}

First note that by our assumptions, $H^{2q-1}(M_i)=0$, and hence $H^{2q-1}(N)\cong H^{2q-1}(U)=0$ also vanishes.  Thus, by the Thom Isomorphism Lemma \ref{lem:Thom}, we have $H^{2q-1}(\mathcal{E},\mathcal{E}_0)\cong H^{2q-1}(M_i,M_i^*)\cong H^{2q-1-2(d-k)}(N)=0$ also vanishes.  Moreover it follows from surjectivity of the restriction maps $\iota_i^*\colon H^*(M_i)\rightarrow H^*(N)$ that the Gysin maps 
\vskip 0.2cm
{\centering{$H^*(M_i,M_i^*)\cong H^{*-2(d-k)}(N)\rightarrow H^*(M_i)$}\vskip 0.2cm\noindent}
are injective, cf. Remark \ref{rem:PDGysin}.  It follows by the commutativity of the diagram that the map $\rho\colon H^{2q}(\mathcal{E},\mathcal{E}_0)\rightarrow H^{2q}(M_1\times_NM_2)$ is also injective.  

By the above analysis, Grid \eqref{eq:LongGrid} breaks into the following grid of exact sequences of cohomology groups:
\
\begin{equation}
\label{eq:LongGrid2}
\xymatrix{& 0 & 0 & 0 & \\
	&H^{2q}(M_1\#_NM_2)\ar[r]^-{\epsilon}\ar[u] & H^{2q}(M_1^*)\oplus H^{2q}(M_2^*)\ar[r]^-{\iota_1^*-\iota_2^*\circ\psi^*}\ar[u] & H^{2q}(U_1^*)\ar[u]& \\
	0\ar[r] & H^{2q}(M_1\times_NM_2)\ar[r]_-{\gamma}\ar@{-->}[u]^-{\theta} & H^{2q}(M_1)\oplus H^{2q}(M_2)\ar[r]^-{\iota_1^*-\iota_2^*\circ\phi^*}\ar[u]_-{\alpha} & H^{2q}(U_1)\ar[u]\ar[r] & 0 \\
	0\ar[r] & H^{2q}(\mathcal{E},\mathcal{E}_0)\ar[r]\ar[u]^-{\rho} & H^{2q}(M_1,M_1^*)\oplus H^{2q}(M_2,M_2^*)\ar[r]^-{\iota_1^*-\iota_2^*\circ\phi^*}\ar[u] & H^{2q}(U_1,U_1^*)\ar[u]\ar[r] & 0\\
	& 0\ar[u] & 0\ar[r]\ar[u] & H^{2q-1}(U_1^*)\ar[u]\ar[lluuu]\\
	&  & & 0\ar[u] & \\}
\end{equation}
Set 
$$A_*=H^{2*}(M_1), \ B_*=H^{2*}(M_2), \ T_*=H^{2*}(N)$$
with surjective projection maps $\pi_A=\iota_1^*\colon A\rightarrow T$ and $\pi_B=\iota_2^*\colon B\rightarrow T$.
Then note that the exact sequence in the middle row of Grid \eqref{eq:LongGrid2} is the analogue of Sequence \eqref{eq:SESFP}, and hence we have a ring isomorphism 
$$\left(A\times_TB\right)_*\cong H^{2*}(M_1\times_NM_2),$$
which is item (2).

Next note that by composing with the Thom isomorphism, the left vertical map $\rho\colon H^{2q}(\mathcal{E},\mathcal{E}_0)\rightarrow H^{2q}(M_1\times_NM_2)$ is really just multiplication by the Thom class 
$$\times\tau\colon H^{2q-2(d-k)}(N)\rightarrow H^{2q}(M_1\times_NM_2)$$
and the restriction of the Thom class to $H^{2(d-k)}(M_i)$ is exactly the Thom class of the normal bundle $\mathcal{E}_i\rightarrow N$, which proves item (1).

Finally we need to justify the existence of a surjective map $\theta\colon H^{2q}(M_1\times_NM_2)\rightarrow H^{2q}(M_1\#_NM_2)$ making the left column a short exact sequence.  
By the commutativity and the exactness of Grid \eqref{eq:LongGrid2}, we claim there exists a (non-canonical) map of $\Q$-vector spaces $\theta\colon H^{2q}(M_1\times_NM_2)\dashrightarrow H^{2q}(M_1\#_NM_2)$ which makes Grid \eqref{eq:LongGrid2} commute.  Indeed, from Grid~\eqref{eq:LongGrid2} we extract the diagram  
$$\xymatrix{H^{2q}(M_1\#_NM_2)\ar[r]^-\epsilon & \operatorname{Im}(\epsilon)\ar[r] & 0\\ H^{2q}(M_1\times_N M_2)\ar@{-->}[u]^-{\theta}\ar[ur]_-{\alpha\circ \gamma} & &}$$ 
and we deduce that $\theta$ exists because $H^{2q}(M_1\times_NM_2)$ is a projective $\Q$-module.

Thus, the leftmost column gives an exact sequence of even cohomology groups analogous to Sequence \eqref{eq:SESCS}:
\begin{equation}
\label{eq:CSSESM}
\xymatrixcolsep{3pc}\xymatrix{0\ar[r] & H^{2*-2(d-k)}(N)\ar[r]^-{\times\tau_1\oplus\times\tau_2} & H^{2*}(M_1\times_NM_2)\ar[r]^-\theta & H^{2*}(M_1\#_NM_2)\ar[r] & 0}
\end{equation} 
and hence we see that as graded vector spaces we have the equality
$$\left(A\#_TB\right)_*\cong H^{2*}(M_1\#_NM_2)$$
which is item (3), and this completes the proof.

\begin{remark}
	\label{rem:CSNP}
	Looking back, we see that Proposition \ref{ex:NonSLP} reflects the fact the connected sum of two projective manifolds need not be a projective manifold.  Indeed consider the topological connected sum of a complex projective space with itself over a projective subspace, say $X=\C\P^{m-1}\#_{\C\P^{t-1}}\C\P^{m-1}$.  Then the cohomology rings are $H^{2*}(\C\P^{m-1},\F)\cong\F[x]/(x^m)=A,B$ and $H(\C\P^{t-1},\F)\cong \F[z]/(z^t)=T$, and Theorem \ref{thm:TopCS} implies that the cohomology ring of the topological connected sum is the algebraic connected sum of cohomology rings:  
	$$H^{2*}(X,\F)\cong A\#_TB\cong \frac{\F[z_1,z_2]}{(z_1^{m-t},z_1^tz_2-z_2^2)}.$$
	On the other hand, we saw in Proposition \ref{ex:NonSLP} that this ring does not have SLP, and hence $X$ cannot be a projective manifold, or even homotopically equivalent to one. 
\end{remark}\noindent
{\bf Question}. There have been studies of the algebraic rational homotopy properties of the connected sum of two manifolds over a point, see for example \cite[\S 3.1.2]{FOT}.
It would be of interest to see what results might extend to our more general setting.\par

\begin{ack}
We thank J. Migliore and U. Nagel for organizing the special session on Lefschetz properties at the March, 2018 AMS Sectional meeting in Columbus, Ohio, where this project was started. We thank J. Watanabe for a correction, Z. Teitler for an example, and P. Macias Marques for comments. We appreciate comments and suggestions of the referee. The second author appreciates his sabbatical support from Endicott College in Spring 2019, and the hospitality of Northeastern University for a visit January - March, 2019.  The third author was supported by NSF grant DMS--1601024 and EpSCOR award OIA--1557417.
\end{ack}

\end{document}